\documentclass[11pt]{article}
\usepackage[margin=1in]{geometry}
\usepackage{hyperref}
\hypersetup{%
    colorlinks=true, 
    urlcolor=blue, 
    citecolor=blue, 
    linkcolor=blue%
}
\usepackage{amsmath, amsthm, amsfonts, mathtools, xfrac}
\usepackage{tikz, pgfplots}
\usepackage[dvipsnames]{xcolor} %For ForestGreen, etc.
\usepackage{tikz-cd}
\usepackage{csquotes} %For lengthy blockquotes
\usetikzlibrary{decorations.markings, calc, intersections, calligraphy}
\usepackage{tcolorbox}
\tcbuselibrary{skins}
\usepackage[numbers, sort&compress]{natbib}
\usepackage{doi} %For goog doi citations
\usepackage{epigraph}
\usepackage{ragged2e} %For \justifying

\makeatletter

\newcommand{\doi@}[1]{\urlstyle{same}\url{https://doi.org/#1}}
\DeclareRobustCommand{\doi}{\hyper@normalise\doi@}
\makeatother

\newcommand{\flushleftright}{\noindent\justifying}

\newcommand{\Equals}{\xrightleftharpoons{\quad}}

\title{A Categorical Approach to Euclidean Ratios and Proportions}
\author{
  Matthew Petersen \\
  \small Department of Mathematics, Whitman College
}
\begin{document}

\maketitle

\begin{abstract}
A categorial approach to the non-metric geometry in Books V and VI of Euclid's \textit{Elements} is presented. Specifically, we introduce a diagrammatic syntax that can be overlaid immediately on his diagrams, thus bridging intuitive presentation with fidelity to Euclid's arguments. This syntax makes complicated definitions like V.5, and indeed the arguments throughout books V and VI, including arguments about similar figures, intuitively clear. We show in an appendix that this syntax can be used to solve a puzzle regarding ancient mathematics. Finally, we offer evidence that this approach to Euclidean diagrams is rooted in the Aristotelian tradition itself, and that a similar syntax was utilized, in antiquity, for related questions of numeric and proportions. Thus the syntax is plausibly faithful to Euclid's own thought-world, and not an outside-imposition.
\end{abstract}

\epigraph{%
\small ``One way of thinking \textit{unlike} Euclid is to use the algebraic approach to interpreting his works.'' \raisebox{0.5ex}{\texttildelow}Sabetai \citet{unguru1979}.%
}

Books V and VI are two of the most difficult of Euclid's books; indeed, to understand them we are often forced to adopt an algebraic symbolism foreign to Euclid.\footnote{This algebra is foreign, not because he did not know algebra, but because pre-modern algebra was a numeric problem-solving technique, not a universal means of treating number and magnitude; see Sialaros and Christianidis  \citep{SialarosChristianidis2016}.} For example, in explaining definition V.5, Heath \citep{Heath1908ElementsVol2}, p. 125, writes ``Now let \( \sfrac{x}{y} \) and \( \sfrac{x'}{y'} \) be equal ratios in Euclid's sense.'' And Fitzpatrick \citep{Fitzpatrick2008Euclid}, p. 131, glosses definition V.1 ``In other words, \( \alpha\) is said to be a part of \(\beta\) if \(\beta = m\alpha\).'' As the epigraph indicates, historians of mathematics have long argued this approach has several problems. For example, this notation presupposes a unit, but Euclid \textit{never} employs a unit, and it obscures Euclid's type system---for Euclid, a ratio is neither a line and nor a number, but the algebraic notation reduces all three to a single variable that can be written \(x\).\footnote{It is true that \( \mathbb{N} \) is not \( \mathbb{R} \), and the use of \(m\) indicates an element of \(\mathbb{N} \). But we usually treat \(\mathbb{N} \) as a subset of \(\mathbb{R}\), rather than as a separate category that has an isomorphic copy canonically embedded in \(\mathbb{R}\), so that \(m\) is not of a different \textit{type} from \(x\). Furthermore, even if mathematicians parse these type-theoretic implications in a Euclidean manner, the general reader (or the student) is \textit{extremely} unlikely to. Fitzpatrick does note that in his glosses, he uses Greek letters to refer to \textit{magnitudes}, but even if his algebra here works immediately with unmeasured magnitudes (as Oaks \citep{Oaks2022} has shown Vieta's and Descartes' algebra did), this is a subtle point, and it remains difficult for modern mathematicians to read such an algebra without re-introducing real-numbers.} But it is difficult to see what other alternatives we have for explaining these challenging texts.

On the other hand, historians of mathematics have know since at least the early 1980's that, in contrast to modern ratios, Euclid's ratios are not objects, but relations.\footnote{Thus, Ian Mueller \protect\citep{mueller1981} ``In general a ratio is a relation between two objects and not itself an object or even a pair of objects, as in the modern foundational definition of fractions as pairs of integers''(\textit{Philosophy of Mathematics and Deductive Structure in
Euclid's \emph{Elements,}} p. 66). And Fabio Acerbi \protect\citep{acerbi2021} ``No object is identified by a binary relation---to repeat what should be a commonplace; a Greek ratio is not a mathematical object'' (\textit{The Logical Syntax of Greek Mathematics,} p. 171). Aristotle \cite{aristotle-meta} is also clear that ratios are relations---to something---in \textit{Metaphysics} V.15.} This suggests a \textit{categorical} treatment  in which ratios are written as arrows. Nevertheless---and quite surprisingly---to the best of our knowledge, Euclid's ratios have never been treated categorically as arrows between objects.\footnote{Zeeman \citep{zeeman1977} argues Euclid should be interpreted Categorically, but his approach differs from ours because he did not treat the ratio itself as an arrow. Perhaps more importantly, though we too want to understand modern math and ancient math together, we here use Category Theory as a tool for highlighting and expositing Euclid's original sense, and distinguishing it from our habitual mathematics, not for reformulating his ideas in modern terms. It may, however, be interesting to consider Zeeman's suggestions of a category intermediary between magnitudes and ratios, in the light of Aristotle's note in \textit{Metaphysics} X.1 that in music the diesis served as unit of measurement. In Pythagorean music theory, intervals are ratios---for example the diesis is the ratio of 531441 to 524288. Once these intervals (that is, ratios) are quantified via measurement with the diesis or other intervals, then we should be able to consider, in a manner at least in the spirit of ancient mathematics, ratios of ratios. For example, the interval we call a major ninth is the double of a perfect fifth, and a semi-tone is (approximately) half a tone.} This in spite the fact that the common Greek term for relation, \textit{pros ti,} toward something, used by Nicomachus to describe ratios, in the standard Euclidean description of a ratio, and still fossilized in English (``as is $A$ \textit{to} $B$ so is $C$ \textit{to} $D$''), points to such a categorical treatment.\footnote{Indeed, though he does not develop them as such, Euclid's ratios nearly form an \textit{allegory}, in the sense of Freyd and Scedrov \protect\citep{freyd-scedrov1990}, \textit{Categories, Allegories}. See Appendix \ref{AppendixAllegories}.}

This suspicion that a categorical treatment is apt is strengthened by the following quote from Thomas de Vio Cajetanus. Standing at the end of a long commentary tradition, in his commentary \citep{cajetan1587} on Aristotle's \textit{Categories} \citep{aristotle-cat}, Cajetan gives the justification for what he takes to be a Platonic neologism for relation, (\textit{pros ti}), toward something---the very language preserved in our language for ratios. Cajetan writes:\footnote{Page 171. Author's translation.}

\begin{quotation}
It should be seen why \textit{to something} is said more than relation or relative. It should therefore be known that (as Porphyry, Alexander and Albert expressly say) Plato invented \textit{to something}, for, in fact, it was called relative by the ancients, as he himself says. The invention of Plato, however, is shown reasonable and necessary by a double reason.

First, from the part of the mode of signifying: Because ``relation" [\textit{relatio}, literally, traced-back-to-ness] signifies relation as conceived, and not as exercised; on account of which, to something is not expressed according to that name, just as neither fathership nor sonship. But ``relative" [\textit{relativus}, literally, traced-back-to-some] signifies relation as exercising more than as exercised. For it imports a certain potency for referring [\textit{referendum}, literally, tracing back to]. But \textit{to something} signifies the very relation as exercised, for it signifies a respect [\textit{respectum}, literally, looking again at] as terminated to another: which is (as thus it may be permitted to speak) the exercise of the relation, since then the relation is posited in truth and its proper exercise, when respecting toward something is posited [\textit{cum ad aliquid respicere ponitur}]. It is agreed, however, that a thing as exercised is signified better for knowing than otherwise: and because of this Plato fittingly invented this name.

Secondly, because of the thing signified. For the very nature of respect or relation is so fittingly and expressly signified through ``to something" that I think it is impossible to be signified more clearly by a unique name. For since relation is between two extremes, namely, fundament and terminus, relation does not have what it is from the act of being that it possesses in the fundament, but from this that it respects the terminus. Hence here (6a37) in the [Platonic] definition of relatives it is said that “this very what they are, [is said to be] of others,” that is, they are of termini. But I do not know how better the respecting of termini is signified than through ``to something", for ``to" shows the respect itself, and ``something" the terminus. Therefore the proper essence of relation which, as will be shown below, more fully consists in the act of being toward another, is best signified for knowing through \textit{to something}. And much better than the name of relation which does not explain the aforesaid, but signifies it as a certain thing; or than the name of relative, which imports it more from the part of the fundament.

\end{quotation}

How better can this ``act of being toward something'' be signified that with an arrow? For the arrow itself signifies the respect or relation, whereas the \textit{terminus} of the relation is signified by the object it points to and the fundament by the object it points from.\footnote{The fundament and terminus are, of course, in categorical terms, the domain and codomain. I make no attempt to change our terminology, but, while the commonalities should indeed be noted, the terminological difference should alert us to real philosophical differences.}

This understanding of ratio was also shared by Isaac Barrow \citep{barrow-ML}, Newton's predecessor as the Lucastrian Chair of Mathematics at the University of Cambridge, whose mathematical philosophy seems to have been adopted by Newton himself. Writing against Wallis and Hobbes and, as he says, nearly all moderns, Barrow argues at length that a ratio is not a quantity but a relation whose very essence is to be to another. Thus it seems plausible that even that for Newton ratios are \textit{to something}, and should be treated categorically.

This paper therefore develops a Euclidean relational calculus that is, we hope, true to Euclid's thought-world, can be overlaid directly on his diagrams, and that offers a non-algebraic means of making notoriously difficult passages of Euclid clear.

\section{The Relational Calculus}

Rather than writing $AB:CD$ as is traditional, we write, in categorical fashion, \(AB\xrightarrow{to} CD\). This notation brings out that ratios are relations, \textit{to something,} rather than objects, and draws considerably more attention to the ratio itself. For now, we write ``\textit{to}'' over the arrow to help legibility, though it is not the name of this particular arrow, and therefore it will be dropped.

But the main advantage of this approach is that the diagram itself can be annotated to show not only the quantities, but their ratios. Thus, when Euclid's diagram includes lines $A$ and $B,$ and we need to consider the ratio of $A$ to $B$, we can add an arrow to the diagram itself, and write:

\centering
\begin{tikzcd}
{\begin{tikzpicture}
\draw (0,0) to
		node[midway, above, inner sep = 2] {$A$}
	  (1.5,0);
\end{tikzpicture}}
\ar[r,
	shift right,
	"to"]
&
	{\begin{tikzpicture}
	\draw (0,0) to
			node[midway, above, inner sep = 2] {$B$}
		  (2.038,0);
	\end{tikzpicture}}
\end{tikzcd}

\flushleftright
From ancient perspective, equality is a ratio,\footnote{Note that this is \textit{not} equality of measures. Euclid's geometry is wholly non-metrical. When we say, for example, that the old national standard meters were all equal to the universal standard in Paris, we certainly do not mean that their measures are equal---since it is this equality that establishes their use as measures. This is the sense that Euclid uses equality: a equality of the extensions that would be the foundation for any metric equality.} as are inequality, greater than and less than, as Nicomachus \citep{nicomachus-EN} says (Book I Chapter XVII),\footnote{``Of quantity \textit{to something}, the highest generic divisions are two, equality and inequality.\ldots The unequal, on the other hand, is split up by subdivisions, and one part of it is the greater, the other the less'' (translation modified on the basis of \citep{nicomachus-GK}). He then goes on, in Chapter XIX to list ratios like ``$3: 2,$ $6: 4,$ $9: 6,$ $12 : 8$''  as examples of inequality (translation modified on the basis of \protect\cite{nicomachus-GK}). This tradition is also found in Thomas Aquinas \protect\citep{aquinas-st} who writes ``It should be said that ratio (\textit{proportio}) is said doubly. In one mode, a certain habitude of one quantity to another, according to which double, triple and equal are species of ratio'' (\textit{ST }I.12a1, my translation). We could perhaps say that greater than and less than are types of ratio, not actually existing ratios. But all ratios (and quantities) considered in Greek mathematics are types. Even the most specific ratios like sesquialter (the ratio of three to two) are more like leaf types than concrete particulars. The actual concrete instances of those types are relations in and toward physical things in the world. Moreover, though the thing in the world is, in this case, more actual than the type, this is not because it is concrete, but because it is a substance, and the relational form in it only an accident. But the relational form which is known, universally, in mathematics, and exists as a particular in the individual. For this distinction, see Avicenna's \protect\citep{avicenna-meta} \textit{Metaphysics,} V.1. Though Acerbi \protect\citep{acerbi2021} does not attempt to understand the philosophy, in \textit{The Logical Syntax of Greek Mathematics,} he argues, I think correctly, that Euclidean mathematics is always universal, and never particular: there is, for example, no universal instantiation in Greek mathematics. Though he works from a specifically Piercian perspective, Claas Lattmann \protect\citep{lattmann2018} is also worth reading here. It seems to me that many of his Piercian concerns can also be stated in an ancient mode. For example, they can be tied to Aristotle's \protect\citep{aristotle-soul} claim (\textit{On the Soul} III.8) that one who contemplates simultaneously contemplates a particular. In the diagrams Euclid provides particulars that exemplify the matter at hand, and that the student can form in their imagination simultaneously with their contemplation of the universal signified by the words.} and thus \textit{to something.} Indeed, Needham \citep{needham2021} notes that Newton continually writes that two quantities ``ultimately have the ratio of equality.''\footnote{Page xx.} This expression only makes sense if, for Newton, equality is a \textit{ratio.}

Because equality is a special ratio, I suggest, in contexts that work in a Euclidean context, we write it with a pair of harpoons, \(\Equals\), which is reminiscent enough of our familiar symbol =,\footnote{We overload = with a call to the `execute' function, but it was not used, even in pre-modern algebra, quite like we use it today; see Oaks \protect\citep{oaks2009}. Specifically, from at least Diophantus through to Vieta, algebraists rigorously distinguished between what we might describe as calls to the evaluate function on a series of operations, and a statement that two quantities have the ratio of equality.} while suggesting a categorical \textit{pointing.} Thus, for example, in \textit{Data} I \citep{euclid-data}, after producing four magnitudes, Euclid says ``equal is $A$ to $C$ and $B$ to $D$'' (my translation). This can be drawn as follows:\footnote{Though note that Euclid does not say $A$ and $B$ are equal to each other, but that $A$ is equal to $B$. That is, though the reciprocal relations are the same, and so can, without harm, both be noted, as I do below, he only explicitly notes the arrow from $A$ to $C$ and from $B$ to $D$.}

\centering%
\begin{tikzcd}
{\begin{tikzpicture}
\draw (0,0) -- 
		node[midway, above] {$A$}
	  (2.71828,0);
\end{tikzpicture}}
\ar[d,
	rightharpoonup,
	shift left = 0.5]
&
	{\begin{tikzpicture}
	\draw (0,0) --
			node[midway, above] {$B$}
		  (3.14159,0);
	\end{tikzpicture}}
	\ar[d,
		rightharpoonup,
		shift left = 0.5]\\
{\begin{tikzpicture}
\draw (0,0) -- 
		node[midway, above] {$C$}
	  (2.71828,0);
\end{tikzpicture}}
\ar[u,
	rightharpoonup,
	shift left = 0.5]
&
	{\begin{tikzpicture}
	\draw (0,0) --
			node[midway, above] {$D$}
		  (3.14159,0);
	\end{tikzpicture}}
	\ar[u,
		rightharpoonup,
		shift left = 0.5]
\end{tikzcd}

\flushleftright%
\subsection{Proportion or Analogy}
Though this relational syntax makes a beginning, Book V and VI of the \textit{Elements} are not built on ratios alone, but on what can be called a proportion or analogy.\footnote{I prefer ``analogy,'' the transliteration of Euclid's word, to ``proportion'' as the latter has strong connotations of direct proportionality in algebra. Though the two concepts are interrelated, they are distinctly different. Specifically, when we say two quantities (e.g., force and acceleration) are directly proportional we compare the two quantities to each other, for example, forces to accelerations. But here, we compare homogeneous quantities, and the link between different genera of quantity is made through the comparison of ratios of homogeneous quantities. In the medieval scholastic terminology, this is an analogy of proper proportionality.} Euclid's treatment of analogies begins in Book V, definition 5, where he defines \textit{same ratio.} For Euclid ratios are never first-class objects independent of the quantities that have the ratio, so, here, he says ``magnitudes are said to be in the same ratio\ldots'' This phrasing could make it sound like sameness of ratio is a quaternary relation of magnitudes. But throughout Book V, including immediately in definition 6, Euclid says that a first magnitude \textit{has} the same ratio to a second as a third \textit{has} to a fourth. Therefore, it seems that ``same ratio'' is a relation of relations.\footnote{Though here, the precise nature of this ratio of ratios is debated. Barrow argues that they are not actually ratios of ratios, but that we treat them as such linguistically. In this he follows the Thomistic position. On the other hand, his opponents follow Scotus. Either way, it is linguistically treated as a ratio of ratios.} This will therefore be written as follows:

\centering
\begin{tikzcd}[row sep = 1.5cm, math mode = false]
{\begin{tikzpicture}
\draw (0,0) to
		node[midway, above, inner sep = 2] {$A$}
	  (1.5,0);
\end{tikzpicture}}
\ar[r,
	shift right,
	"$to$" {font = \tiny},
	""' {name = AB, pos = 0.49}]
&
	{\begin{tikzpicture}
	\draw (0,0) to
			node[midway, above, inner sep = 2] {$B$}
		  (2.038,0);
	\end{tikzpicture}}\\
%SEcond Row
{\begin{tikzpicture}
\draw (0,0) to
		node[midway, above, inner sep = 2] {$C$}
	  (1.75,0);
\end{tikzpicture}}
\ar[r,
	shift right,
	"$to$" {font = \tiny, name = CD, pos = 0.51}]
&
	{\begin{tikzpicture}
	\draw (0,0) to
			node[midway, above, inner sep = 2] {$D$}
		  (2.378,0);
	\end{tikzpicture}}.
\ar[from = AB,
	to = CD,
	Leftrightarrow,
	gray,
	"same" {sloped, font = \small},
	"ratio"' {sloped, font = \small}]
\end{tikzcd}

\flushleftright%
Euclid also says that one quantity can have a greater ratio to a second than a third has to a fourth. This relation of ratios is illustrated similarly, though obviously, the arrow does not point both directions.\footnote{Euclid doesn't treat these as functors that relate whole categories, but as relations of relations. Also, though, as noted below occasionally it is helpful to label a ratios rather than writing ``same ratio''; in books V and VI, Euclid does not name the ratios but always names them in terms of the magnitudes. For this reason, and because it better matches ``greater ratio than'' it is better to write this as a relation of relations: this practice allows more clear annotations on the diagrams, and allows the notation to be more consistent with ``greater ratio than''.}

In the diagrams, this arrow is used to connect ratios that cannot be labeled with a species of ratio. But Euclid's ratios seem to be belong to a Julia-like \citep{juliatypes2026} type-system in which abstract types of ratio can be reasoned about, but cannot be instantiated except in and through a concrete leaf-type ratio. The exact status of these relations of relations was debated in the middle ages.\footnote{See John of St. Thomas  \protect\citep{johnofstthomas1985}, Second Part of the Logic q17a3.} But they seem to refer to the concrete leaf-type, in and through which an abstract type is instantiated, even when the leaf-type of the ratios cannot be known. So for example, below we will see this arrow connecting two arrows labeled ``multiple of'', as in the following diagram:

\centering
\begin{tikzcd}[font = \scriptsize, row sep = 1.5cm, column sep = 2cm]
{\begin{tikzpicture}
\draw (0,0) to node[midway, above, inner sep = 2] {E} (3*0.5,0);
\foreach \n in {1,2} \draw (\n*0.5,-0.05) -- ++(0,0.1);
\end{tikzpicture}}
\ar[r,
	black!60,
	shift right,
	"\text{multiple of}",
	""' {name = EA}]
&
	{\begin{tikzpicture}
	\draw (0,0) to node[midway, above, inner sep = 2] {A} (0.5,0);
	\end{tikzpicture}}\\
	%Second Row
{\begin{tikzpicture}
\draw (0,0) to node[midway, above, inner sep = 2] {F} (3*0.4,0);
\foreach \n in {1,2} \draw (\n*0.4,-0.05) -- ++(0,0.1);
\end{tikzpicture}}
\ar[r,
	black!60,
	shift right,
	"\text{multiple of}"',
	""{name = FD}]
&
	{\begin{tikzpicture}
	\draw (0,0) to  node[midway, above, inner sep = 2] {C} (0.4,0);
	\end{tikzpicture}},
\ar[from = EA,
	to = FD,
	Leftrightarrow,
	gray!75,
	"\text{same}" {sloped, font = \small},
	"\text{ratio}"' {sloped, font = \small}]
\end{tikzcd}

\flushleftright%
The claim is not the trivial claim that ``multiple of'' is the same abstract ratio type as ``multiple of'' (that is, that both ratios have the same genus, multiple of). Rather it makes a claim about the the concrete, leaf ratio type (that is, the most specific species of ratio) that would actually be instantiated. \textit{These} are, in fact the same.\footnote{For this reason, it is occasionally helpful to label the ratio of two quantities with the ratio of different quantities that paradigmatically have that ratio. For example, since 3 and 2 are relatively prime, we can write $6\xrightarrow{\;3\text{-to-}2\;}4$. This expresses the same fact as the display-style notation with an arrow joining the ratios.}

Though for simplicity I have drawn all the quantities in the diagrams above as lines, the quantities on the two rows are not necessarily of the same genus. Therefore, since Euclid does not allow ratios between quantities of different genera, there may be no first-order relation between the quantities on different rows.\footnote{Barrow also fiercely defends the claim that there cannot be non-homogeneous ratios. Euclid's position, here, is similar to our sense that a claim like 2 feet equals (or is the double of or half of) 2 Joules is poorly typed and makes no sense. Later Greek mathematicians, like Heron, would consider operations from one genus to another, but, in our terminology, these belong to an entirely different category than Euclid's ratios do. Euclid's rule that ratios must be between homogeneous quantities functions something like a type-checking: If, after a series of calculations, you find yourself saying that 2J = 2 ft, or that 6J is the double of 2C, your statement is dimensionally inconsistent and so doesn't type-check correctly.} However, in definition 6, Euclid says that magnitudes that have the same ratio are called analogous. Therefore, we can add the following quasi-relation to the diagram:

\centering
\begin{tikzcd}[row sep = 2cm, math mode = false]
{\begin{tikzpicture}
\draw (0,0) to
		node[midway, above, inner sep = 2] {$A$}
	  (1.5,0);
\end{tikzpicture}}
\ar[r,
	shift right,
	"$to$" {font = \tiny},
	""' {name = AB, pos = 0.49}]
\ar[d,
	Leftrightarrow,
	gray,
	densely dotted,
	"analogous" {sloped, font = \small},
	"quantities"' {sloped, font = \small}]
&
	{\begin{tikzpicture}
	\draw (0,0) to
			node[midway, above, inner sep = 2] {$B$}
		  (2.038,0);
	\end{tikzpicture}}\\
%SEcond Row
{\begin{tikzpicture}
\draw (0,0) to
		node[midway, above, inner sep = 2] {$C$}
	  (1.75,0);
\end{tikzpicture}}
\ar[r,
	shift right,
	"$to$" {font = \tiny, name = CD, pos = 0.51}]
&
	{\begin{tikzpicture}
	\draw (0,0) to
			node[midway, above, inner sep = 2] {$D$}
		  (2.378,0);
	\end{tikzpicture}}.
\ar[from = AB,
	to = CD,
	Leftrightarrow,
	gray,
	"same" {sloped, font = \small},
	"ratio"' {sloped, font = \small}]
\end{tikzcd}

\flushleftright%
In the \textit{Categories} \citep{aristotle-cat}, which Euclid seems to presuppose, Aristotle argues that reciprocal relations are simultaneous,\footnote{In the terms of Freyd and Scedrov, we may say that for every ratio \(A\xrightarrow{to}B\) there is a reciprocal relation \((A\xrightarrow{to}B)^{\circ} = B\xrightarrow{to}A.\) } and in the corollary following Proposition 7, Euclid shows that if \(A\) is analogous to \(C\), then \(B\) is also analogous to \(D\). We can therefore, without much harm, add that analogy into the diagram.

\centering
\begin{tikzcd}[row sep = 2cm, math mode = false]
{\begin{tikzpicture}
\draw (0,0) to
		node[midway, above, inner sep = 2] {$A$}
	  (1.5,0);
\end{tikzpicture}}
\ar[r,
	shift right,
	"$to$" {font = \tiny},
	""' {name = AB, pos = 0.49}]
\ar[d,
	Leftrightarrow,
	gray,
	densely dotted,
	"analogous" {sloped, font = \small},
	"quantities"' {sloped, font = \small}]
&
	{\begin{tikzpicture}
	\draw (0,0) to
			node[midway, above, inner sep = 2] {$B$}
		  (2.038,0);
	\end{tikzpicture}}
	\ar[d,
		Leftrightarrow,
		gray,
		densely dotted,
		"analogous" {sloped, font = \small},
		"quantities"' {sloped, font = \small}]\\
%SEcond Row
{\begin{tikzpicture}
\draw (0,0) to
		node[midway, above, inner sep = 2] {$C$}
	  (1.75,0);
\end{tikzpicture}}
\ar[r,
	shift right,
	"$to$" {font = \tiny, name = CD, pos = 0.51}]
&
	{\begin{tikzpicture}
	\draw (0,0) to
			node[midway, above, inner sep = 2] {$D$}
		  (2.378,0);
	\end{tikzpicture}}.
\ar[from = AB,
	to = CD,
	Leftrightarrow,
	gray,
	"same" {sloped, font = \small},
	"ratio"' {sloped, font = \small}]
\end{tikzcd}

\flushleftright%

This diagram looks like a 2-cell in a double category---and it is these squares that Euclid reasons about---but readers should note the real differences between a double category and what Euclid is doing here. Specifically, the solid gray arrow founds the light gray arrows, and so the light arrows are not, to use our language, the morphisms in Euclid's category. For this reason, same-ratio acts more like a simple 2-morphism in a 2-category.

Indeed, it is interesting to note that Elements V.11 and 13 show a sort of vertical composition, whereas Elements V.22 shows a sort of horizontal composition.

It may be possible to attempt to push the Categorical understanding of Euclid further. For example, if we restrict our attention to the most specific species of ratio (as opposed to genera), then we could say that Euclid works in a family of isolated codiscrete categories. But these squares provide a fundamental unity to his mathematical world. 

But we should be careful not to place Euclid on the Procrustean bed of modern category theory, or to see him as merely anticipating our modern notions. Rather, it seems better to use the categorical diagrams as a tool for illuminating his arguments, and the significant differences between his mathematics and ours.

Euclid's relational logic, and the corresponding notation developed here, is powerful enough that we are able to use this notation to diagram every step in every argument in Books V and VI of the \textit{Elements,} making the arguments visually legible, to modern mathematics students, in a manner that is mathematically and philosophically faithful to Euclid himself. I have already prepared such diagrams both books, in their entirety, and for parts of Book VII. Here, I will merely present some highlights demonstrating the viability of this notation.\footnote{I believe that similar diagrams can be used to illustrate some points in ancient and medieval philosophy. For example, when Aristotle \citep{aristotle-physics} says that as ``as the bronze is to the statue, [and] the wood to the bed,'' so is matter to form (\textit{Physics} I.7), he is using ``is to'' in the same manner as Euclid is. It also can illuminate claims in the heights of philosophy and theology in authors from ibn al-Arabi and al-Ghazali to Thomas Aquinas and Domingo Ba\~nez. This, however, is not the place for these considerations.} But first, we must head-off a potential misunderstanding.

\section{Heading off a Misunderstanding}

In modern mathematics it is common (though not necessary) to think of morphisms as \textit{actions} or \textit{transformations.} Thus, for example, Reihl \citep{riehl2017} (page xi) writes ``A \textit{category} is a context for the study of a particular class of mathematical objects. Importantly, a category\ldots has both ``nouns'' and ``verbs,'' containing specified collections of objects and transformations, called \textit{morphisms,} between them.'' Ancient mathematics from Heron to Diophantus explicitly consider mathematical transformations or actions, and these can aptly receive a categorical treatment. However, these are completely separate, even stylistically, from Euclid's ratios. Euclid's ratios are beings that exist in and toward existing quantities. We can aptly say that, in contrast to the dynamism of the actions described in Heron and Diophantus, Euclid's ratios are static, provided we strip the connotations of value in the terms ``dynamic'' and ``static''. For a Greek, a dynamic coming to be is less ``real'' and actual than a static being. We can perhaps capture something of that if we think of the difference between mathematicians standing perfectly silent and still as they stare at a board contemplating a new proof, and undergraduates who approach a novel proof ``dynamically'', that is, without standing so still, and without the long contemplative silences.

Some examples of dynamic, functional, Greek mathematics can perhaps help make this distinction clear between kinds of morphism clear.

In Book I Proposition 2 of Heron's \textit{Metrica}, specifically in the Method for finding the area and hypotenuse of a right-triangle given its legs (here exemplified with a 3-4-5 right triangle), Heron writes: ``Multiplying the 3 by the 4, take their half. It makes 6. Of so much is the area of the triangle. And its hypotenuse: Multiplying the 3 by itself and similarly the 4 by itself, add. And they make 25. And, taking a side of these, have the hypotenuse of the triangle.''\footnote{My translation of the Greek text found in Acerbi and Vitrac \protect\citep{acerbi-vitrac2014}.} Heron's text can be annotated relatively easily with mapsto arrows:

\begin{tikzcd}[math mode = false, font = \small,
			  column sep = 5cm,
			  row sep = 1.5ex,
			  /tikz/column 1/.append style={anchor=base east},
			  /tikz/column 2/.append style={anchor=base east},
			  /tikz/column 3/.append style={anchor=base west}]
\textit{Multiplying} the 3 by the 4 \( \left(\mapsto{} 12\right)\),\\[-2.5ex]
\textbf{take} their half.
	\ar[rr, maps to,"\tiny It makes"]
	&[-5cm]	& 6;\\[-2.5ex]
	&	& \hspace{-1.25em}$=$ \,the area of\\[-2.5ex]
	&	& \hspace{-1.25em}\phantom{$=$} \,the triangle.\\[2ex]
\textit{Multiplying} the 3 by itself \( \left(\mapsto{} 9\right)\)\\[-2.5ex]
	& \textbf{add.}
	  \ar[r, maps to,"\tiny And they make"]
		& 25;\\[-2.5ex]
and, the 4 by itself \( \left(\mapsto{} 16\right),\)\\[2ex]
\textit{taking} a side, \textbf{have}
	\ar[rr, maps to]
	&	& {\color{gray}5,}\\[-2.5ex]
	&	& \hspace{-1.25em}$=$ \,the hypotenuse\\[-2.5ex]
	&	& \hspace{-1.25em}\phantom{$=$} \,of the triangle.
\ar[from = 5-1,
	to = 7-1,
	-,
	start anchor = north east,
	end anchor = south east,
	decorate,
	decoration = {calligraphic brace}]
\end{tikzcd}

As in modern category theory, the arrows represent \textit{actions} or transformations of measures of one mathematical object to obtain the measures of a different object.

While Heron presents this function in what Acerbi \citep{acerbi2021} calls an algorithmic language, and so it is aptly written with mapsto arrows, the corresponding procedure\footnote{The language of algorithm and procedure is taken from Acerbi \citep{acerbi2021}, pages 12--23. Following Acerbi's usage, an algorithm operates with exemplary particulars, like ``the 3'' and ``the 4'' here. A procedure describes the operations and the kinds of objects operated on. These seem to roughly correspond to our distinction between maps-to notation, e.g., \(\left(\text{the }3, \text{the 4}\right)\xmapsto{\times} 12 \) and functional notation, e.g., \(\text{Base}\times\text{Height}\xrightarrow{\times} \text{Area}\).} can be written with functional arrows as a pair of commutative diagrams. However, the complexity of the procedures in Greek mathematics make wiring diagrams\footnote{As in Fong and Spivak \citep{fong-spivak2019}.} much more helpful, generally.

Though I have prepared wiring diagrams for Heron's subsequent Propositions measuring triangles, the power of this approach becomes clear when applied to following two procedures from Diophantus' \textit{On Polygonal Numbers,} the first of which can be used to calculate a polygonal number given its side and the number of vertices, and the second, its inverse, can be used to calculate the side of a polygonal number:

\begin{quotation}
In fact, taking the side of the polygonal always doubling we shall subtract a unit, and multiplying the remainder by the [number] less by a dyad than the multiplicity of the angles we shall always add a dyad to the result, and taking the square on the result we shall subtract from it the [square] on the [number] less by a tetrad than the multiplicity of the angles, and dividing the remainder by the octuple of the [number] less by a dyad than the multiplicity of the angles, we shall find the sought polygonal.

And again, the polygonal itself being given, we shall find the side as follows: multiplying it by the octuple of the [number] less by a dyad than the multiplicity of the angles and adding to the result the square on the [number] less by a tetrad than the multiplicity of the angles we shall find a square---whenever the assigned [number] be really polygonal---and always subtracting from the side of this square a dyad we shall divide the remainder by the [number] less by a dyad than the multiplicity of the angles, and adding to the result a unit and taking half of the result we shall have the side of the sought polygonal.\footnote{Quoted in Acerbi \protect\citep{acerbi2021}, p. 13.}
\end{quotation}

This dense text becomes clear when annotated with wiring diagrams. The second is the inverse of the first, so I have maintained the geometry, and it should be read right-to-left:

\centering
\includegraphics[scale=1]{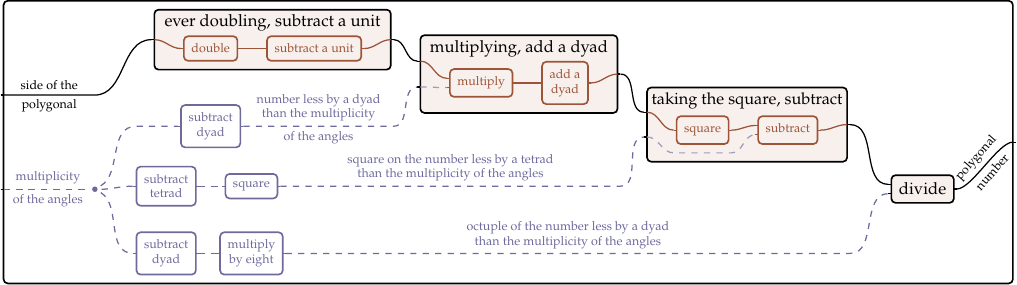}\\
\includegraphics[scale=1]{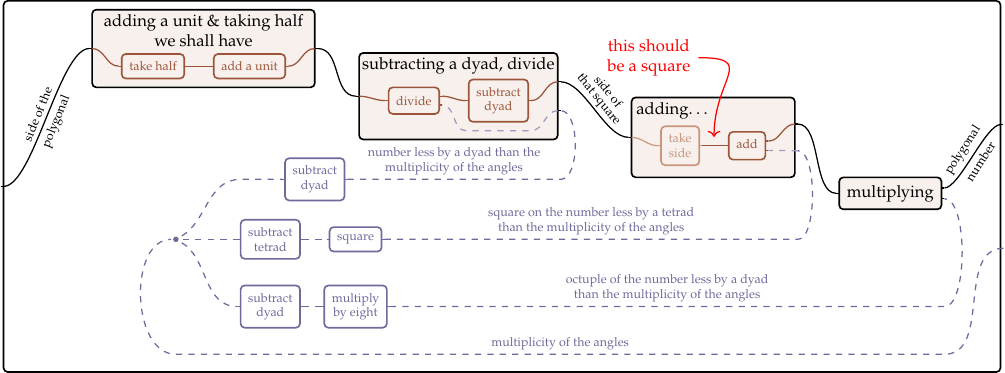}

\flushleftright
The red text and arrow ``this should be a square'' reflects Diophantus' ``We shall find a square---whenever the assigned [number] be really polygonal'' and refusal to actually take its side. And of course, the rhythms of Diophantus' text in the second version could be better captured with a diagram that isn't drawn to show that the second procedure is the inverse of the first.

However, Euclid's ratios do not produce a new line or number, they relate already existing ones---and this is even marked stylistically.\footnote{See Acerbi \citep{acerbi2011} A similar distinction between operations that produce something new and static beings can be found in Greek and Arabic algebra, as Oaks has shown repeatedly \citep{Oaks2022, Oaks2018, Oaks2017, Oaks2010, oaks2009, ChristianidisOaks2024}. }

We need to keep these two categories strictly separate in our minds. Indeed, if we restrict our attention to operations like taking the half or doubling quantities, Euclid's category is like the categorical \textit{dual} of this category of operations. That is, whereas we double three to make six, or halve six to make three; for Euclid, six is the double of three and three is the half of six---and these arrows point from the double to the half, or from the half to the double:\footnote{In the diagram on the right I have not written ``maps to'' arrows because, for the ancients, \textit{a} three or \textit{a} six, at least as considered in Euclid, is a multitude of six things that are internally coherent and bodily distinct from other things. In the science of mathematics, however, these particular multitudes are known universally, according to the very nature common to every three or six. Thus this or that three is, roughly, a term of type three. See, for example Avicenna \protect\citep{avicenna-meta} \textit{The Metaphysics of the Healing,} I.3 and III.2--6. For further discussion of this issue, Klein \protect\citep{klein1968} and Paulus Soncinas \protect\citep{Soncinas1586}, Book X, \textit{Quaestio} 3. In his Elizabethean translation and commentary on Euclid's \textit{Elements,} Billingsley \protect\citep{billingsley1570}, the first person to translate Euclid into English, also clearly holds that Euclid's \textit{monas} is a principle of unity and distinction \textit{in the world.} (See his comments on Book VII definition 1.) That there is no unique 3 for Greek mathematicians is clearly stated by Mueller \citep{mueller1981}, p.59 ``In Greek arithmetic there are indefinitely many units and indefinitely many ways of combining them into multitudes. Clearly then, there is no unique 2 or 3; any pair of units is a 2, for example.''}

\centering%
\begin{tikzcd}[column sep = 2cm]
3
\ar[r,
	maps to,
	"\text{doubles to}",
	bend left]
&
	6
	\ar[l,
		maps to,
		"\text{halves to}",
		bend left]
\end{tikzcd}\qquad
\begin{tikzcd}[column sep = 2cm]
3
\ar[r,
	"\text{half of}"',
	bend right]
&
	6
	\ar[l,
		"\text{double of}"',
		bend right]
\end{tikzcd}

\flushleftright
We should therefore be careful not to consider the ratios as \textit{actions} or scalings of lines. In \textit{The Joy of Abstraction,}\footnote{p. 53.} Eugenia Cheng \citep{cheng2022} introduces commutative diagrams concretely by drawing diagrams of familial relations. For example, 

\centering%
\begin{tikzcd}[column sep = 1.5cm]
A
\ar[r,
	"\text{sister of}"]
\ar[rrr,
	bend right,
	"\text{great aunt of}"']
&
	B
	\ar[r,
		"\text{mother of}"]
	&
		C
		\ar[r,
			"\text{mother of}"]
		&
			D
\end{tikzcd} 

\flushleftright%
These relations are certainly \textit{not} transformations of the persons, but inhere in the persons by virtue of natality. Similarly, Euclid's ratios and analogies are not transformations---though the ratios are opposite to multiplicative transformations---but inhere in the existing quantities themselves.

\section{Diagramming Euclid}
How, then, can we apply these diagrammatic principles, in detail, to Euclid? Let's begin with Euclid's definition of ``Same same ratio.''

\begin{quote}
Magnitudes are said to be \textit{in the same ratio,} the first to the second and the third to the fourth, when, if any equimultiples whatever are taken of the first and third, and any equimultiples whatever of the second and fourth, the former equimultiples alike exceed, are alike equal to, or alike fall short of, the latter equimultiples respectively taken in corresponding order.
\end{quote}

(Translation by David E. Joyce \citep{joyce-elements}.)

We can diagram the property defined---being in the same ratio---as above:

\centering%
\begin{tikzcd}[font = \scriptsize, row sep = 1.5cm]
{\begin{tikzpicture}
\draw (0,0) to node[midway, above, inner sep = 2] {A} (0.5,0);
\end{tikzpicture}}
\ar[r,
	black!60,
	shift right,
	"\text{ratio to}",
	""' {name = AB}]
&
{\begin{tikzpicture}
\draw (0,0) to  node[midway, above, inner sep = 2] {B} (0.75,0);
\end{tikzpicture}}\\
%Second Row
{\begin{tikzpicture}
\draw (0,0) to  node[midway, above, inner sep = 2] {C} (0.4,0);
\end{tikzpicture}}
\ar[r,
	black!60,
	shift right,
	"\text{ratio to}"',
	"" {name = CD}]
&
{\begin{tikzpicture}
\draw (0,0) to  node[midway, above, inner sep = 2] {D} (0.6,0);
\end{tikzpicture}}
\ar[from = AB,
	to = CD,
	Leftrightarrow,
	Red!50,
	"\text{same}" sloped,
	"\text{ratio}"' sloped]
\end{tikzcd}

\flushleftright%
I have written ``same ratio'' in red to highlight that this is the property we are attempting to define, and have drawn the arrows slightly lighter than the lines so the eye isn't overwhelmed by lines.

Multiple is a species of ratio---as Nicomachus says \citep{nicomachus-EN}---and ``equimultiple'' is just ``same ratio'' or ``analogous'' when the ratio in question is a multiple. So the situation, after equi-multiples have been taken, can be diagrammed as follows:

\centering%
\begin{tikzcd}[font = \scriptsize, row sep = 2cm, column sep = 2cm]
{\begin{tikzpicture}
\draw (0,0) to node[midway, above, inner sep = 2] {E} (3*0.5,0);
\foreach \n in {1,2} \draw (\n*0.5,-0.05) -- ++(0,0.1);
\end{tikzpicture}}
\ar[r,
	black!60,
	shift right,
	"\text{multiple of}",
	""' {name = EA}]
&
	{\begin{tikzpicture}
	\draw (0,0) to node[midway, above, inner sep = 2] {A} (0.5,0);
	\end{tikzpicture}}
	\ar[r,
		black!60,
		"\text{ratio to}",
		shift right,
		""' {name = AB}]
	&
	{\begin{tikzpicture}
	\draw (0,0) to  node[midway, above, inner sep = 2] {B} (0.75,0);
	\end{tikzpicture}}
	&
		{\begin{tikzpicture}
		\draw (0,0) to  node[midway, above, inner sep = 2] {G} (4*0.75,0);
		\foreach \n in {1,2,3}
			\draw (\n*0.75,-0.05) -- ++(0,0.1);
		\end{tikzpicture}}
		\ar[l,
			black!60,
			shift left,
			"\rotatebox{180}{\tiny multiple of}" {name = GC}]\\
	%Second Row
{\begin{tikzpicture}
\draw (0,0) to node[midway, above, inner sep = 2] {F} (3*0.4,0);
\foreach \n in {1,2} \draw (\n*0.4,-0.05) -- ++(0,0.1);
\end{tikzpicture}}
\ar[r,
	black!60,
	shift right,
	"\text{multiple of}"',
	""{name = FD}]
&
	{\begin{tikzpicture}
	\draw (0,0) to  node[midway, above, inner sep = 2] {C} (0.4,0);
	\end{tikzpicture}}
	\ar[r,
		black!60,
		"\text{ratio to}"',
		shift right,
		"" {name = CD}]
	&
	{\begin{tikzpicture}
	\draw (0,0) to  node[midway, above, inner sep = 2] {D} (0.6,0);
	\end{tikzpicture}}
		&
		{\begin{tikzpicture}
		\draw (0,0) to  node[midway, above, inner sep = 2] {H} (4*0.6,0);
		\foreach \n in {1,2,3}
			\draw (\n*0.6,-0.05) -- ++(0,0.1);
		\end{tikzpicture}}
		\ar[l,
			black!60,
			shift left,
			"\rotatebox{180}{\tiny multiple of}"' {name = HB}]
\ar[from = AB,
	to = CD,
	Leftrightarrow,
	Red!50,
	"\text{same}" sloped,
	"\text{ratio}"' sloped]
\ar[from = EA,
	to = FD,
	Leftrightarrow,
	gray!75,
	"\text{same}" sloped,
	"\text{multiple}"' sloped]
\ar[from = GC,
	to = HB,
	Leftrightarrow,
	gray!75,
	"\text{same}" sloped,
	"\text{multiple}"' sloped]
\end{tikzcd}

\flushleftright%
Exceeding and falling short are genera of ratio---as Nicomachus \citep{nicomachus-EN} says---and being equal to is itself a ratio. So the condition ``alike exceed, are alike equal to, or alike fall short of'' can also be diagrammed with arrows:

\centering%
\begin{tikzcd}[font = \scriptsize, row sep = 2cm, column sep = 2cm]
{\begin{tikzpicture}
\draw (0,0) to node[midway, above, inner sep = 2] {E} (3*0.5,0);
\foreach \n in {1,2} \draw (\n*0.5,-0.05) -- ++(0,0.1);
\end{tikzpicture}}
\ar[r,
	black!60,
	shift right,
	"\text{multiple of}",
	""' {name = EA}]
\ar[rrr,
	blue,
	bend left = 1.3cm,
	"\text{greater than}"]
\ar[rrr,
	ForestGreen,
	bend left = 0.9cm,
	"\text{equal to}"]
\ar[rrr,
	Orange,
	bend left = 0.5cm,
	"\text{less than}"]
&
	{\begin{tikzpicture}
	\draw (0,0) to node[midway, above, inner sep = 2] {A} (0.5,0);
	\end{tikzpicture}}
	\ar[r,
		black!60,
		"\text{ratio to}",
		shift right,
		""' {name = AB}]
	&
	{\begin{tikzpicture}
	\draw (0,0) to  node[midway, above, inner sep = 2] {B} (0.75,0);
	\end{tikzpicture}}
	&
		{\begin{tikzpicture}
		\draw (0,0) to  node[midway, above, inner sep = 2] {G} (4*0.75,0);
		\foreach \n in {1,2,3}
			\draw (\n*0.75,-0.05) -- ++(0,0.1);
		\end{tikzpicture}}
		\ar[l,
			black!60,
			shift left,
			"\rotatebox{180}{\tiny multiple of}" {name = GC}]\\
	%Second Row
{\begin{tikzpicture}
\draw (0,0) to node[midway, above, inner sep = 2] {F} (3*0.4,0);
\foreach \n in {1,2} \draw (\n*0.4,-0.05) -- ++(0,0.1);
\end{tikzpicture}}
\ar[r,
	black!60,
	shift right,
	"\text{multiple of}"',
	""{name = FD}]
\ar[rrr,
	blue,
	shift right = 2,
	bend right = 0.5cm,
	"\text{greater than}"]
\ar[rrr,
	ForestGreen,
	shift right =2,
	bend right = 0.9cm,
	"\text{equal to}"]
\ar[rrr,
	Orange,
	shift right =2,
	bend right = 1.3cm,
	"\text{less than}"]
&
	{\begin{tikzpicture}
	\draw (0,0) to  node[midway, above, inner sep = 2] {C} (0.4,0);
	\end{tikzpicture}}
	\ar[r,
		black!60,
		"\text{ratio to}"',
		shift right,
		"" {name = CD}]
	&
	{\begin{tikzpicture}
	\draw (0,0) to  node[midway, above, inner sep = 2] {D} (0.6,0);
	\end{tikzpicture}}
		&
		{\begin{tikzpicture}
		\draw (0,0) to  node[midway, above, inner sep = 2] {H} (4*0.6,0);
		\foreach \n in {1,2,3}
			\draw (\n*0.6,-0.05) -- ++(0,0.1);
		\end{tikzpicture}}
		\ar[l,
			black!60,
			shift left,
			"\rotatebox{180}{\tiny multiple of}"' {name = HB}]
\ar[from = AB,
	to = CD,
	Leftrightarrow,
	Red!50,
	"\text{same}" sloped,
	"\text{ratio}"' sloped]
\ar[from = EA,
	to = FD,
	Leftrightarrow,
	gray!75,
	"\text{same}" sloped,
	"\text{multiple}"' sloped]
\ar[from = GC,
	to = HB,
	Leftrightarrow,
	gray!75,
	"\text{same}" sloped,
	"\text{multiple}"' sloped]
\end{tikzcd}

\flushleftright%
The relations of greater than, equal to, and less than are certainly not simultaneously pointing from $E$ to $G$ or $F$ to $D.$ Rather, they illustrate the three types of ratio that could exist from $E$ to $G$, and the corresponding ratios that must be shown to exist from $F$ to $H$ if $A$ is in the same ratio to $B$ as $C$ is to $D$. Here, of course, we can see, by visual inspection, that $E$ is less than $G$ and $F$ than $H$. But we need to show that the ratios on colored ratios on the top row entail the colored ratios on the bottom row. (Or vice-versa, though Euclid doesn't argue in that direction.)

If we want all the arrows to point in the same direction, we could take the reciprocal of the ``multiple of'' ratio from $G$ to $B$ and from $H$ to $D$, though Euclid does not generally do this, and it breaks the symmetry between the two sides.\footnote{This looks like the composition of multiple of, ratio of and part of. But these species of ratios compose to give the highest genus, namely, \textit{ratio}---which is therefore uninformative. Rather, Eudoxius has broken this relation into its three species, greater than, equal to, and less than; and requires that these species coordinate on the two levels.}

\centering%
\begin{tikzcd}[font = \scriptsize, row sep = 2cm, column sep = 2cm]
{\begin{tikzpicture}
\draw (0,0) to node[midway, above, inner sep = 2] {E} (3*0.5,0);
\foreach \n in {1,2} \draw (\n*0.5,-0.05) -- ++(0,0.1);
\end{tikzpicture}}
\ar[r,
	black!60,
	shift right,
	"\text{multiple of}",
	""' {name = EA}]
\ar[rrr,
	blue,
	bend left = 1.3cm,
	"\text{greater than}"]
\ar[rrr,
	ForestGreen,
	bend left = 0.9cm,
	"\text{equal to}"]
\ar[rrr,
	Orange,
	bend left = 0.5cm,
	"\text{less than}"]
&
	{\begin{tikzpicture}
	\draw (0,0) to node[midway, above, inner sep = 2] {A} (0.5,0);
	\end{tikzpicture}}
	\ar[r,
		black!60,
		"\text{ratio to}",
		shift right,
		""' {name = AB}]
	&
	{\begin{tikzpicture}
	\draw (0,0) to  node[midway, above, inner sep = 2] {B} (0.75,0);
	\end{tikzpicture}}
	&
		{\begin{tikzpicture}
		\draw (0,0) to  node[midway, above, inner sep = 2] {G} (4*0.75,0);
		\foreach \n in {1,2,3}
			\draw (\n*0.75,-0.05) -- ++(0,0.1);
		\end{tikzpicture}}
		\ar[l,
			leftarrow,
			black!60,
			shift left,
			"\text{part of}"',
			"" {name = GC}]\\
	%Second Row
{\begin{tikzpicture}
\draw (0,0) to node[midway, above, inner sep = 2] {F} (3*0.4,0);
\foreach \n in {1,2} \draw (\n*0.4,-0.05) -- ++(0,0.1);
\end{tikzpicture}}
\ar[r,
	black!60,
	shift right,
	"\text{multiple of}"',
	""{name = FD}]
\ar[rrr,
	blue,
	shift right = 2,
	bend right = 0.5cm,
	"\text{greater than}"]
\ar[rrr,
	ForestGreen,
	shift right =2,
	bend right = 0.9cm,
	"\text{equal to}"]
\ar[rrr,
	Orange,
	shift right =2,
	bend right = 1.3cm,
	"\text{less than}"]
&
	{\begin{tikzpicture}
	\draw (0,0) to  node[midway, above, inner sep = 2] {C} (0.4,0);
	\end{tikzpicture}}
	\ar[r,
		black!60,
		"\text{ratio to}"',
		shift right,
		"" {name = CD}]
	&
	{\begin{tikzpicture}
	\draw (0,0) to  node[midway, above, inner sep = 2] {D} (0.6,0);
	\end{tikzpicture}}
		&
		{\begin{tikzpicture}
		\draw (0,0) to  node[midway, above, inner sep = 2] {H} (4*0.6,0);
		\foreach \n in {1,2,3}
			\draw (\n*0.6,-0.05) -- ++(0,0.1);
		\end{tikzpicture}}
		\ar[l,
			leftarrow,
			black!60,
			shift left,
			"\text{part of}",
			""' {name = HB}]
\ar[from = AB,
	to = CD,
	Leftrightarrow,
	Red!50,
	"\text{same}" sloped,
	"\text{ratio}"' sloped]
\ar[from = EA,
	to = FD,
	Leftrightarrow,
	gray!75,
	"\text{same}" sloped,
	"\text{multiple}"' sloped]
\ar[from = GC,
	to = HB,
	Leftrightarrow,
	gray!75,
	"\text{same}" sloped,
	"\text{part}"' sloped]
\end{tikzcd}

\flushleftright%

\section{Vertical Composition}
In Proposition 11, Euclid demonstrates something like vertical composition of 2-cells. That is ``Ratios which are the same with the same ratio are also the same with one another.'' (David E. Joyce, translator \citep{joyce-elements}.) 

Diagramatically this can be stated as follows:

\centering%
\begin{tcolorbox}[%
	colbacktitle=white,
	coltitle = black,
	center title,
	halign upper=center,
	halign lower=center,
	fonttitle = \bfseries,
	colbacklower = white,
	colback = white,
	width = \linewidth,
	lower separated=true,
	sidebyside,
	boxrule=.2mm]
\begin{tikzcd}[ampersand replacement = \&, row sep = 8ex, column sep = 3em, math mode = false]
{\begin{tikzpicture}
\draw[|-|] (0,0) -- (0.5,0);
\node[inner sep = 2, anchor = east] (Letter) at (0,0) {\tiny $A$};
\node[inner sep = 2, anchor = east] at (Letter.west) {\tiny As};
\end{tikzpicture}}
\ar[r,
	shift right,
	"\tiny to",
	""' {name = AB}]
\ar[d,
	Leftrightarrow,
	gray!50,
	densely dotted,
	shift left = 3,
	"\rotatebox{90}{\tiny Analogous}"',
	"\rotatebox{90}{\tiny Quantities}"]
\&
	{\begin{tikzpicture}
	\draw[|-|] (0,0) -- (0.3,0);
	\node[inner sep = 2, anchor = east] at (0,0) {\tiny $B$};
	\node[inner sep = 2, anchor = west] at (0.3,0) {\tiny ,};
	\end{tikzpicture}}
	\ar[d,
		Leftrightarrow,
		gray!50,
		densely dotted,
		shift left,
		"\rotatebox{90}{\tiny Analogous}"',
		"\rotatebox{90}{\tiny Quantities}"]\\
%
%
% Second Row
{\begin{tikzpicture}
\draw[|-|] (0,0) -- (0.667,0);
\node[inner sep = 2, anchor = east] (Letter) at (0,0) {\tiny $C$};
\node[inner sep = 2, anchor = east] at (Letter.west) {\tiny so};
\end{tikzpicture}}
\ar[r,
	shift right,
	"\tiny to",
	"\phantom{\tiny to}" {name = CD}]
\&
	{\begin{tikzpicture}
	\draw[|-|] (0,0) -- (0.4,0);
	\node[inner sep = 2, anchor = east] at (0,0) {\tiny $D$};
	\node[inner sep = 2, anchor = west] at (0.4,0) {\tiny ;};
	\end{tikzpicture}}\\[-9ex]
%
%
% Third Row
{\begin{tikzpicture}
\draw[|-|] (0,0) -- (0.667,0);
\node[inner sep = 2, anchor = east] (Letter) at (0,0) {\tiny $C$};
\node[inner sep = 2, anchor = east] at (Letter.west) {\tiny as};
\end{tikzpicture}}
\ar[r,
	shift right,
	"\tiny to",
	""' {name = CD2, pos = 0.505}]
\ar[d,
	Leftrightarrow,
	gray!50,
	densely dotted,
	shift left = 3,
	"\rotatebox{90}{\tiny Analogous}"',
	"\rotatebox{90}{\tiny Quantities}"]
\&
	{\begin{tikzpicture}
	\draw[|-|] (0,0) -- (0.4,0);
	\node[inner sep = 2, anchor = east] at (0,0) {\tiny $D$};
	\node[inner sep = 2, anchor = west] at (0.4,0) {\tiny ,};
	\end{tikzpicture}}
	\ar[d,
		Leftrightarrow,
		gray!50,
		densely dotted,
		shift left,
		"\rotatebox{90}{\tiny Analogous}"',
		"\rotatebox{90}{\tiny Quantities}"]\\
%
%
% Fourth Row
{\begin{tikzpicture}
\draw[|-|] (0,0) -- (1,0);
\node[inner sep = 2, anchor = east] (Letter) at (0,0) {\tiny $E$};
\node[inner sep = 2, anchor = east] at (Letter.west) {\tiny so};
\end{tikzpicture}}
\ar[r,
	shift right,
	"\tiny to",
	"\phantom{\tiny to}" {name = EF, pos = 0.495}]
\&
	{\begin{tikzpicture}
	\draw[|-|] (0,0) -- (0.6,0);
	\node[inner sep = 2, anchor = east] at (0,0) {\tiny $F$};
	\node[inner sep = 2, anchor = west] at (0.6,0) {\tiny .};
	\end{tikzpicture}}
\ar[	from = AB,
	to = CD,
	Leftrightarrow,
	gray,
	"\rotatebox{90}{\scriptsize Same}"',
	"\rotatebox{90}{\scriptsize Ratio}"]
\ar[	from = CD2,
	to = EF,
	Leftrightarrow,
	gray,
	"\rotatebox{90}{\scriptsize Same}"',
	"\rotatebox{90}{\scriptsize Ratio}"]
\end{tikzcd}	

\tcblower

\flushleftright%
The claim is:

\centering%
\begin{tikzcd}[ampersand replacement = \&, row sep = 8ex, column sep = 3em, math mode = false]
{\begin{tikzpicture}
\draw[|-|] (0,0) -- (0.5,0);
\node[inner sep = 2, anchor = east] (Letter) at (0,0) {\tiny $A$};
\node[inner sep = 2, anchor = east] at (Letter.west) {\tiny As};
\end{tikzpicture}}
\ar[r,
	shift right,
	"\tiny to",
	""' {name = AB, pos = 0.51}]
\ar[d,
	Leftrightarrow,
	gray!50,
	densely dotted,
	shift left = 3,
	"\rotatebox{90}{\tiny Analogous}"',
	"\rotatebox{90}{\tiny Quantities}"]
\&
	{\begin{tikzpicture}
	\draw[|-|] (0,0) -- (0.3,0);
	\node[inner sep = 2, anchor = east] at (0,0) {\tiny $B$};
	\node[inner sep = 2, anchor = west] at (0.3,0) {\tiny ,};
	\end{tikzpicture}}
	\ar[d,
		Leftrightarrow,
		gray!50,
		densely dotted,
		shift left,
		"\rotatebox{90}{\tiny Analogous}"',
		"\rotatebox{90}{\tiny Quantities}"]\\
%
%
% Second Row
{\begin{tikzpicture}
\draw[|-|] (0,0) -- (1,0);
\node[inner sep = 2, anchor = east] (Letter) at (0,0) {\tiny $E$};
\node[inner sep = 2, anchor = east] at (Letter.west) {\tiny so};
\end{tikzpicture}}
\ar[r,
	shift right,
	"\tiny to",
	"\phantom{\tiny to}" {name = EF, pos = 0.49}]
\&
	{\begin{tikzpicture}
	\draw[|-|] (0,0) -- (0.6,0);
	\node[inner sep = 2, anchor = east] at (0,0) {\tiny $F$};
	\node[inner sep = 2, anchor = west] at (0.6,0) {\tiny .};
	\end{tikzpicture}}
\ar[	from = AB,
	to = EF,
	Leftrightarrow,
	gray,
	"\rotatebox{90}{\scriptsize Same}"',
	"\rotatebox{90}{\scriptsize Ratio}"]
\end{tikzcd}	
\end{tcolorbox}

\flushleftright%
It is clear, diagrammatically, that this is something like vertical composition.

Euclid argues this by first taking equimultiples of $A,$ $C$ and $E$ and other equimultiples of $B,$ $D$ and $F$:

\centering%
\begin{tikzcd}[ampersand replacement = \&, row sep = 8ex, column sep = 2em, math mode = false]
{\begin{tikzpicture}
\draw[|-|] (0,0) -- (1,0);
\node[inner sep = 2, anchor = east] (Letter) at (0,0) {\tiny $G$};
\end{tikzpicture}}
\ar[r,
	end anchor = {[xshift = 1ex]west},
	gray,
	shift right,
	""' {name = GA}]
\ar[d,
	Leftrightarrow,
	densely dotted,
	shift left = 1.5,
	gray,
	"\rotatebox{90}{\tiny Equimultiples}"',
	"\rotatebox{90}{\tiny of \scalebox{0.95}{$C;$ $A$}}"]
\ar[dd,
	Leftrightarrow,
	densely dotted,
	shift right,
	start anchor = {[xshift = -2ex]south},
	end anchor = {[xshift = -2ex]north},
	bend right = 1.75cm,
	gray,
	"\rotatebox{90}{\tiny Equimultiples}"',
	"\rotatebox{90}{\tiny of \scalebox{0.95}{$E;$ $A$}}"]
\&
	{\begin{tikzpicture}
	\draw[|-|] (0,0) -- (0.5,0);
	\node[inner sep = 2, anchor = east] (Letter) at (0,0) {\tiny $A$};
	\end{tikzpicture}}
	\ar[r,
		end anchor = {[xshift = 1ex]west},
		gray,
		shift right,
		""' {name = AB}]
	\ar[d,
		Leftrightarrow,
		gray!50,
		densely dotted,
		shift left,
		"\rotatebox{90}{\tiny Analogous}"',
		"\rotatebox{90}{\tiny Quantities}"]
	\&
		{\begin{tikzpicture}
		\draw[|-|] (0,0) -- (0.3,0);
		\node[inner sep = 2, anchor = east] at (0,0) {\tiny $B$};
		\end{tikzpicture}}
		\ar[r,
			leftarrow,
			gray,
			end anchor = {[xshift = 1ex]west},
			shift right,
			"" {name = LB}]
		\ar[d,
			Leftrightarrow,
			gray!50,
			densely dotted,
			shift left,
			"\rotatebox{90}{\tiny Analogous}"',
			"\rotatebox{90}{\tiny Quantities}"]
		\&
			{\begin{tikzpicture}
			\draw[|-|] (0,0) -- (0.9,0);
			\node[inner sep = 2, anchor = east] at (0,0) {\tiny $L$};
			\end{tikzpicture}}
			\ar[d,
				Leftrightarrow,
				densely dotted,
				shift left = 1.0,
				gray,
				"\rotatebox{-90}{\tiny Equimultiples}",
				"\rotatebox{-90}{\tiny of \scalebox{0.95}{$B;$ $D$}}"']
			\ar[dd,
				Leftrightarrow,
				densely dotted,
				shift left,
				start anchor = {[xshift = 2.25ex]south},
				end anchor = {[xshift = 2.25ex]north},
				bend left = 1.75cm,
				gray,
				"\rotatebox{-90}{\tiny Equimultiples}",
				"\rotatebox{-90}{\tiny of \scalebox{0.95}{$B;$ $F$}}"']\\
%
%
% Second Row
{\begin{tikzpicture}
\draw[|-|] (0,0) -- (1.333,0);
\node[inner sep = 2, anchor = east] (Letter) at (0,0) {\tiny $H$};
\end{tikzpicture}}
\ar[r,
	end anchor = {[xshift = 1ex]west},
	gray,
	shift right,
	"" {name = HC}]
\ar[d,
	Leftrightarrow,
	densely dotted,
	gray,
	shift left = 1.5,
	"\rotatebox{90}{\tiny Equimultiples}"',
	"\rotatebox{90}{\tiny of \scalebox{0.95}{$E;$ $C$}}"]
\&
	{\begin{tikzpicture}
	\draw[|-|] (0,0) -- (0.667,0);
	\node[inner sep = 2, anchor = east] (Letter) at (0,0) {\tiny $C$};
	\end{tikzpicture}}
	\ar[r,
		gray,
		shift right,
%		bend left,
		start anchor = {[yshift = 0.1ex]east},
		end anchor = {[xshift = 1ex, yshift = 0.5]west},
		"" {name = CD},
		""' {name = CD2}]
	\ar[d,
		Leftrightarrow,
		gray!50,
		densely dotted,
		shift left,
		"\rotatebox{90}{\tiny Analogous}"',
		"\rotatebox{90}{\tiny Quantities}"]
	\&
		{\begin{tikzpicture}
		\draw[|-|] (0,0) -- (0.4,0);
		\node[inner sep = 2, anchor = east] at (0,0) {\tiny $D$};
		\end{tikzpicture}}
		\ar[r,
			leftarrow,
			gray,
			end anchor = {[xshift = 1ex]west},
			shift right,
			"" {name = MD}]
		\ar[d,
			Leftrightarrow,
			gray!50,
			densely dotted,
			shift left,
			"\rotatebox{90}{\tiny Analogous}"',
			"\rotatebox{90}{\tiny Quantities}"]
		\&
			{\begin{tikzpicture}
			\draw[|-|] (0,0) -- (1.2,0);
			\node[inner sep = 2, anchor = east] at (0,0) {\tiny $M$};
			\end{tikzpicture}}
			\ar[d,
				Leftrightarrow,
				densely dotted,
				shift left = 1.0,
				gray,
				"\rotatebox{-90}{\tiny Equimultiples}",
				"\rotatebox{-90}{\tiny of \scalebox{0.95}{$D;$ $F$}}"']\\
%
%
% Third Row
{\begin{tikzpicture}
\draw[|-|] (0,0) -- (2,0);
\node[inner sep = 2, anchor = east] (Letter) at (0,0) {\tiny $K$};
\end{tikzpicture}}
\ar[r,
	gray,
	end anchor = {[xshift = 1ex]west},
	shift right,
	"" {name = KE}]
\&
	{\begin{tikzpicture}
	\draw[|-|] (0,0) -- (1,0);
	\node[inner sep = 2, anchor = east] (Letter) at (0,0) {\tiny $E$};
	\end{tikzpicture}}
	\ar[r,
		gray,
		end anchor = {[xshift = 1ex]west},
		shift right,
		"" {name = EF}]
	\&
		{\begin{tikzpicture}
		\draw[|-|] (0,0) -- (0.6,0);
		\node[inner sep = 2, anchor = east] at (0,0) {\tiny $F$};
		\end{tikzpicture}}
		\ar[r,
			leftarrow,
			gray,
			end anchor = {[xshift = 1ex]west},
			shift right,
			"" {name = NF}]
		\&
			{\begin{tikzpicture}
			\draw[|-|] (0,0) -- (1.8,0);
			\node[inner sep = 2, anchor = east] at (0,0) {\tiny $N$};
			\end{tikzpicture}}
\ar[	from = AB,
	to = CD,
	Leftrightarrow,
	gray,
	"\rotatebox{90}{\scriptsize Same}"',
	"\rotatebox{90}{\scriptsize Ratio}"]
\ar[	from = CD2,
	to = EF,
	Leftrightarrow,
	gray,
	"\rotatebox{90}{\scriptsize Same}"',
	"\rotatebox{90}{\scriptsize Ratio}"]
\end{tikzcd}	

\flushleftright%
Because it is easier to cross-over when I label the lines than when I label the arrows, and to avoid diagrammatic clutter, I have only labeled the lines as equimultiples, not the ratios. Nevertheless, that the ratios are the same should be understood---and is guaranteed by construction.

Euclid then reasons that since $A$ is to $B$ as $C$ is to $D$, equimultiples of $G$ and $H$ simultaneously exceed, are equal to, or are less than, equimultiples of $B$ and $D$. (This is, again, the definition of ``same ratio.'')

\centering%
\begin{tikzcd}[ampersand replacement = \&, row sep = 8ex, column sep = 2em, math mode = false]
{\begin{tikzpicture}
\draw[|-|] (0,0) -- (1,0);
\node[inner sep = 2, anchor = east] (Letter) at (0,0) {\tiny $G$};
\node[inner sep = 2, anchor = east] at (Letter.west) {\footnotesize If};
\end{tikzpicture}}
\ar[r,
	end anchor = {[xshift = 1ex]west},
	gray!50,
	shift right,
	""' {name = GA}]
\ar[d,
	Leftrightarrow,
	densely dotted,
	shift left = 3,
	gray,
	"\rotatebox{90}{\tiny Equimultiples}"',
	"\rotatebox{90}{\tiny of \scalebox{0.95}{$C;$ $A$}}"]
\ar[dd,
	Leftrightarrow,
	densely dotted,
	shift right,
	start anchor = {[xshift = -3ex]south},
	end anchor = {[xshift = -3ex]north},
	out = -135,
	in = 135,
	distance = 1.75cm,
	gray!50,
	"\rotatebox{90}{\tiny Equimultiples}"',
	"\rotatebox{90}{\tiny of \scalebox{0.95}{$E;$ $A$}}"]
\&
	{\begin{tikzpicture}
	\draw[|-|, black!60] (0,0) -- (0.5,0);
	\node[inner sep = 2, anchor = east, black!60] (Letter) at (0,0) {\tiny $A$};
	\end{tikzpicture}}
	\ar[r,
		end anchor = {[xshift = 1ex]west},
		gray!50,
		shift right,
		""' {name = AB}]
	\&
		{\begin{tikzpicture}
		\draw[|-|, black!60] (0,0) -- (0.3,0);
		\node[inner sep = 2, anchor = east, black!60] at (0,0) {\tiny $B$};
		\end{tikzpicture}}
		\ar[r,
			leftarrow,
			gray!50,
			end anchor = {[xshift = 1ex]west},
			shift right,
			"" {name = LB}]
		\&
			{\begin{tikzpicture}
			\draw[|-|] (0,0) -- (0.9,0);
			\node[inner sep = 2, anchor = east] at (0,0) {\tiny $L$};
			\node[anchor = west, inner sep = 2] at (0.9,0) {\footnotesize ,};
			\end{tikzpicture}}
			\ar[d,
				Leftrightarrow,
				densely dotted,
				shift left = 0,
				gray,
				"\rotatebox{-90}{\tiny Equimultiples}",
				"\rotatebox{-90}{\tiny of \scalebox{0.95}{$B;$ $D$}}"']
			\ar[dd,
				Leftrightarrow,
				densely dotted,
				shift left,
				start anchor = {[xshift = 2ex]south},
				end anchor = {[xshift = 2ex]north},
				bend left = 1.75cm,
				gray!50,
				"\rotatebox{-90}{\tiny Equimultiples}",
				"\rotatebox{-90}{\tiny of \scalebox{0.95}{$B;$ $F$}}"']\\
%
%
% Second Row
{\begin{tikzpicture}
\draw[|-|] (0,0) -- (1.333,0);
\node[inner sep = 2, anchor = east] (Letter) at (0,0) {\tiny $H$};
\node[inner sep = 2, anchor = east] at (Letter.west) {\footnotesize then};
\end{tikzpicture}}
\ar[r,
	end anchor = {[xshift = 1ex]west},
	gray!50,
	shift right,
	"" {name = HC}]
\ar[d,
	Leftrightarrow,
	densely dotted,
	gray,
	shift left = 3,
	"\rotatebox{90}{\tiny Equimultiples}"',
	"\rotatebox{90}{\tiny of \scalebox{0.95}{$E;$ $C$}}"]
\&
	{\begin{tikzpicture}
	\draw[|-|, black!60] (0,0) -- (0.667,0);
	\node[inner sep = 2, anchor = east, black!60] (Letter) at (0,0) {\tiny $C$};
	\end{tikzpicture}}
	\ar[r,
		gray!50,
		shift right = 0.5,
		bend left,
		start anchor = {[yshift = 0.1ex]east},
		end anchor = {[xshift = 1ex, yshift = 0.5]west},
		"" {name = CD}]
	\ar[r,
		gray!50,
		shift right = 1.5,
		bend right,
		start anchor = {[yshift = -0.1ex]east},
		end anchor = {[xshift = 1ex, yshift = -0.5]west},
		""' {name = CD2}]
	\ar[d,
		Leftrightarrow,
		gray!50,
		densely dotted,
		shift left,
		"\rotatebox{90}{\tiny Analogous}"',
		"\rotatebox{90}{\tiny Quantities}"]
	\&
		{\begin{tikzpicture}
		\draw[|-|, black!60] (0,0) -- (0.4,0);
		\node[inner sep = 2, anchor = east, black!60] at (0,0) {\tiny $D$};
		\end{tikzpicture}}
		\ar[r,
			leftarrow,
			gray!50,
			end anchor = {[xshift = 1ex]west},
			shift right,
			"" {name = MD}]
		\ar[d,
			Leftrightarrow,
			gray!50,
			densely dotted,
			shift left,
			"\rotatebox{90}{\tiny Analogous}"',
			"\rotatebox{90}{\tiny Quantities}"]
		\&
			{\begin{tikzpicture}
			\draw[|-|] (0,0) -- (1.2,0);
			\node[inner sep = 2, anchor = east] at (0,0) {\tiny $M$};
			\node[anchor = west, inner sep = 2] at (1.2,0) {\footnotesize .};
			\end{tikzpicture}}
			\ar[d,
				Leftrightarrow,
				densely dotted,
				shift left = 0.5,
				gray,
				"\rotatebox{-90}{\tiny Equimultiples}",
				"\rotatebox{-90}{\tiny of \scalebox{0.95}{$D;$ $F$}}"']\\
%
%
% Third Row
{\begin{tikzpicture}
\draw[|-|] (0,0) -- (2,0);
\node[inner sep = 2, anchor = east] (Letter) at (0,0) {\tiny $K$};
\end{tikzpicture}}
\ar[r,
	gray,
	end anchor = {[xshift = 1ex]west},
	shift right,
	"" {name = KE}]
\&
	{\begin{tikzpicture}
	\draw[|-|] (0,0) -- (1,0);
	\node[inner sep = 2, anchor = east] (Letter) at (0,0) {\tiny $E$};
	\end{tikzpicture}}
	\ar[r,
		gray,
		end anchor = {[xshift = 1ex]west},
		shift right,
		"" {name = EF}]
	\&
		{\begin{tikzpicture}
		\draw[|-|] (0,0) -- (0.6,0);
		\node[inner sep = 2, anchor = east] at (0,0) {\tiny $F$};
		\end{tikzpicture}}
		\ar[r,
			leftarrow,
			gray,
			end anchor = {[xshift = 1ex]west},
			shift right,
			"" {name = NF}]
		\&
			{\begin{tikzpicture}
			\draw[|-|] (0,0) -- (1.8,0);
			\node[inner sep = 2, anchor = east] at (0,0) {\tiny $N$};
			\end{tikzpicture}}
\ar[	from = AB,
	to = CD,
	Leftrightarrow,
	gray!25,
	"\rotatebox{90}{\scriptsize Same}"',
	"\rotatebox{90}{\scriptsize Ratio}"]
\ar[	from = CD2,
	to = EF,
	Leftrightarrow,
	gray!75,
	"\rotatebox{90}{\scriptsize Same}"',
	"\rotatebox{90}{\scriptsize Ratio}"]
\ar[from = 1-1,
	to = 1-4,
	bend left = 0.55cm,
	blue,
	start anchor = north east,
	end anchor = {[xshift = 1ex]north west},
	"\text{\scriptsize exceeds}"]
\ar[from = 1-1,
	to = 1-4,
	crossing over,
	bend left = 0.3cm,
	ForestGreen,
	start anchor = east,
	end anchor = {[xshift = 1ex]west},
	shift right = 0.5ex,
	"\text{\scriptsize is equal to}"]
\ar[from = 1-1,
	to = 1-4,
	crossing over,
	bend right = 0.3cm,
	start anchor = south east,
	end anchor = {[xshift = 1ex]south west},
	Orange,
	"\text{\scriptsize is less than}"]
\ar[from = 2-1,
	to = 2-4,
	bend left = 0.65cm,
	blue,
	crossing over,
	start anchor = north east,
	end anchor = {[xshift = 1ex]north west},
	"\text{\scriptsize exceeds}"]
\ar[from = 2-1,
	crossing over,
	to = 2-4,
	ForestGreen,
	start anchor = east,
	bend left = 0.4cm,
	end anchor = {[xshift = 1ex]west},
	shift right = 0.5ex,
	"\text{\scriptsize is equal to}"]
\ar[from = 2-1,
	to = 2-4,
	crossing over,
	bend right = 0.3cm,
	start anchor = south east,
	end anchor = {[xshift = 1ex]south west},
	Orange,
	"\text{\scriptsize is less than}"]
\end{tikzcd}	

\flushleftright%
For the same reason, if $H$ exceeds $M$, $K$ exceeds $N$ (and if equal, equal; and if less than, less than).

\centering
\begin{tikzcd}[ampersand replacement = \&, row sep = 8ex, column sep = 2em, math mode = false]
{\begin{tikzpicture}
\draw[|-|, gray] (0,0) -- (1,0);
\node[inner sep = 2, anchor = east, gray] (Letter) at (0,0) {\tiny $G$};
\node[inner sep = 2, anchor = east] at (Letter.west) {\phantom{\footnotesize If}};
\end{tikzpicture}}
\ar[r,
	end anchor = {[xshift = 1ex]west},
	gray!25,
	shift right,
	""' {name = GA}]
\ar[d,
	Leftrightarrow,
	densely dotted,
	shift left = 3,
	gray!50,
	"\rotatebox{90}{\tiny Equimultiples}"',
	"\rotatebox{90}{\tiny of \scalebox{0.95}{$C;$ $A$}}"]
\ar[dd,
	Leftrightarrow,
	densely dotted,
	shift right,
	start anchor = {[xshift = -2ex]south},
	end anchor = {[xshift = -2ex]north},
	bend right = 2cm,
	gray!50,
	"\rotatebox{90}{\tiny Equimultiples}"',
	"\rotatebox{90}{\tiny of \scalebox{0.95}{$E;$ $A$}}"]
\&
	{\begin{tikzpicture}
	\draw[|-|, gray!50] (0,0) -- (0.5,0);
	\node[inner sep = 2, anchor = east, gray!50] (Letter) at (0,0) {\tiny $A$};
	\end{tikzpicture}}
	\ar[r,
		end anchor = {[xshift = 1ex]west},
		gray!25,
		shift right,
		""' {name = AB}]
	\&
		{\begin{tikzpicture}
		\draw[|-|, gray!50] (0,0) -- (0.3,0);
		\node[inner sep = 2, anchor = east, gray!50] at (0,0) {\tiny $B$};
		\end{tikzpicture}}
		\ar[r,
			leftarrow,
			gray!25,
			end anchor = {[xshift = 1ex]west},
			shift right,
			"" {name = LB}]
		\&
			{\begin{tikzpicture}
			\draw[|-|, gray] (0,0) -- (0.9,0);
			\node[inner sep = 2, anchor = east, gray] at (0,0) {\tiny $L$};
			\end{tikzpicture}}
			\ar[d,
				Leftrightarrow,
				densely dotted,
				shift left = 0.5,
				gray!50,
				"\rotatebox{-90}{\tiny Equimultiples}",
				"\rotatebox{-90}{\tiny of \scalebox{0.95}{$B;$ $D$}}"']
			\ar[dd,
				Leftrightarrow,
				densely dotted,
				shift left,
				start anchor = {[xshift = 2ex]south},
				end anchor = {[xshift = 2ex]north},
				bend left = 1.75cm,
				gray!50,
				"\rotatebox{-90}{\tiny Equimultiples}",
				"\rotatebox{-90}{\tiny of \scalebox{0.95}{$B;$ $F$}}"']\\
%
%
% Second Row
{\begin{tikzpicture}
\draw[|-|] (0,0) -- (1.333,0);
\node[inner sep = 2, anchor = east] (Letter) at (0,0) {\tiny $H$};
\node[inner sep = 2, anchor = east] at (Letter.west) {\footnotesize If};
\end{tikzpicture}}
\ar[r,
	end anchor = {[xshift = 1ex]west},
	gray!50,
	shift right,
	"" {name = HC}]
\ar[d,
	Leftrightarrow,
	densely dotted,
	gray,
	shift left = 3,
	"\rotatebox{90}{\tiny Equimultiples}"',
	"\rotatebox{90}{\tiny of \scalebox{0.95}{$E;$ $C$}}"]
\&
	{\begin{tikzpicture}
	\draw[|-|, black!60] (0,0) -- (0.667,0);
	\node[inner sep = 2, anchor = east, black!60] (Letter) at (0,0) {\tiny $C$};
	\end{tikzpicture}}
	\ar[r,
		gray!50,
		shift right = 0.5,
		bend left,
		start anchor = {[yshift = 0.1ex]east},
		end anchor = {[xshift = 1ex, yshift = 0.5]west},
		"" {name = CD}]
	\ar[r,
		gray!50,
		shift right = 1.5,
		bend right,
		start anchor = {[yshift = -0.1ex]east},
		end anchor = {[xshift = 1ex, yshift = -0.5]west},
		""' {name = CD2}]
	\&
		{\begin{tikzpicture}
		\draw[|-|, black!60] (0,0) -- (0.4,0);
		\node[inner sep = 2, anchor = east, black!60] at (0,0) {\tiny $D$};
		\end{tikzpicture}}
		\ar[r,
			leftarrow,
			gray!50,
			end anchor = {[xshift = 1ex]west},
			shift right,
			"" {name = MD}]
		\&
			{\begin{tikzpicture}
			\draw[|-|] (0,0) -- (1.2,0);
			\node[inner sep = 2, anchor = east] at (0,0) {\tiny $M$};
			\node[anchor = west, inner sep = 2] at (1.2,0) {\footnotesize ,};
			\end{tikzpicture}}
			\ar[d,
				Leftrightarrow,
				densely dotted,
				shift left = 0.5,
				gray,
				"\rotatebox{-90}{\tiny Equimultiples}",
				"\rotatebox{-90}{\tiny of \scalebox{0.95}{$D;$ $F$}}"']\\
%
%
% Third Row
{\begin{tikzpicture}
\draw[|-|] (0,0) -- (2,0);
\node[inner sep = 2, anchor = east] (Letter) at (0,0) {\tiny $K$};
\node[inner sep = 2, anchor = east] at (Letter.west) {\footnotesize then};
\end{tikzpicture}}
\ar[r,
	gray!50,
	end anchor = {[xshift = 1ex]west},
	shift right,
	"" {name = KE}]
\&
	{\begin{tikzpicture}
	\draw[|-|, black!60] (0,0) -- (1,0);
	\node[inner sep = 2, anchor = east, black!60] (Letter) at (0,0) {\tiny $E$};
	\end{tikzpicture}}
	\ar[r,
		gray!50,
		end anchor = {[xshift = 1ex]west},
		shift right,
		"" {name = EF}]
	\&
		{\begin{tikzpicture}
		\draw[|-|, black!60] (0,0) -- (0.6,0);
		\node[inner sep = 2, anchor = east, black!60] at (0,0) {\tiny $F$};
		\end{tikzpicture}}
		\ar[r,
			leftarrow,
			gray!50,
			end anchor = {[xshift = 1ex]west},
			shift right,
			"" {name = NF}]
		\&
			{\begin{tikzpicture}
			\draw[|-|] (0,0) -- (1.8,0);
			\node[inner sep = 2, anchor = east] at (0,0) {\tiny $N$};
			\node[anchor = west, inner sep = 2] at (1.8,0) {\footnotesize .};
			\end{tikzpicture}}
\ar[	from = AB,
	to = CD,
	Leftrightarrow,
	gray!25,
	"\rotatebox{90}{\scriptsize Same}"',
	"\rotatebox{90}{\scriptsize Ratio}"]
\ar[	from = CD2,
	to = EF,
	Leftrightarrow,
	gray!25,
	"\rotatebox{90}{\scriptsize Same}"',
	"\rotatebox{90}{\scriptsize Ratio}"]
\ar[from = 1-1,
	to = 1-4,
	bend left = 0.55cm,
	blue!25,
	start anchor = north east,
	end anchor = {[xshift = 1ex]north west},
	"\text{\scriptsize exceeds}"]
\ar[from = 1-1,
	to = 1-4,
	crossing over,
	bend left = 0.3cm,
	ForestGreen!25,
	start anchor = east,
	end anchor = {[xshift = 1ex]west},
	shift right = 0.5ex,
	"\text{\scriptsize is equal to}"]
\ar[from = 1-1,
	to = 1-4,
	crossing over,
	bend right = 0.3cm,
	start anchor = south east,
	end anchor = {[xshift = 1ex]south west},
	Orange!25,
	"\text{\scriptsize is less than}"]
\ar[from = 2-1,
	to = 2-4,
	bend left = 0.65cm,
	blue,
	crossing over,
	start anchor = north east,
	end anchor = {[xshift = 1ex]north west},
	"\text{\scriptsize exceeds}"]
\ar[from = 2-1,
	crossing over,
	to = 2-4,
	ForestGreen,
	start anchor = east,
	bend left = 0.4cm,
	end anchor = {[xshift = 1ex]west},
	shift right = 0.5ex,
	"\text{\scriptsize is equal to}"]
\ar[from = 2-1,
	to = 2-4,
	crossing over,
	bend right = 0.3cm,
	start anchor = south east,
	end anchor = {[xshift = 1ex]south west},
	Orange,
	"\text{\scriptsize is less than}"]
\ar[from = 3-1,
	to = 3-4,
	bend left = 0.65cm,
	blue,
	crossing over,
	start anchor = north east,
	end anchor = {[xshift = 1ex]north west},
	"\text{\scriptsize exceeds}"]
\ar[from = 3-1,
	crossing over,
	to = 3-4,
	ForestGreen,
	start anchor = east,
	bend left = 0.4cm,
	end anchor = {[xshift = 1ex]west},
	shift right = 0.5ex,
	"\text{\scriptsize is equal to}"]
\ar[from = 3-1,
	to = 3-4,
	crossing over,
	bend right = 0.3cm,
	start anchor = south east,
	end anchor = {[xshift = 1ex]south west},
	Orange,
	"\text{\scriptsize is less than}"]
\end{tikzcd}	

\flushleftright
And therefore, if $G$ exceeds $L$, $H$ exceeds $M$, and so $K$ exceeds $N$ (and if equal, equal; and if less than, less than). Hence, $A$ is to $B$ as $E$ is to $F$.

\begin{tcolorbox}[%
	colbacktitle=white,
	coltitle = black,
	center title,
	fonttitle = \bfseries,
	colbacklower = white,
	colback = white,
	width = 1.05\linewidth,
	lower separated=true,
	boxrule=.2mm]

\centering%
\begin{tikzcd}[ampersand replacement = \&, row sep = 8ex, column sep = 2em, math mode = false]
{\begin{tikzpicture}
\draw[|-|] (0,0) -- (1,0);
\node[inner sep = 2, anchor = east] (Letter) at (0,0) {\tiny $G$};
\node[inner sep = 2, anchor = east] at (Letter.west) {\footnotesize If};
\end{tikzpicture}}
\ar[r,
	end anchor = {[xshift = 1ex]west},
	gray,
	shift right,
	""' {name = GA}]
\ar[d,
	Leftrightarrow,
	densely dotted,
	shift left = 4.5,
	gray!50,
	"\rotatebox{90}{\tiny Equimultiples}"',
	"\rotatebox{90}{\tiny of \scalebox{0.95}{$C;$ $A$}}"]
\ar[dd,
	Leftrightarrow,
	densely dotted,
	shift right,
	start anchor = {[xshift = -2ex]south},
	bend right = 2.5cm,
	gray!50,
	"\rotatebox{90}{\tiny Equimultiples}"' {anchor = east},
	"\rotatebox{90}{\tiny of \scalebox{0.95}{$C;$ $A$}}" {anchor = west}]
\&
	{\begin{tikzpicture}
	\draw[|-|, black!60] (0,0) -- (0.5,0);
	\node[inner sep = 2, anchor = east, black!60] (Letter) at (0,0) {\tiny $A$};
	\end{tikzpicture}}
	\ar[r,
		end anchor = {[xshift = 1ex]west},
		gray!50,
		shift right,
		""' {name = AB}]
	\&
		{\begin{tikzpicture}
		\draw[|-|, black!60] (0,0) -- (0.3,0);
		\node[inner sep = 2, anchor = east, black!60] at (0,0) {\tiny $B$};
		\end{tikzpicture}}
		\ar[r,
			leftarrow,
			gray!50,
			end anchor = {[xshift = 1ex]west},
			shift right,
			"" {name = LB}]
		\&
			{\begin{tikzpicture}
			\draw[|-|] (0,0) -- (0.9,0);
			\node[inner sep = 2, anchor = east] at (0,0) {\tiny $L$};
			\node[inner sep = 2, anchor = west] at (0.9,0) {\footnotesize ,};
			\end{tikzpicture}}
			\ar[d,
				Leftrightarrow,
				densely dotted,
				shift left = 0,
				gray!50,
				"\rotatebox{-90}{\tiny Equimultiples}",
				"\rotatebox{-90}{\tiny of \scalebox{0.95}{$B;$ $D$}}"']
			\ar[dd,
				Leftrightarrow,
				densely dotted,
				shift left,
				start anchor = {[xshift = 2ex]south},
				end anchor = {[xshift = 2ex]north},
				bend left = 1.75cm,
				gray!50,
				"\rotatebox{-90}{\tiny Equimultiples}",
				"\rotatebox{-90}{\tiny of \scalebox{0.95}{$B;$ $F$}}"']\\
%
%
% Second Row
{\begin{tikzpicture}
\draw[|-|] (0,0) -- (1.333,0);
\node[inner sep = 2, anchor = east] (Letter) at (0,0) {\tiny $H$};
\node[inner sep = 2, anchor = east] at (Letter.west) {\footnotesize then};
\end{tikzpicture}}
\ar[r,
	end anchor = {[xshift = 1ex]west},
	gray!25,
	shift right,
	"" {name = HC}]
\ar[d,
	Leftrightarrow,
	densely dotted,
	gray!50,
	shift left = 3,
	"\rotatebox{90}{\tiny Equimultiples}"',
	"\rotatebox{90}{\tiny of \scalebox{0.95}{$E;$ $C$}}"]
\&
	{\begin{tikzpicture}
	\draw[|-|, gray!50] (0,0) -- (0.667,0);
	\node[inner sep = 2, anchor = east, gray!50] (Letter) at (0,0) {\tiny $C$};
	\end{tikzpicture}}
	\ar[r,
		gray!25,
		shift right = 0.5,
		bend left,
		start anchor = {[yshift = 0.1ex]east},
		end anchor = {[xshift = 1ex, yshift = 0.5]west},
		"" {name = CD}]
	\ar[r,
		gray!25,
		shift right = 1.5,
		bend right,
		start anchor = {[yshift = -0.1ex]east},
		end anchor = {[xshift = 1ex, yshift = -0.5]west},
		""' {name = CD2}]
	\&
		{\begin{tikzpicture}
		\draw[|-|, gray!50] (0,0) -- (0.4,0);
		\node[inner sep = 2, anchor = east, gray!50] at (0,0) {\tiny $D$};
		\end{tikzpicture}}
		\ar[r,
			leftarrow,
			gray!50,
			end anchor = {[xshift = 1ex]west},
			shift right,
			"" {name = MD}]
		\&
			{\begin{tikzpicture}
			\draw[|-|] (0,0) -- (1.2,0);
			\node[inner sep = 2, anchor = east] at (0,0) {\tiny $M$};
			\node[anchor = west, inner sep = 2] at (1.2,0) {\footnotesize ,};
			\end{tikzpicture}}
			\ar[d,
				Leftrightarrow,
				densely dotted,
				shift left = 0,
				gray!50,
				"\rotatebox{-90}{\tiny Equimultiples}",
				"\rotatebox{-90}{\tiny of \scalebox{0.95}{$D;$ $F$}}"']\\
%
%
% Third Row
{\begin{tikzpicture}
\draw[|-|] (0,0) -- (2,0);
\node[inner sep = 2, anchor = east] (Letter) at (0,0) {\tiny $K$};
\node[inner sep = 2, anchor = east] at (Letter.west) {\footnotesize and so};
\end{tikzpicture}}
\ar[r,
	gray!50,
	end anchor = {[xshift = 1ex]west},
	shift right,
	"" {name = KE}]
\&
	{\begin{tikzpicture}
	\draw[|-|, black!60] (0,0) -- (1,0);
	\node[inner sep = 2, anchor = east, black!60] (Letter) at (0,0) {\tiny $E$};
	\end{tikzpicture}}
	\ar[r,
		gray!50,
		end anchor = {[xshift = 1ex]west},
		shift right,
		"" {name = EF}]
	\&
		{\begin{tikzpicture}
		\draw[|-|, black!60] (0,0) -- (0.6,0);
		\node[inner sep = 2, anchor = east, black!60] at (0,0) {\tiny $F$};
		\end{tikzpicture}}
		\ar[r,
			leftarrow,
			gray!50,
			end anchor = {[xshift = 1ex]west},
			shift right,
			"" {name = NF}]
		\&
			{\begin{tikzpicture}
			\draw[|-|] (0,0) -- (1.8,0);
			\node[inner sep = 2, anchor = east] at (0,0) {\tiny $N$};
			\node[anchor = west, inner sep = 2] at (1.8,0) {\footnotesize .};
			\end{tikzpicture}}
\ar[	from = AB,
	to = CD,
	Leftrightarrow,
	gray!25,
	"\rotatebox{90}{\scriptsize Same}"',
	"\rotatebox{90}{\scriptsize Ratio}"]
\ar[	from = CD2,
	to = EF,
	Leftrightarrow,
	gray!25,
	"\rotatebox{90}{\scriptsize Same}"',
	"\rotatebox{90}{\scriptsize Ratio}"]
\ar[from = 1-1,
	to = 1-4,
	bend left = 0.55cm,
	blue,
	start anchor = north east,
	end anchor = {[xshift = 1ex]north west},
	"\text{\scriptsize exceeds}"]
\ar[from = 1-1,
	to = 1-4,
	crossing over,
	bend left = 0.3cm,
	ForestGreen,
	start anchor = east,
	end anchor = {[xshift = 1ex]west},
	shift right = 0.5ex,
	"\text{\scriptsize is equal to}"]
\ar[from = 1-1,
	to = 1-4,
	crossing over,
	bend right = 0.3cm,
	start anchor = south east,
	end anchor = {[xshift = 1ex]south west},
	Orange,
	"\text{\scriptsize is less than}"]
\ar[from = 2-1,
	to = 2-4,
	bend left = 0.65cm,
	blue,
	crossing over,
	start anchor = north east,
	end anchor = {[xshift = 1ex]north west},
	"\text{\scriptsize exceeds}"]
\ar[from = 2-1,
	crossing over,
	to = 2-4,
	ForestGreen,
	start anchor = east,
	bend left = 0.4cm,
	end anchor = {[xshift = 1ex]west},
	shift right = 0.5ex,
	"\text{\scriptsize is equal to}"]
\ar[from = 2-1,
	to = 2-4,
	crossing over,
	bend right = 0.3cm,
	start anchor = south east,
	end anchor = {[xshift = 1ex]south west},
	Orange,
	"\text{\scriptsize is less than}"]
\ar[from = 3-1,
	to = 3-4,
	bend left = 0.65cm,
	blue,
	crossing over,
	start anchor = north east,
	end anchor = {[xshift = 1ex]north west},
	"\text{\scriptsize exceeds}"]
\ar[from = 3-1,
	crossing over,
	to = 3-4,
	ForestGreen,
	start anchor = east,
	bend left = 0.4cm,
	end anchor = {[xshift = 1ex]west},
	shift right = 0.5ex,
	"\text{\scriptsize is equal to}"]
\ar[from = 3-1,
	to = 3-4,
	crossing over,
	bend right = 0.3cm,
	start anchor = south east,
	end anchor = {[xshift = 1ex]south west},
	Orange,
	"\text{\scriptsize is less than}"]
\end{tikzcd}	
\tcblower
Therefore,

\centering%
\begin{tikzcd}[ampersand replacement = \&, row sep = 8ex, column sep = 2em, math mode = false]
{\begin{tikzpicture}
\draw[|-|, gray!50] (0,0) -- (1,0);
\node[inner sep = 2, anchor = east, gray!50] (Letter) at (0,0) {\tiny $G$};
\end{tikzpicture}}
\ar[r,
	end anchor = {[xshift = 1ex]west},
	gray!37.5,
	shift right,
	""' {name = GA}]
\ar[d,
	Leftrightarrow,
	densely dotted,
	shift left = 1.5,
	gray!25,
	"\rotatebox{90}{\tiny Equimultiples}"',
	"\rotatebox{90}{\tiny of \scalebox{0.95}{$K;$ $A$}}"]
\&
	{\begin{tikzpicture}
	\draw[|-|] (0,0) -- (0.5,0);
	\node[inner sep = 2, anchor = east] (Letter) at (0,0) {\tiny $A$};
	\node[inner sep = 2, anchor = east] at (Letter.west) {\footnotesize As};
	\end{tikzpicture}}
	\ar[r,
		end anchor = {[xshift = 1ex]west},
		gray,
		shift right,
		"\tiny to" {black},
		""' {name = AB, pos = 0.51}]
	\&
		{\begin{tikzpicture}
		\draw[|-|] (0,0) -- (0.3,0);
		\node[inner sep = 2, anchor = east] at (0,0) {\tiny $B$};
		\node[inner sep = 2, anchor = west] at (0.3,0) {\footnotesize ,};
		\end{tikzpicture}}
		\ar[r,
			leftarrow,
			gray!37.5,
			end anchor = {[xshift = 1ex]west},
			shift right,
			"" {name = LB}]
		\&
			{\begin{tikzpicture}
			\draw[|-|, gray!50] (0,0) -- (0.9,0);
			\node[inner sep = 2, anchor = east, gray!50] at (0,0) {\tiny $L$};
			\end{tikzpicture}}
			\ar[d,
				Leftrightarrow,
				densely dotted,
				shift left = 1.0,
				gray!25,
				"\rotatebox{-90}{\tiny Equimultiples}",
				"\rotatebox{-90}{\tiny of \scalebox{0.95}{$B;$ $F$}}"']\\
%
%
% Second Row
{\begin{tikzpicture}
\draw[|-|, gray!50] (0,0) -- (2,0);
\node[inner sep = 2, anchor = east, gray!50] (Letter) at (0,0) {\tiny $K$};
\end{tikzpicture}}
\ar[r,
	gray!25,
	end anchor = {[xshift = 1ex]west},
	shift right,
	"" {name = KE}]
\&
	{\begin{tikzpicture}
	\draw[|-|] (0,0) -- (1,0);
	\node[inner sep = 2, anchor = east] (Letter) at (0,0) {\tiny $E$};
	\node[inner sep = 2, anchor = east] at (Letter.west) {\footnotesize so};
	\end{tikzpicture}}
	\ar[r,
		gray,
		end anchor = {[xshift = 1ex]west},
		shift right,
		"\tiny to" {black},
		"\phantom{\tiny to}" {name = EF, pos = 0.49}]
	\&
		{\begin{tikzpicture}
		\draw[|-|] (0,0) -- (0.6,0);
		\node[inner sep = 2, anchor = east] at (0,0) {\tiny $F$};
		\node[inner sep = 2, anchor = west] at (0.6,0) {\footnotesize .};
		\end{tikzpicture}}
		\ar[r,
			leftarrow,
			gray!25,
			end anchor = {[xshift = 1ex]west},
			shift right,
			"" {name = NF}]
		\&
			{\begin{tikzpicture}
			\draw[|-|, gray!50] (0,0) -- (1.8,0);
			\node[inner sep = 2, anchor = east, gray!50] at (0,0) {\tiny $N$};
			\end{tikzpicture}}
\ar[from = 1-1,
	to = 1-4,
	bend left = 0.55cm,
	blue!25,
	start anchor = north east,
	end anchor = {[xshift = 1ex]north west},
	"\text{\scriptsize exceeds}"]
\ar[from = 1-1,
	to = 1-4,
	bend left = 0.3cm,
	ForestGreen!25,
	start anchor = east,
	end anchor = {[xshift = 1ex]west},
	shift right = 0.5ex,
	"\text{\scriptsize is equal to}"]
\ar[from = 1-1,
	to = 1-4,
	bend right = 0.3cm,
	start anchor = south east,
	end anchor = {[xshift = 1ex]south west},
	Orange!12.5,
	"\text{\scriptsize is less than}"]
\ar[from = 2-1,
	to = 2-4,
	bend left = 0.65cm,
	blue!12.5,
	start anchor = north east,
	end anchor = {[xshift = 1ex]north west},
	"\text{\scriptsize exceeds}"]
\ar[from = 2-1,
	to = 2-4,
	ForestGreen!12.5,
	start anchor = east,
	bend left = 0.4cm,
	end anchor = {[xshift = 1ex]west},
	shift right = 0.5ex,
	"\text{\scriptsize is equal to}"]
\ar[from = 2-1,
	to = 2-4,
	bend right = 0.3cm,
	start anchor = south east,
	end anchor = {[xshift = 1ex]south west},
	Orange!25,
	"\text{\scriptsize is less than}"]
\ar[	from = AB,
	to = EF,
	Leftrightarrow,
	crossing over,
	gray,
	"\rotatebox{90}{\scriptsize Same}"',
	"\rotatebox{90}{\scriptsize Ratio}"]
\ar[from = 1-2,
	to = 2-2,
	Leftrightarrow,
	crossing over,
	gray!50,
	densely dotted,
	shift left = 4,
	"\rotatebox{90}{\tiny Analogous}"',
	"\rotatebox{90}{\tiny Quantities}"]
\ar[from = 1-3,
	to = 2-3,
	Leftrightarrow,
	crossing over,
	gray!50,
	densely dotted,
	shift left,
	"\rotatebox{90}{\tiny Analogous}"',
	"\rotatebox{90}{\tiny Quantities}"]
\end{tikzcd}	
\end{tcolorbox}

\flushleftright
\subsection{Proposition 13}
Proposition 13 is very similar to Proposition 11, but here, $A$ has the same ratio to $B$ as $C$ has to $D$, but $C$ has a greater ratio to $D$ than $E$ has to $F$.

\begin{tcolorbox}[%
	colbacktitle=white,
	coltitle = black,
	center title,
	halign upper=center,
	halign lower=center,
	fonttitle = \bfseries,
	colbacklower = white,
	colback = white,
	width = \linewidth,
	lower separated=true,
	sidebyside,
	boxrule=.2mm]
\begin{tikzcd}[ampersand replacement = \&, row sep = 8ex, column sep = 3em, math mode = false]
{\begin{tikzpicture}
\draw[|-|] (0,0) -- (0.6,0);
\node[inner sep = 2, anchor = east] (Letter) at (0,0) {\tiny $A$};
\node[inner sep = 2, anchor = east] at (Letter.west) {\tiny As};
\node[inner sep = 2, anchor = west] at (0.7,0) {\tiny is};
\end{tikzpicture}}
\ar[r,
	shift right,
	"\tiny to",
	""' {name = AB}]
\ar[d,
	Leftrightarrow,
	gray!50,
	densely dotted,
	shift left = 2,
	"\rotatebox{90}{\tiny Analogous}"',
	"\rotatebox{90}{\tiny Quantities}"]
\&
	{\begin{tikzpicture}
	\draw[|-|] (0,0) -- (0.2,0);
	\node[inner sep = 2, anchor = east] at (0,0) {\tiny $B$};
	\node[inner sep = 2, anchor = west] at (0.2,0) {\tiny ,};
	\end{tikzpicture}}
	\ar[d,
		Leftrightarrow,
		gray!50,
		densely dotted,
		shift left,
		"\rotatebox{90}{\tiny Analogous}"',
		"\rotatebox{90}{\tiny Quantities}"]\\
%
%
% Second Row
{\begin{tikzpicture}
\draw[|-|] (0,0) -- (0.9,0);
\node[inner sep = 2, anchor = east] (Letter) at (0,0) {\tiny $C$};
\node[inner sep = 2, anchor = east] at (Letter.west) {\tiny so};
\node[inner sep = 2, anchor = west] at (1,0) {\tiny is};

\end{tikzpicture}}
\ar[r,
	shift right,
	"\tiny to",
	"\phantom{\tiny to}" {name = CD}]
\&
	{\begin{tikzpicture}
	\draw[|-|] (0,0) -- (0.3,0);
	\node[inner sep = 2, anchor = east] at (0,0) {\tiny $D$};
	\node[inner sep = 2, anchor = west] at (0.3,0) {\tiny .};
	\end{tikzpicture}}\\[-8ex]
%
%
% Third Row
{\begin{tikzpicture}
\draw[|-|] (0,0) -- (0.9,0);
\node[inner sep = 2, anchor = east] (Letter) at (0,0) {\tiny $C$};
\node[inner sep = 2, anchor = east] at (Letter.west) {\tiny But};
\node[inner sep = 2, anchor = west] at (0.9,0) {\tiny has,};
\end{tikzpicture}}
\ar[r,
	shift right = 2,
	start anchor = {[xshift = -1ex]east},
	end anchor = {[xshift = 1ex]west},
	"\tiny \textbf{to}",
	""' {name = CD2, pos = 0.44}]
\&
	{\begin{tikzpicture}
	\draw[|-|] (0,0) -- (0.3,0);
	\node[inner sep = 2, anchor = east] at (0,0) {\tiny $D$};
	\node[inner sep = 2, anchor = west] at (0.3,0) {\tiny ,};
	\end{tikzpicture}}\\
%
%
% Fourth Row
{\begin{tikzpicture}
\draw[|-|] (0,0) -- (0.5,0);
\node[inner sep = 2, anchor = east] (Letter) at (0,0) {\tiny $E$};
\node[inner sep = 2, anchor = west] at (0.5,0) {\tiny has,};
\end{tikzpicture}}
\ar[r,
	shift right = 2,
	start anchor = {[xshift = 1ex]east},
	end anchor = {[xshift = -1]west},
	thin,
	gray!75,
	"\scalebox{0.9}{\tiny \textit{to}}",
	"\phantom{\tiny to}" {name = EF, pos = 0.56}]
\&
	{\begin{tikzpicture}
	\draw[|-|] (0,0) -- (0.4,0);
	\node[inner sep = 2, anchor = east] at (0,0) {\tiny $F$};
	\node[inner sep = 2, anchor = west] at (0.4,0) {\tiny .};
	\end{tikzpicture}}
\ar[	from = AB,
	to = CD,
	Leftrightarrow,
	gray,
	"\rotatebox{90}{\scriptsize Same}"',
	"\rotatebox{90}{\scriptsize Ratio}"]
\ar[from = CD2,
	to = EF,
	Rightarrow,
	gray,
	"\tiny a Greater" sloped,
	"\tiny Ratio than"' sloped]
\end{tikzcd}	

\tcblower

\flushleftright%
The claim is:

\centering%
\begin{tikzcd}[ampersand replacement = \&, row sep = 8ex, column sep = 3em, math mode = false]
{\begin{tikzpicture}
\draw[|-|] (0,0) -- (0.6,0);
\node[inner sep = 2, anchor = east] (Letter) at (0,0) {\tiny $A$};
\node[inner sep = 2, anchor = west] at (0.6,0) {\tiny has,};
\end{tikzpicture}}
\ar[r,
	shift right = 2,
	start anchor = {[xshift = -1ex]east},
	end anchor = {[xshift = 1ex]west},
	"\tiny \textbf{to}",
	""' {name = CD2, pos = 0.505}]
\&
	{\begin{tikzpicture}
	\draw[|-|] (0,0) -- (0.2,0);
	\node[inner sep = 2, anchor = east] at (0,0) {\tiny $B$};
	\node[inner sep = 2, anchor = west] at (0.2,0) {\tiny ,};
	\end{tikzpicture}}\\
%
%
% Fourth Row
{\begin{tikzpicture}
\draw[|-|] (0,0) -- (0.5,0);
\node[inner sep = 2, anchor = east] (Letter) at (0,0) {\tiny $E$};
\node[inner sep = 2, anchor = west] at (0.5,0) {\tiny has,};
\end{tikzpicture}}
\ar[r,
	shift right = 2,
	start anchor = {[xshift = 1ex]east},
	end anchor = {[xshift = -1]west},
	thin,
	gray!75,
	"\scalebox{0.9}{\tiny \textit{to}}",
	"\phantom{\tiny to}" {name = EF, pos = 0.495}]
\&
	{\begin{tikzpicture}
	\draw[|-|] (0,0) -- (0.4,0);
	\node[inner sep = 2, anchor = east] at (0,0) {\tiny $F$};
	\node[inner sep = 2, anchor = west] at (0.4,0) {\tiny .};
	\end{tikzpicture}}
\ar[from = CD2,
	to = EF,
	Rightarrow,
	gray,
	"\tiny a Greater" sloped,
	"\tiny Ratio than"' sloped]
\end{tikzcd}	
\end{tcolorbox}

\flushleftright%
The argument is very similar to before. This time, since $C$ has a greater ratio to $D$ than $E$ has to $F$, we can take equimultiples of $C$ and $E$ and other equal multiples of $D$ and $F$ such that the multiple of $E$ is less than the multiple of $F$ but the multiple of $C$ exceeds the multiple of $D$.

\centering%
\begin{tikzcd}[ampersand replacement = \&, row sep = 8ex, column sep = 2em, math mode = false]
\&
{\begin{tikzpicture}
\draw[|-|] (0,0) -- (0.6,0);
\node[inner sep = 2, anchor = east]  at (0,0) {\tiny $A$};
\end{tikzpicture}}
\ar[r,
	shift right,
	"\tiny to",
	""' {name = AB, pos = 0.55}]
\ar[d,
	Leftrightarrow,
	gray!50,
	densely dotted,
	shift left = 2,
	"\rotatebox{90}{\tiny Analogous}"',
	"\rotatebox{90}{\tiny Quantities}"]
\&
	{\begin{tikzpicture}
	\draw[|-|] (0,0) -- (0.2,0);
	\node[inner sep = 2, anchor = east] at (0,0) {\tiny $B$};
	\node[inner sep = 2, anchor = west] at (0.2,0) {\tiny ,};
	\end{tikzpicture}}
	\ar[d,
		Leftrightarrow,
		gray!50,
		densely dotted,
		shift left,
		"\rotatebox{90}{\tiny Analogous}"',
		"\rotatebox{90}{\tiny Quantities}"]\\
%%Real Second Row
{\begin{tikzpicture}
\draw[|-|] (0,0) -- (1.8,0);
\node[inner sep = 2, anchor = east] (Letter) at (0,0) {\tiny $G$};
\end{tikzpicture}}
\ar[r,
	shift right,
	""' {name = GC}]
\ar[d,
	gray,
	Leftrightarrow,
	densely dotted,
	shift left,
	"\rotatebox{90}{\tiny Equimultiples}"',
	"\rotatebox{90}{\tiny of \scalebox{0.95}{$E;$ $C$}}"]
\&
	{\begin{tikzpicture}
	\draw[|-|] (0,0) -- (0.9,0);
	\node[inner sep = 2, anchor = east] (Letter) at (0,0) {\tiny $C$};
	\end{tikzpicture}}
	\ar[r,
		shift right,
		start anchor = {[xshift = -0.5ex]east},
		end anchor = {[xshift = 1ex]west},
		"\tiny \textbf{to}" {name = CD},
		""' {name = CD2, pos = 0.445}]
	\&
		{\begin{tikzpicture}
		\draw[|-|] (0,0) -- (0.3,0);
		\node[inner sep = 2, anchor = east] at (0,0) {\tiny $D$};
		\end{tikzpicture}}
		\ar[r,
			leftarrow,
			shift right,
			""' {name = KD}]
		\&
			{\begin{tikzpicture}
			\draw[|-|] (0,0) -- (1.2,0);
			\node[inner sep = 2, anchor = east] at (0,0) {\tiny $K$};
			\end{tikzpicture}}
			\ar[d,
				gray,
				Leftrightarrow,
				densely dotted,
				"\rotatebox{-90}{\tiny Equimultiples}",
				"\rotatebox{-90}{\tiny of \scalebox{0.95}{$K;$ $L$}}"']\\
%
%
% Second Row
{\begin{tikzpicture}
\draw[|-|] (0,0) -- (1,0);
\node[inner sep = 2, anchor = east] (Letter) at (0,0) {\tiny $H$};
\end{tikzpicture}}
\ar[rrr,
	bend right = 0.5cm,
	start anchor = {[yshift = 0.5ex]south east},
	end anchor = {[yshift = 0.5ex, xshift = 1ex]south west},
	Orange,
	"\tiny Does not Exceed"]
\ar[r,
	shift right,
	"" {name = HE}]
\&
	{\begin{tikzpicture}
	\draw[|-|] (0,0) -- (0.5,0);
	\node[inner sep = 2, anchor = east] (Letter) at (0,0) {\tiny $E$};
	\end{tikzpicture}}
	\ar[r,
		shift right,
		start anchor = {[xshift = 0.5ex]east},
		end anchor = {[xshift = 0]west},
		thin,
		gray!75,
		"\scalebox{0.9}{\tiny \textit{to}}",
		"\phantom{\tiny to}" {name = EF, pos = 0.555}]
	\&
		{\begin{tikzpicture}
		\draw[|-|] (0,0) -- (0.4,0);
		\node[inner sep = 2, anchor = east] at (0,0) {\tiny $F$};
		\end{tikzpicture}}
		\ar[r,
			leftarrow,
			shift right,
			"" {name = LF}]
		\&
			{\begin{tikzpicture}
			\draw[|-|] (0,0) -- (1.2,0);
			\node[inner sep = 2, anchor = east] at (0,0) {\tiny $L$};
			\end{tikzpicture}}
\ar[	from = AB,
	to = CD,
	Leftrightarrow,
	gray!50,
	"\rotatebox{90}{\scriptsize Same}"',
	"\rotatebox{90}{\scriptsize Ratio}"]
\ar[from = CD2,
	to = EF,
	Rightarrow,
	gray,
	"\tiny Greater" sloped,
	"\tiny Ratio than"' sloped]
\ar[from = 2-1,
	to = 2-4,
	bend left = 0.5cm,
	start anchor = {[yshift = -0.5ex]north east},
	end anchor = {[yshift = -0.5ex, xshift = 1ex]north west},
	white,
	-,
	line width = 1.25mm]
\ar[from = 2-1,
	to = 2-4,
	bend left = 0.5cm,
	start anchor = {[yshift = -0.5ex]north east},
	end anchor = {[yshift = -0.5ex, xshift = 1ex]north west},
	blue,
	"\tiny Exceeds" {fill = white, inner sep = 2, yshift = 1.33pt, fill opacity = 0.7, text = blue, text opacity = 1}]
\end{tikzcd}	

\flushleftright%
We then take an equimultiple of $A$ and another equimultiple of $B$---each the same as the corresponding multiples of $C$ and $E$, or of $D$ and $F$. Since $G$ exceeds $K$ and the ratios are the same, the multiple of $A$ exceeds the multiple of $B$. And therefore the multiple of $A$ exceeds the multiple of $B$, but $H$ does not exceed $L$, and so $A$ has, to $B$ a greater ratio than $H$ has to $L$.

\begin{tcolorbox}[%
	colbacktitle=white,
	coltitle = black,
	center title,
	halign upper=center,
	fonttitle = \bfseries,
	colbacklower = white,
	colback = white,
	width = 1.05\linewidth,
	lower separated=true,
	boxrule=.2mm]

\hspace{-0.375cm}%
\begin{tikzcd}[ampersand replacement = \&, row sep = 8ex, column sep = 2em, math mode = false]
{\begin{tikzpicture}
\draw[|-|] (0,0) -- (1.2,0);
\node[inner sep = 2, anchor = east] (Letter) at (0,0) {\tiny $M$};
\end{tikzpicture}}
\ar[r,
	shift right,
	""' {name = MA}]
\ar[d,
	gray,
	Leftrightarrow,
	densely dotted,
	shift left,
	"\rotatebox{90}{\tiny Equimultiples}"',
	"\rotatebox{90}{\tiny of \scalebox{0.95}{$C;$ $A$}}"]
\ar[dd,
	Leftrightarrow,
	densely dotted,
	shift right,
	start anchor = {[xshift = -2.5ex]south},
	end anchor = {[xshift = -2.5ex]north},
	out = -135,
	in = 135,
	distance = 1.5cm,
	gray,
	"\rotatebox{90}{\tiny Equimultiples}"',
	"\rotatebox{90}{\tiny of \scalebox{0.95}{$E;$ $A$}}"]
\&
	{\begin{tikzpicture}
	\draw[|-|] (0,0) -- (0.6,0);
	\node[inner sep = 2, anchor = east] (Letter) at (0,0) {\tiny $A$};
	\end{tikzpicture}}
	\ar[r,
		shift right,
		""' {name = AB, pos = 0.55}]
	\&
		{\begin{tikzpicture}
		\draw[|-|] (0,0) -- (0.2,0);
		\node[inner sep = 2, anchor = east] at (0,0) {\tiny $B$};
		\end{tikzpicture}}
		\ar[r,
			leftarrow,
			shift right,
			""' {name = NB}]
		\&
			{\begin{tikzpicture}
			\draw[|-|] (0,0) -- (0.6,0);
			\node[inner sep = 2, anchor = east] at (0,0) {\tiny $N$};
			\end{tikzpicture}}
			\ar[d,
				gray,
				Leftrightarrow,
				densely dotted,
				"\rotatebox{-90}{\tiny Equimultiples}",
				"\rotatebox{-90}{\tiny of \scalebox{0.95}{$B;$ $D$}}"']
			\ar[dd,
				Leftrightarrow,
				densely dotted,
				shift left,
				start anchor = {[xshift = 2ex]south},
				end anchor = {[xshift = 2ex]north},
				bend left = 1.75cm,
				gray,
				"\rotatebox{-90}{\tiny Equimultiples}",
				"\rotatebox{-90}{\tiny of \scalebox{0.95}{$B;$ $F$}}"']\\
%
%
% Second Row
{\begin{tikzpicture}
\draw[|-|] (0,0) -- (1.8,0);
\node[inner sep = 2, anchor = east] (Letter) at (0,0) {\tiny $G$};
\end{tikzpicture}}

\ar[r,
	shift right,
	""' {name = GC}]
\ar[d,
	gray,
	Leftrightarrow,
	densely dotted,
	shift left,
	"\rotatebox{90}{\tiny Equimultiples}"',
	"\rotatebox{90}{\tiny of \scalebox{0.95}{$E;$ $C$}}"]
\&
	{\begin{tikzpicture}
	\draw[|-|] (0,0) -- (0.9,0);
	\node[inner sep = 2, anchor = east] (Letter) at (0,0) {\tiny $C$};
	\end{tikzpicture}}
	\ar[r,
		shift right,
		start anchor = {[xshift = -0.5ex]east},
		end anchor = {[xshift = 1ex]west},
		"\tiny \textbf{to}",
		"\phantom{\tiny \textbf{to}}" {name = CD, pos = 0.48},
		""' {name = CD2, pos = 0.445}]
	\&
		{\begin{tikzpicture}
		\draw[|-|] (0,0) -- (0.3,0);
		\node[inner sep = 2, anchor = east] at (0,0) {\tiny $D$};
		\end{tikzpicture}}
		\ar[r,
			leftarrow,
			shift right,
			""' {name = KD}]
		\&
			{\begin{tikzpicture}
			\draw[|-|] (0,0) -- (1.2,0);
			\node[inner sep = 2, anchor = east] at (0,0) {\tiny $K$};
			\end{tikzpicture}}
			\ar[d,
				gray,
				Leftrightarrow,
				densely dotted,
				"\rotatebox{-90}{\tiny Equimultiples}",
				"\rotatebox{-90}{\tiny of \scalebox{0.95}{$K;$ $L$}}"']\\
%
%
% Second Row
{\begin{tikzpicture}
\draw[|-|] (0,0) -- (1,0);
\node[inner sep = 2, anchor = east] (Letter) at (0,0) {\tiny $H$};
\end{tikzpicture}}
\ar[r,
	shift right,
	"" {name = HE}]
\&
	{\begin{tikzpicture}
	\draw[|-|] (0,0) -- (0.5,0);
	\node[inner sep = 2, anchor = east] (Letter) at (0,0) {\tiny $E$};
	\end{tikzpicture}}
	\ar[r,
		shift right,
		start anchor = {[xshift = 0.5ex]east},
		end anchor = {[xshift = 0]west},
		thin,
		gray!75,
		"\scalebox{0.9}{\tiny \textit{to}}",
		"\phantom{\tiny to}" {name = EF, pos = 0.555}]
	\&
		{\begin{tikzpicture}
		\draw[|-|] (0,0) -- (0.4,0);
		\node[inner sep = 2, anchor = east] at (0,0) {\tiny $F$};
		\end{tikzpicture}}
		\ar[r,
			leftarrow,
			shift right,
			"" {name = LF}]
		\&
			{\begin{tikzpicture}
			\draw[|-|] (0,0) -- (1.2,0);
			\node[inner sep = 2, anchor = east] at (0,0) {\tiny $L$};
			\end{tikzpicture}}
\ar[from = CD2,
	to = EF,
	Rightarrow,
	gray!50,
	"\tiny Greater" sloped,
	"\tiny Ratio than"' sloped]
\ar[from = 1-3,
	to = 2-3,
	Leftrightarrow,
	crossing over,
	gray!25,
	densely dotted,
	shift left,
	"\rotatebox{90}{\tiny Analogous}"',
	"\rotatebox{90}{\tiny Quantities}"]
\ar[from = 1-2,
	to = 2-2,
	Leftrightarrow,
	crossing over,
	gray!25,
	densely dotted,
	shift left,
	"\rotatebox{90}{\tiny Analogous}"',
	"\rotatebox{90}{\tiny Quantities}"]
\ar[	from = AB,
	to = CD,
	Leftrightarrow,
	crossing over,
	gray!50,
	"\rotatebox{90}{\scriptsize Same}"',
	"\rotatebox{90}{\scriptsize Ratio}"]
\ar[from = 3-1,
	to = 3-4,
	bend right = 0.5cm,
	start anchor = {[yshift = 0.5ex]south east},
	end anchor = {[yshift = 0.5ex, xshift = 1ex]south west},
	Orange,
	"\tiny Does not Exceed"]
\ar[from = 2-1,
	to = 2-4,
	bend left = 0.5cm,
	start anchor = {[yshift = -0.5ex]north east},
	end anchor = {[yshift = -0.5ex, xshift = 1ex]north west},
	crossing over,
	blue,
	"\tiny Exceeds"]
\ar[from = 1-1,
	to = 1-4,
	bend left = 0.5cm,
	start anchor = {[yshift = -0.5ex]north east},
	end anchor = {[yshift = -0.5ex, xshift = 1ex]north west},
	crossing over,
	blue,
	"\tiny Exceeds"]
\end{tikzcd}	
\tcblower
Therefore,

\centering
\begin{tikzcd}[ampersand replacement = \&, row sep = 8ex, column sep = 2em, math mode = false]
{\begin{tikzpicture}
\draw[|-|, gray!50] (0,0) -- (1.2,0);
\node[inner sep = 2, anchor = east, gray!50] (Letter) at (0,0) {\tiny $M$};
\end{tikzpicture}}
\ar[r,
	gray!25,
	shift right,
	""' {name = MA}]
\ar[d,
	gray!25,
	Leftrightarrow,
	densely dotted,
	shift left,
	"\rotatebox{90}{\tiny Equimultiples}"',
	"\rotatebox{90}{\tiny of \scalebox{0.95}{$C;$ $A$}}"]
\&
	{\begin{tikzpicture}
	\draw[|-|] (0,0) -- (0.6,0);
	\node[inner sep = 2, anchor = east] (Letter) at (0,0) {\tiny $A$};
	\end{tikzpicture}}
	\ar[r,
		shift right = 2,
		"\tiny \textbf{to}",
		""' {name = AB, pos = 0.53}]
	\&
		{\begin{tikzpicture}
		\draw[|-|] (0,0) -- (0.2,0);
		\node[inner sep = 2, anchor = east] at (0,0) {\tiny $B$};
		\end{tikzpicture}}
		\ar[r,
			gray!25,
			leftarrow,
			shift right,
			""' {name = NB}]
		\&
			{\begin{tikzpicture}
			\draw[|-|, gray!50] (0,0) -- (0.6,0);
			\node[inner sep = 2, anchor = east, gray!50] at (0,0) {\tiny $N$};
			\end{tikzpicture}}
			\ar[d,
				gray!25,
				Leftrightarrow,
				densely dotted,
				"\rotatebox{-90}{\tiny Equimultiples}",
				"\rotatebox{-90}{\tiny of \scalebox{0.95}{$B;$ $D$}}"']\\
%
%
% Second Row
{\begin{tikzpicture}
\draw[|-|, gray!50] (0,0) -- (1,0);
\node[inner sep = 2, anchor = east, gray!50] (Letter) at (0,0) {\tiny $H$};
\end{tikzpicture}}
\ar[r,
	gray!25,
	shift right,
	"" {name = HE}]
\&
	{\begin{tikzpicture}
	\draw[|-|] (0,0) -- (0.5,0);
	\node[inner sep = 2, anchor = east] (Letter) at (0,0) {\tiny $E$};
	\end{tikzpicture}}
	\ar[r,
		shift right,
		start anchor = {[xshift = 0.5ex]east},
		end anchor = {[xshift = 0]west},
		thin,
		gray!75,
		"\scalebox{0.9}{\tiny \textit{to}}",
		"\phantom{\tiny to}" {name = EF, pos = 0.555}]
	\&
		{\begin{tikzpicture}
		\draw[|-|] (0,0) -- (0.4,0);
		\node[inner sep = 2, anchor = east] at (0,0) {\tiny $F$};
		\end{tikzpicture}}
		\ar[r,
			gray!25,
			leftarrow,
			shift right,
			"" {name = LF}]
		\&
			{\begin{tikzpicture}
			\draw[|-|, gray!50] (0,0) -- (1.2,0);
			\node[inner sep = 2, anchor = east, gray!50] at (0,0) {\tiny $L$};
			\end{tikzpicture}}
\ar[from = AB,
	to = EF,
	Rightarrow,
	gray,
	"\tiny Greater" sloped,
	"\tiny Ratio than"' sloped]
\ar[from = 2-1,
	to = 2-4,
	bend right = 0.5cm,
	start anchor = {[yshift = 0.5ex]south east},
	end anchor = {[yshift = 0.5ex, xshift = 1ex]south west},
	Orange!25,
	"\tiny Does not Exceed"]
\ar[from = 1-1,
	to = 1-4,
	bend left = 0.5cm,
	start anchor = {[yshift = -0.5ex]north east},
	end anchor = {[yshift = -0.5ex, xshift = 1ex]north west},
	crossing over,
	blue!25,
	"\tiny Exceeds"]
\end{tikzcd}	
\end{tcolorbox}

\flushleftright%
Though Euclid could have proved more theorems that concern what we call vertical composition, he did not.\footnote{I know that this isn't quite what we usually mean by vertical composition. Since Euclid works in codiscrete categories, he cannot have three functors from one category to a second and natural transformations joining them pairwise.} I have selected these because they highlight the syntax, and can make some of the relations to category theory clear, and because the arguments are relatively short. The diagrams in many other proofs have a deep beauty that is missed in these proofs, but the arguments are also much longer. I turn therefore to what could be described as horizontal composition, though is also akin to a proof of functoriality.

\section{Proposition 22: Horizontal Composition}
In Proposition 22 of Book V, Euclid demonstrates that ``If there are any number of magnitudes whatever, and others equal to them in multitude, which taken two and two together are in the same ratio, then they are also in the same ratio ex aequali'' (David E. Joyce, trans. \citep{joyce-elements}) This is difficult to understand, but an annotated diagram, drawn using Euclid's labels and exemplary three pairs of magnitudes, makes it clear:

\centering%
\begin{tcolorbox}[%
	colbacktitle=white,
	coltitle = black,
	center title,
	halign upper=center,
	fonttitle = \bfseries,
	colbacklower = white,
	colback = white,
	width = 0.6\linewidth,
	lower separated=true,
	boxrule=.2mm]
\begin{tikzcd}[ampersand replacement = \&, math mode = false, row sep = 9ex]
{\begin{tikzpicture}
\draw[|-|] (0,0) -- (1.414,0);
\node[inner sep = 0, xshift = -3, anchor = east] at (0,0) {\tiny $A$};
\end{tikzpicture}}
\ar[r,
	shift right,
	gray,
	""' {name = AB, pos = 0.485}]
\ar[d,
	Leftrightarrow,
	densely dotted,
	gray!50,
	shift left,
	"\rotatebox{90}{\tiny Analogous}"',
	"\rotatebox{90}{\tiny Magnitudes}"]
\&%
	{\begin{tikzpicture}
	\draw[|-|] (0,0) -- (0.318,0);
	\node[inner sep = 0, xshift = -3, anchor = east] at (0,0) {\tiny $B$};
	\end{tikzpicture}}
	\ar[r,
		shift right,
		gray,
		""' {name = BC, pos = 0.507}]
	\ar[d,
		Leftrightarrow,
		densely dotted,
		gray!50,
		shift left,
		"\rotatebox{90}{\tiny Analogous}"',
		"\rotatebox{90}{\tiny Magnitudes}"]
	\&%
		{\begin{tikzpicture}
		\draw[|-|] (0,0) -- (1,0);
		\node[inner sep = 0, xshift = -3, anchor = east] at (0,0) {\tiny $C$};
		\end{tikzpicture}}
		\ar[d,
			Leftrightarrow,
			densely dotted,
			gray!50,
			shift left,
			"\rotatebox{90}{\tiny Analogous}"',
			"\rotatebox{90}{\tiny Magnitudes}"]\\
%
%
% Second Row
{\begin{tikzpicture}
\draw[|-|] (0,0) -- (0.8*1.414,0);
\node[inner sep = 0, xshift = -3, anchor = east] at (0,0) {\tiny $D$};
\end{tikzpicture}}
\ar[r,
	shift right,
	gray,
	"" {name = DE, pos = 0.515}]
\&%
	{\begin{tikzpicture}
	\draw[|-|] (0,0) -- (0.8*0.318,0);
	\node[inner sep = 0, xshift = -3, anchor = east] at (0,0) {\tiny $E$};
	\end{tikzpicture}}
	\ar[r,
		shift right,
		gray,
		"" {name = EF, pos = 0.493}]
	\&%
		{\begin{tikzpicture}
		\draw[|-|] (0,0) -- (0.8*1,0);
		\node[inner sep = 0, xshift = -3, anchor = east] at (0,0) {\tiny $F$};
		\end{tikzpicture}}
\ar[from = AB,
	to = DE,
	Leftrightarrow,
	gray!50,
	"\rotatebox{90}{\tiny Same}"',
	"\rotatebox{90}{\tiny Ratio}"]
\ar[from = BC,
	to = EF,
	Leftrightarrow,
	gray!50,
	"\rotatebox{90}{\tiny Same}"',
	"\rotatebox{90}{\tiny Ratio}"]
\end{tikzcd}

\tcblower
I say that,

\centering
\begin{tikzcd}[ampersand replacement = \&, math mode = false, row sep = 9ex]
{\begin{tikzpicture}
\draw[|-|] (0,0) -- (1.414,0);
\node[inner sep = 0, xshift = -3, anchor = east] at (0,0) {\tiny $A$};
\end{tikzpicture}}
\ar[r,
	shift right,
	gray!25,
	""' {name = AB, pos = 0.49}]
\ar[rr,
	bend right = 0.5cm,
	start anchor = {[xshift = -0.5ex]south east},
	end anchor = {[xshift = 0.5ex]south west},
	Plum!50!Orange,
	shift left,
	""' {name = AC, pos = 0.498}]
\ar[d,
	Leftrightarrow,
	densely dotted,
	gray!50,
	shift left,
	"\rotatebox{90}{\tiny Analogous}"',
	"\rotatebox{90}{\tiny Magnitudes}"]
\&%
	{\begin{tikzpicture}
	\draw[|-|, gray!50] (0,0) -- (0.318,0);
	\node[inner sep = 0, xshift = -3, anchor = east, gray!50] at (0,0) {\tiny $B$};
	\end{tikzpicture}}
	\ar[r,
		shift right,
		gray!50,
		""' {name = BC, pos = 0.505}]
	\&%
		{\begin{tikzpicture}
		\draw[|-|] (0,0) -- (1,0);
		\node[inner sep = 0, xshift = -3, anchor = east] at (0,0) {\tiny $C$};
		\end{tikzpicture}}
		\ar[d,
			Leftrightarrow,
			densely dotted,
			gray!50,
			shift left,
			"\rotatebox{90}{\tiny Analogous}"',
			"\rotatebox{90}{\tiny Magnitudes}"]\\
%
%
% Second Row
{\begin{tikzpicture}
\draw[|-|] (0,0) -- (0.8*1.414,0);
\node[inner sep = 0, xshift = -3, anchor = east] at (0,0) {\tiny $D$};
\end{tikzpicture}}
\ar[r,
	shift right,
	gray!50,
	"" {name = DE, pos = 0.51}]
\ar[rr,
	bend left = 0.5cm,
	start anchor = {[xshift = -0.5ex]north east},
	end anchor = {[xshift = 0.5ex]north west},
	shift right,
	Plum!50!Orange,
	"" {name = DF, pos = 0.502}]
\&%
	{\begin{tikzpicture}
	\draw[|-|, gray!50] (0,0) -- (0.8*0.318,0);
	\node[inner sep = 0, xshift = -3, anchor = east, gray!50] at (0,0) {\tiny $E$};
	\end{tikzpicture}}
	\ar[r,
		shift right,
		gray!50,
		"" {name = EF, pos = 0.495}]
	\&%
		{\begin{tikzpicture}
		\draw[|-|] (0,0) -- (0.8*1,0);
		\node[inner sep = 0, xshift = -3, anchor = east] at (0,0) {\tiny $F$};
		\end{tikzpicture}}
\ar[from = AC,
	to = DF,
	Leftrightarrow,
	Plum!50!Orange!75,
	"\rotatebox{90}{\tiny Same}"',
	"\rotatebox{90}{\tiny Ratio}"]
\end{tikzcd}
\end{tcolorbox}

\flushleftright%
I should add that though this is perhaps only somewhat like horizontal composition (since the ratios aren't parallel), the composition of the ratios on the top and bottom rows absolutely \textit{is} composition of arrows. Just as Eugenia Cheng's diagrams of biological relations, drawn above, show actual relational composition, so this diagram shows actual relational composition (though, to keep the diagram from being cluttered I have grayed out the composed relations in the bottom diagram).\footnote{As we will see below, mathematicians from Euclid up to and including Newton saw this relational composition as a form of addition.}

Euclid argues by taking equimultiples of all the quantities. It is impossible to draw them all horizontally (since the equimultiple for $E$ would be on top of $D$ or $F$) so I have drawn the equimultiples vertically.

\centering%
\begin{tcolorbox}[%
	colbacktitle=white,
	coltitle = black,
	center title,
	halign upper=center,
	fonttitle = \bfseries,
	colbacklower = white,
	colback = white,
	width = 0.8\linewidth,
	lower separated=true,
	boxrule=.2mm]
\begin{tikzcd}[ampersand replacement = \&, math mode = false, row sep = 8ex]
{} \& {} \& {}\\[-2ex]
{\begin{tikzpicture}
\draw[|-|, black!75] (0,0) -- (2*1.414,0);
\node[inner sep = 0, xshift = -3, anchor = east, black!75] at (0,0) {\tiny $G$};
\end{tikzpicture}}
\ar[d,
	OliveGreen!50,
	shift left,
	"\rotatebox{-90}{\tiny Multiple of}"]
\ar[u,
	Leftarrow,
	densely dotted,
	shift right,
	OliveGreen!25,
	"\rotatebox{90}{\tiny Equi-}",
	"\rotatebox{90}{\tiny of}"']
\&%
	{\begin{tikzpicture}
	\draw[|-|, black!75] (0,0) -- (4*0.318,0);
	\node[inner sep = 0, xshift = -3, anchor = east, black!75] at (0,0) {\tiny $K$};
	\end{tikzpicture}}
	\ar[d,
		Blue!50,
		shift left,
		"\rotatebox{-90}{\tiny Multiple of}"]
	\ar[u,
		Leftarrow,
		densely dotted,
		shift right,
		Blue!25,
		"\rotatebox{90}{\tiny Equi-}",
		"\rotatebox{90}{\tiny of}"']
	\&%
		{\begin{tikzpicture}
		\draw[|-|, black!75] (0,0) -- (3*1,0);
		\node[inner sep = 0, xshift = -3, anchor = east, black!75] at (0,0) {\tiny $M$};
		\end{tikzpicture}}
		\ar[d,
			Sepia!50,
			shift left,
			"\rotatebox{-90}{\tiny Multiple of}"]
		\ar[u,
			Leftarrow,
			densely dotted,
			shift right,
			Sepia!25,
			"\rotatebox{90}{\tiny Equi-}",
			"\rotatebox{90}{\tiny of}"']\\
%
%
% Second Row
{\begin{tikzpicture}
\draw[|-|] (0,0) -- (1.414,0);
\node[inner sep = 0, xshift = -3, anchor = east] at (0,0) {\tiny $A$};
\end{tikzpicture}}
\ar[r,
	shift right,
	gray,
	""' {name = AB, pos = 0.49}]
\ar[d,
	Leftrightarrow,
	densely dotted,
	gray!50,
	shift left]
\&%
	{\begin{tikzpicture}
	\draw[|-|] (0,0) -- (0.318,0);
	\node[inner sep = 0, xshift = -3, anchor = east] at (0,0) {\tiny $B$};
	\end{tikzpicture}}
	\ar[r,
		shift right,
		gray,
		""' {name = BC, pos = 0.505}]
	\ar[d,
		Leftrightarrow,
		densely dotted,
		gray!50,
		shift left]
	\&%
		{\begin{tikzpicture}
		\draw[|-|] (0,0) -- (1,0);
		\node[inner sep = 0, xshift = -3, anchor = east] at (0,0) {\tiny $C$};
		\end{tikzpicture}}
		\ar[d,
			Leftrightarrow,
			densely dotted,
			gray!50,
			shift left]\\[-3ex]
%
%
% Third Row
{\begin{tikzpicture}
\draw[|-|] (0,0) -- (0.8*1.414,0);
\node[inner sep = 0, xshift = -3, anchor = east] at (0,0) {\tiny $D$};
\end{tikzpicture}}
\ar[r,
	shift right,
	gray,
	"" {name = DE, pos = 0.51}]
\ar[d,
	leftarrow,
	shift left,
	OliveGreen!50,
	"\rotatebox{90}{\tiny Same}"',
	"\rotatebox{90}{\tiny Multiple of}"]
\&%
	{\begin{tikzpicture}
	\draw[|-|] (0,0) -- (0.8*0.318,0);
	\node[inner sep = 0, xshift = -3, anchor = east] at (0,0) {\tiny $E$};
	\end{tikzpicture}}
	\ar[r,
		shift right,
		gray,
		"" {name = EF, pos = 0.495}]
	\ar[d,
		leftarrow,
		shift left,
		Blue!50,
		"\rotatebox{90}{\tiny Same}"',
		"\rotatebox{90}{\tiny Multiple of}"]
	\&%
		{\begin{tikzpicture}
		\draw[|-|] (0,0) -- (0.8*1,0);
		\node[inner sep = 0, xshift = -3, anchor = east] at (0,0) {\tiny $F$};
		\end{tikzpicture}}
		\ar[d,
			leftarrow,
			shift left,
			Sepia!50,
			"\rotatebox{90}{\tiny Same}"',
			"\rotatebox{90}{\tiny Multiple of}"]\\
%
%
% Fourth Row
{\begin{tikzpicture}
\draw[|-|, black!75] (0,0) -- (2*0.8*1.414,0);
\node[inner sep = 0, xshift = -3, anchor = east, black!75] at (0,0) {\tiny $H$};
\end{tikzpicture}}
\ar[d,
	Leftarrow,
	OliveGreen!25,
	densely dotted,
	shift left,
	"\rotatebox{90}{\tiny Multiples}"',
	"\rotatebox{90}{\scalebox{0.95}{\tiny $D;$ $A$}}"]
\&%
	{\begin{tikzpicture}
	\draw[|-|, black!75] (0,0) -- (4*0.8*0.318,0);
	\node[inner sep = 0, xshift = -3, anchor = east, black!75] at (0,0) {\tiny $L$};
	\end{tikzpicture}}
	\ar[d,
		Leftarrow,
		Blue!25,
		densely dotted,
		shift left,
		"\rotatebox{90}{\tiny Multiples}"',
		"\rotatebox{90}{\scalebox{0.95}{\tiny $E;$ $B$}}"]
	\&%
		{\begin{tikzpicture}
		\draw[|-|, black!75] (0,0) -- (3*0.8*1,0);
		\node[inner sep = 0, xshift = -3, anchor = east, black!75] at (0,0) {\tiny $N$};
		\end{tikzpicture}}
		\ar[d,
			Leftarrow,
			Sepia!25,
			densely dotted,
			shift left,
			"\rotatebox{90}{\tiny Multiples}"',
			"\rotatebox{90}{\scalebox{0.95}{\tiny $C;$ $F$}}"]\\[-2ex]
{} \& {} \& {}
\ar[from = AB,
	to = DE,
	Leftrightarrow,
	gray!50,
	"\rotatebox{90}{\tiny Same}"',
	"\rotatebox{90}{\tiny Ratio}"]
\ar[from = BC,
	to = EF,
	Leftrightarrow,
	gray!50,
	"\rotatebox{90}{\tiny Same}"',
	"\rotatebox{90}{\tiny Ratio}"]
\end{tikzcd}
\end{tcolorbox}

\flushleftright
Though the equimultiples are labeled, as on a cylinder, the same multiples are color-coded.

Since equimultiples of equimultiples are equimultiples of the original quantities (V.3), it follows (V.4) that if two pairs of quantities have the same ratio, equimultiples of the antecedents (the domain of the original ratio) have the same ratio to the equimultiples of the consequents (the codomain of the original ratios), though, in general, not the same ratio as the original quantities have. Thus, $G$ has the same ratio to $K$ that $H$ has to $L$, and $K$ has the same ratio to $M$ that $L$ has to $N$.

\centering%
\begin{tcolorbox}[%
	colbacktitle=white,
	coltitle = black,
	center title,
	halign upper=center,
	fonttitle = \bfseries,
	colbacklower = white,
	colback = white,
	width = 0.8\linewidth,
	lower separated=true,
	boxrule=.2mm]
\begin{tikzcd}[ampersand replacement = \&, math mode = false, row sep = 8ex]
{} \& {} \& {}\\[-2ex]
{\begin{tikzpicture}
\draw[|-|] (0,0) -- (2*1.414,0);
\node[inner sep = 0, xshift = -3, anchor = east] at (0,0) {\tiny $G$};
\end{tikzpicture}}
\ar[r,
	shift right,
	gray,
	""' {name = GK, pos = 0.49}]
\ar[d,
	OliveGreen!50,
	shift left,
	"\rotatebox{-90}{\tiny Multiple of}"]
\ar[u,
	Leftarrow,
	densely dotted,
	shift right,
	OliveGreen!25,
	"\rotatebox{90}{\tiny Equi-}",
	"\rotatebox{90}{\tiny of}"']
\&%
	{\begin{tikzpicture}
	\draw[|-|] (0,0) -- (4*0.318,0);
	\node[inner sep = 0, xshift = -3, anchor = east] at (0,0) {\tiny $K$};
	\end{tikzpicture}}
	\ar[r,
		shift right,
		gray,
		""' {name = KM, pos = 0.515}]
	\ar[d,
		Blue!50,
		shift left,
		"\rotatebox{-90}{\tiny Multiple of}"]
	\ar[u,
		Leftarrow,
		densely dotted,
		shift right,
		Blue!25,
		"\rotatebox{90}{\tiny Equi-}",
		"\rotatebox{90}{\tiny of}"']
	\&%
		{\begin{tikzpicture}
		\draw[|-|] (0,0) -- (3*1,0);
		\node[inner sep = 0, xshift = -3, anchor = east] at (0,0) {\tiny $M$};
		\end{tikzpicture}}
		\ar[d,
			Sepia!50,
			shift left,
			"\rotatebox{-90}{\tiny Multiple of}"]
		\ar[u,
			Leftarrow,
			densely dotted,
			shift right,
			Sepia!25,
			"\rotatebox{90}{\tiny Equi-}",
			"\rotatebox{90}{\tiny of}"']\\
%
%
% Second Row
{\begin{tikzpicture}
\draw[|-|, gray] (0,0) -- (1.414,0);
\node[inner sep = 0, xshift = -3, anchor = east, gray] at (0,0) {\tiny $A$};
\end{tikzpicture}}
\ar[d,
	Leftrightarrow,
	densely dotted,
	gray!50,
	shift left]
\ar[r,
	shift right,
	gray!25,
	""' {name = AB, pos = 0.49}]
\&%
	{\begin{tikzpicture}
	\draw[|-|, gray] (0,0) -- (0.318,0);
	\node[inner sep = 0, xshift = -3, anchor = east, gray] at (0,0) {\tiny $B$};
	\end{tikzpicture}}
	\ar[r,
		shift right,
		gray!25,
		""' {name = BC, pos = 0.503}]
	\ar[d,
		Leftrightarrow,
		densely dotted,
		gray!25,
		shift left]
	\&%
		{\begin{tikzpicture}
		\draw[|-|, gray] (0,0) -- (1,0);
		\node[inner sep = 0, xshift = -3, anchor = east, gray] at (0,0) {\tiny $C$};
		\end{tikzpicture}}
		\ar[d,
			Leftrightarrow,
			densely dotted,
			gray!50,
			shift left]\\[-3ex]
%
%
% Third Row
{\begin{tikzpicture}
\draw[|-|, gray] (0,0) -- (0.8*1.414,0);
\node[inner sep = 0, xshift = -3, anchor = east, gray] at (0,0) {\tiny $D$};
\end{tikzpicture}}
\ar[r,
	shift right,
	gray!25,
	"" {name = DE, pos = 0.51}]
\ar[d,
	leftarrow,
	shift left,
	OliveGreen!50,
	"\rotatebox{90}{\tiny Same}"',
	"\rotatebox{90}{\tiny Multiple of}"]
\&%
	{\begin{tikzpicture}
	\draw[|-|, gray] (0,0) -- (0.8*0.318,0);
	\node[inner sep = 0, xshift = -3, anchor = east, gray] at (0,0) {\tiny $E$};
	\end{tikzpicture}}
	\ar[r,
		shift right,
		gray!25,
		"" {name = EF, pos = 0.495}]
	\ar[d,
		leftarrow,
		shift left,
		Blue!50,
		"\rotatebox{90}{\tiny Same}"',
		"\rotatebox{90}{\tiny Multiple of}"]
	\&%
		{\begin{tikzpicture}
		\draw[|-|, gray] (0,0) -- (0.8*1,0);
		\node[inner sep = 0, xshift = -3, anchor = east, gray] at (0,0) {\tiny $F$};
		\end{tikzpicture}}
		\ar[d,
			leftarrow,
			shift left,
			Sepia!50,
			"\rotatebox{90}{\tiny Same}"',
			"\rotatebox{90}{\tiny Multiple of}"]\\
%
%
% Fourth Row
{\begin{tikzpicture}
\draw[|-|] (0,0) -- (2*0.8*1.414,0);
\node[inner sep = 0, xshift = -3, anchor = east] at (0,0) {\tiny $H$};
\end{tikzpicture}}
\ar[r,
	shift right,
	gray,
	"" {name = HL, pos = 0.51}]
\ar[d,
	Leftarrow,
	OliveGreen!25,
	densely dotted,
	shift left,
	"\rotatebox{90}{\tiny Multiples}"',
	"\rotatebox{90}{\scalebox{0.95}{\tiny $D;$ $A$}}"]
\&%
	{\begin{tikzpicture}
	\draw[|-|] (0,0) -- (4*0.8*0.318,0);
	\node[inner sep = 0, xshift = -3, anchor = east] at (0,0) {\tiny $L$};
	\end{tikzpicture}}
	\ar[r,
		shift right,
		gray,
		"" {name = LN, pos = 0.47}]
	\ar[d,
		Leftarrow,
		densely dotted,
		Blue!25,
		shift left,
		"\rotatebox{90}{\tiny Multiples}"',
		"\rotatebox{90}{\scalebox{0.95}{\tiny $E;$ $B$}}"]
	\&%
		{\begin{tikzpicture}
		\draw[|-|] (0,0) -- (3*0.8*1,0);
		\node[inner sep = 0, xshift = -3, anchor = east] at (0,0) {\tiny $N$};
		\end{tikzpicture}}
		\ar[d,
			Leftarrow,
			densely dotted,
			Sepia!25,
			shift left,
			"\rotatebox{90}{\tiny Multiples}"',
			"\rotatebox{90}{\scalebox{0.95}{\tiny $C;$ $F$}}"]\\[-2ex]
{} \& {} \& {}
\ar[from = GK,
	to = HL,
	Leftrightarrow,
	crossing over,
	gray,
	"\rotatebox{90}{Same}"',
	"\rotatebox{90}{Ratio}"]
\ar[from = KM,
	to = LN,
	Leftrightarrow,
	crossing over,
	gray,
	"\rotatebox{90}{Same}"',
	"\rotatebox{90}{Ratio}"]
\end{tikzcd}
\end{tcolorbox}

\flushleftright%
But now we have two pairs of three multitudes and the first, in each pair, have the same ratio to the second; and the second have the same ratio to the third. Therefore (V.20), if the first in one pair of three magnitudes is greater than the third in that pair, then the first in the other pair is also greater than the third in that pair---and if the one is equal the other is equal, and if it is lesser, the other is lesser.

\centering%
\begin{tcolorbox}[%
	colbacktitle=white,
	coltitle = black,
	center title,
	halign upper=center,
	fonttitle = \bfseries,
	colbacklower = white,
	colback = white,
	width = 0.8\linewidth,
	lower separated=true,
	boxrule=.2mm]
\begin{tikzcd}[ampersand replacement = \&, math mode = false, row sep = 8ex]
{} \& {} \& {}\\[-2ex]
{\begin{tikzpicture}
\draw[|-|] (0,0) -- (2*1.414,0);
\node[inner sep = 0, xshift = -3, anchor = east] at (0,0) {\tiny $G$};
\end{tikzpicture}}
\ar[r,
	shift right,
	gray!50,
	""' {name = GK, pos = 0.49}]
\ar[d,
	OliveGreen!25,
	shift left,
	"\rotatebox{-90}{\tiny Multiple of}"]
\ar[u,
	Leftarrow,
	densely dotted,
	shift right,
	OliveGreen!12.5,
	"\rotatebox{90}{\tiny Equi-}",
	"\rotatebox{90}{\tiny of}"']
\&%
	{\begin{tikzpicture}
	\draw[|-|, gray] (0,0) -- (4*0.318,0);
	\node[inner sep = 0, xshift = -3, anchor = east, gray] at (0,0) {\tiny $K$};
	\end{tikzpicture}}
	\ar[r,
		shift right,
		gray!50,
		""' {name = KM, pos = 0.515}]
	\ar[d,
		Blue!25,
		shift left,
		"\rotatebox{-90}{\tiny Multiple of}"]
	\&%
		{\begin{tikzpicture}
		\draw[|-|] (0,0) -- (3*1,0);
		\node[inner sep = 0, xshift = -3, anchor = east] at (0,0) {\tiny $M$};
		\end{tikzpicture}}
		\ar[d,
			Sepia!25,
			shift left,
			"\rotatebox{-90}{\tiny Multiple of}"]
		\ar[u,
			Leftarrow,
			densely dotted,
			shift right,
			Sepia!12.5,
			"\rotatebox{90}{\tiny Equi-}",
			"\rotatebox{90}{\tiny of}"']\\
%
%
% Second Row
{\begin{tikzpicture}
\draw[|-|, gray!50] (0,0) -- (1.414,0);
\node[inner sep = 0, xshift = -3, anchor = east, gray!50] at (0,0) {\tiny $A$};
\end{tikzpicture}}
\ar[d,
	Leftrightarrow,
	densely dotted,
	gray!50,
	shift left]
\ar[r,
	shift right,
	gray!25,
	""' {name = AB, pos = 0.49}]
\&%
	{\begin{tikzpicture}
	\draw[|-|, gray!50] (0,0) -- (0.318,0);
	\node[inner sep = 0, xshift = -3, anchor = east, gray!50] at (0,0) {\tiny $B$};
	\end{tikzpicture}}
	\ar[r,
		shift right,
		gray!25,
		""' {name = BC, pos = 0.503}]
	\ar[d,
		Leftrightarrow,
		densely dotted,
		gray!25,
		shift left]
	\&%
		{\begin{tikzpicture}
		\draw[|-|, gray!50] (0,0) -- (1,0);
		\node[inner sep = 0, xshift = -3, anchor = east, gray!50] at (0,0) {\tiny $C$};
		\end{tikzpicture}}
		\ar[d,
			Leftrightarrow,
			densely dotted,
			gray!50,
			shift left]\\[-3ex]
%
%
% Third Row
{\begin{tikzpicture}
\draw[|-|, gray!50] (0,0) -- (0.8*1.414,0);
\node[inner sep = 0, xshift = -3, anchor = east, gray!50] at (0,0) {\tiny $D$};
\end{tikzpicture}}
\ar[r,
	shift right,
	gray!25,
	"" {name = DE, pos = 0.51}]
\ar[d,
	leftarrow,
	shift left,
	OliveGreen!25,
	"\rotatebox{90}{\tiny Same}"',
	"\rotatebox{90}{\tiny Multiple of}"]
\&%
	{\begin{tikzpicture}
	\draw[|-|, gray!50] (0,0) -- (0.8*0.318,0);
	\node[inner sep = 0, xshift = -3, anchor = east, gray!50] at (0,0) {\tiny $E$};
	\end{tikzpicture}}
	\ar[r,
		shift right,
		gray!25,
		"" {name = EF, pos = 0.495}]
	\ar[d,
		leftarrow,
		shift left,
		Blue!25,
		"\rotatebox{90}{\tiny Same}"',
		"\rotatebox{90}{\tiny Multiple of}"]
	\&%
		{\begin{tikzpicture}
		\draw[|-|, gray!50] (0,0) -- (0.8*1,0);
		\node[inner sep = 0, xshift = -3, anchor = east, gray!50] at (0,0) {\tiny $F$};
		\end{tikzpicture}}
		\ar[d,
			leftarrow,
			shift left,
			Sepia!25,
			"\rotatebox{90}{\tiny Same}"',
			"\rotatebox{90}{\tiny Multiple of}"]\\
%
%
% Fourth Row
{\begin{tikzpicture}
\draw[|-|] (0,0) -- (2*0.8*1.414,0);
\node[inner sep = 0, xshift = -3, anchor = east] at (0,0) {\tiny $H$};
\end{tikzpicture}}
\ar[r,
	shift right,
	gray!50,
	"" {name = HL, pos = 0.51}]
\ar[d,
	Leftarrow,
	densely dotted,
	OliveGreen!12.5,
	shift left,
	"\rotatebox{90}{\tiny Multiples}"',
	"\rotatebox{90}{\scalebox{0.95}{\tiny $D;$ $A$}}"]
\&%
	{\begin{tikzpicture}
	\draw[|-|, gray] (0,0) -- (4*0.8*0.318,0);
	\node[inner sep = 0, xshift = -3, anchor = east, gray] at (0,0) {\tiny $L$};
	\end{tikzpicture}}
	\ar[r,
		shift right,
		gray!50,
		"" {name = LN, pos = 0.47}]
	\&%
		{\begin{tikzpicture}
		\draw[|-|] (0,0) -- (3*0.8*1,0);
		\node[inner sep = 0, xshift = -3, anchor = east] at (0,0) {\tiny $N$};
		\end{tikzpicture}}
		\ar[d,
			Leftarrow,
			densely dotted,
			Sepia!12.5,
			shift left,
			"\rotatebox{90}{\tiny Multiples}"',
			"\rotatebox{90}{\scalebox{0.95}{\tiny $C;$ $F$}}"]\\[-2ex]
{} \& {} \& {}
\ar[from = GK,
	to = HL,
	Leftrightarrow,
	crossing over,
	gray!50,
	"\rotatebox{90}{Same}"',
	"\rotatebox{90}{Ratio}"]
\ar[from = KM,
	to = LN,
	Leftrightarrow,
	crossing over,
	gray!50,
	"\rotatebox{90}{Same}"',
	"\rotatebox{90}{Ratio}"]
\ar[from = 2-1,
	to = 2-3,
	start anchor = {[xshift = -1.25ex]north east},
	end anchor = {[xshift = 1.25ex]north west},
	blue,
	bend left = 1.6cm,
	"\scriptsize Exceeds"]
\ar[from = 2-1,
	to = 2-3,
	start anchor = {[xshift = -0.75ex, yshift = -0.25ex]north east},
	end anchor = {[xshift = 0.75ex, yshift = -0.25ex]north west},
	ForestGreen,
	bend left = 0.9cm,
	"\scriptsize Equals"]
\ar[from = 2-1,
	to = 2-3,
	start anchor = {[xshift = -0.25ex, yshift = -0.5ex]north east},
	end anchor = {[xshift = 0.25ex, yshift = -0.5ex]north west},
	Orange,
	bend left = 0.3cm,
	"\scriptsize Less than"]
% Bottom arrows
\ar[from = 5-1,
	to = 5-3,
	start anchor = {[xshift = -1.25ex]south east},
	end anchor = {[xshift = 1.25ex]south west},
	Orange,
	bend right = 1.45cm,
	"\scriptsize Less than"]
\ar[from = 5-1,
	to = 5-3,
	start anchor = {[xshift = -0.75ex, yshift = 0.25ex]south east},
	end anchor = {[xshift = 0.75ex, yshift = 0.25ex]south west},
	ForestGreen,
	bend right = 1cm,
	"\scriptsize Equals"]
\ar[from = 5-1,
	to = 5-3,
	start anchor = {[xshift = -0.25ex, yshift = 0.5ex]south east},
	end anchor = {[xshift = 0.25ex, yshift = 0.5ex]south west},
	blue,
	bend right = 0.4cm,
	"\scriptsize Exceeds"]
\end{tikzcd}
\end{tcolorbox}

\flushleftright%
And therefore---as we can plainly see---$A$ has the same ratio to $C$ as $D$ has to $F$.

\centering%
\begin{tcolorbox}[%
	colbacktitle=white,
	coltitle = black,
	center title,
	halign upper=center,
	fonttitle = \bfseries,
	colbacklower = white,
	colback = white,
	width = 0.8\linewidth,
	lower separated=true,
	boxrule=.2mm]
\begin{tikzcd}[ampersand replacement = \&, math mode = false, row sep = 6ex]
{} \& {} \& {}\\[-2ex]
{\begin{tikzpicture}
\draw[|-|, gray!50] (0,0) -- (2*1.414,0);
\node[inner sep = 0, xshift = -3, anchor = east, gray!50] at (0,0) {\tiny $G$};
\end{tikzpicture}}
\ar[r,
	shift right,
	gray!15,
	""' {name = GK, pos = 0.49}]
\ar[d,
	OliveGreen!33,
	shift left,
	"\rotatebox{-90}{\scalebox{0.8}{\tiny Multiple of}}"]
\ar[u,
	Leftarrow,
	densely dotted,
	shift right,
	OliveGreen!33,
	"\rotatebox{90}{\scalebox{0.8}{\tiny Equi-}}",
	"\rotatebox{90}{\scalebox{0.8}{\tiny of}}"']
\&%
	{\begin{tikzpicture}
	\draw[|-|, gray!25] (0,0) -- (4*0.318,0);
	\node[inner sep = 0, xshift = -3, anchor = east, gray!25] at (0,0) {\tiny $K$};
	\end{tikzpicture}}
	\ar[r,
		shift right,
		gray!15,
		""' {name = KM, pos = 0.515}]
	\ar[d,
		Blue!12.5,
		shift left,
		"\rotatebox{-90}{\scalebox{0.8}{\tiny Multiple of}}"]
	\&%
		{\begin{tikzpicture}
		\draw[|-|, gray!50] (0,0) -- (3*1,0);
		\node[inner sep = 0, xshift = -3, anchor = east, gray!50] at (0,0) {\tiny $M$};
		\end{tikzpicture}}
		\ar[d,
			Sepia!25,
			shift left,
			"\rotatebox{-90}{\scalebox{0.8}{\tiny Multiple of}}"]
		\ar[u,
			Leftarrow,
			densely dotted,
			shift right,
			Sepia!25,
			"\rotatebox{90}{\scalebox{0.8}{\tiny Equi-}}",
			"\rotatebox{90}{\scalebox{0.8}{\tiny of}}"']\\
%
%
% Second Row
{\begin{tikzpicture}
\draw[|-|] (0,0) -- (1.414,0);
\node[inner sep = 0, xshift = -3, anchor = east] at (0,0) {\tiny $A$};
\end{tikzpicture}}
\ar[rr,
	bend right = 0.25cm,
	start anchor = {[xshift = -0.5ex, yshift = 0.5ex]south east},
	end anchor = {[xshift = 0.5ex, yshift = 0.5ex]south west},
	""' {name = AC}]
\ar[r,
	shift right,
	gray!50,
	""' {name = AB, pos = 0.49}]
\ar[d,
	Leftrightarrow,
	densely dotted,
	gray,
	shift left,
	"\rotatebox{90}{\tiny Analogous}"',
	"\rotatebox{90}{\tiny Magnitudes}"]
\&%
	{\begin{tikzpicture}
	\draw[|-|, gray] (0,0) -- (0.318,0);
	\node[inner sep = 0, xshift = -3, anchor = east, gray] at (0,0) {\tiny $B$};
	\end{tikzpicture}}
	\ar[r,
		shift right,
		gray!25,
		""' {name = BC, pos = 0.503}]
	\&%
		{\begin{tikzpicture}
		\draw[|-|] (0,0) -- (1,0);
		\node[inner sep = 0, xshift = -3, anchor = east] at (0,0) {\tiny $C$};
		\end{tikzpicture}}
		\ar[d,
			Leftrightarrow,
			densely dotted,
			gray,
			shift left,
			"\rotatebox{90}{\tiny Analogous}"',
			"\rotatebox{90}{\tiny Magnitudes}"]\\[3ex]
%
%
% Third Row
{\begin{tikzpicture}
\draw[|-|] (0,0) -- (0.8*1.414,0);
\node[inner sep = 0, xshift = -3, anchor = east] at (0,0) {\tiny $D$};
\end{tikzpicture}}
\ar[rr,
	bend left = 0.25cm,
	start anchor = {[xshift = -0.25ex, yshift = -0.5ex]north east},
	end anchor = {[xshift = 0.25ex, yshift = -0.5ex]north west},
	"" {name = DF}]
\ar[r,
	shift right,
	gray!50,
	"" {name = DE, pos = 0.51}]
\ar[d,
	leftarrow,
	shift left,
	OliveGreen!33,
	"\rotatebox{90}{\scalebox{0.8}{\tiny Same}}"',
	"\rotatebox{90}{\scalebox{0.8}{\tiny Multiple of}}"]
\&%
	{\begin{tikzpicture}
	\draw[|-|, gray] (0,0) -- (0.8*0.318,0);
	\node[inner sep = 0, xshift = -3, anchor = east, gray!50] at (0,0) {\tiny $E$};
	\end{tikzpicture}}
	\ar[r,
		shift right,
		gray!50,
		"" {name = EF, pos = 0.495}]
	\ar[d,
		leftarrow,
		shift left,
		Blue!12.5,
		"\rotatebox{90}{\scalebox{0.8}{\tiny Same}}"',
		"\rotatebox{90}{\scalebox{0.8}{\tiny Multiple of}}"]
	\&%
		{\begin{tikzpicture}
		\draw[|-|] (0,0) -- (0.8*1,0);
		\node[inner sep = 0, xshift = -3, anchor = east] at (0,0) {\tiny $F$};
		\end{tikzpicture}}
		\ar[d,
			leftarrow,
			shift left,
			Sepia!25,
			"\rotatebox{90}{\scalebox{0.8}{\tiny Same}}"',
			"\rotatebox{90}{\scalebox{0.8}{\tiny Multiple of}}"]\\
%
%
% Fourth Row
{\begin{tikzpicture}
\draw[|-|, gray!50] (0,0) -- (2*0.8*1.414,0);
\node[inner sep = 0, xshift = -3, anchor = east, gray!50] at (0,0) {\tiny $H$};
\end{tikzpicture}}
\ar[r,
	shift right,
	gray!15,
	"" {name = HL, pos = 0.51}]
\ar[d,
	Leftarrow,
	densely dotted,
	OliveGreen!33,
	shift left,
	"\rotatebox{90}{\scalebox{0.8}{\tiny Multiples}}"',
	"\rotatebox{90}{\scalebox{0.76}{\tiny $D;$ $A$}}"]
\&%
	{\begin{tikzpicture}
	\draw[|-|, gray!25] (0,0) -- (4*0.8*0.318,0);
	\node[inner sep = 0, xshift = -3, anchor = east, gray!25] at (0,0) {\tiny $L$};
	\end{tikzpicture}}
	\ar[r,
		shift right,
		gray!15,
		"" {name = LN, pos = 0.47}]
	\&%
		{\begin{tikzpicture}
		\draw[|-|, gray!50] (0,0) -- (3*0.8*1,0);
		\node[inner sep = 0, xshift = -3, anchor = east, gray!50] at (0,0) {\tiny $N$};
		\end{tikzpicture}}
		\ar[d,
			Leftarrow,
			densely dotted,
			Sepia!25,
			shift left,
			"\rotatebox{90}{\scalebox{0.8}{\tiny Multiples}}"',
			"\rotatebox{90}{\scalebox{0.76}{\tiny $C;$ $F$}}"]\\[-2ex]
{} \& {} \& {}
\ar[from = 2-1,
	to = 2-3,
	start anchor = {[xshift = -1.25ex]north east},
	end anchor = {[xshift = 1.25ex]north west},
	blue!17,
	bend left = 1.6cm,
	"\scriptsize Exceeds"]
\ar[from = 2-1,
	to = 2-3,
	start anchor = {[xshift = -0.75ex, yshift = -0.25ex]north east},
	end anchor = {[xshift = 0.75ex, yshift = -0.25ex]north west},
	ForestGreen!17,
	bend left = 0.9cm,
	"\scriptsize Equals"]
\ar[from = 2-1,
	to = 2-3,
	start anchor = {[xshift = -0.25ex, yshift = -0.5ex]north east},
	end anchor = {[xshift = 0.25ex, yshift = -0.5ex]north west},
	Orange!17,
	bend left = 0.3cm,
	"\scriptsize Less than"]
% Bottom arrows
\ar[from = 5-1,
	to = 5-3,
	start anchor = {[xshift = -1.25ex]south east},
	end anchor = {[xshift = 1.25ex]south west},
	Orange!17,
	bend right = 1.45cm,
	"\scriptsize Less than"]
\ar[from = 5-1,
	to = 5-3,
	start anchor = {[xshift = -0.75ex, yshift = 0.25ex]south east},
	end anchor = {[xshift = 0.75ex, yshift = 0.25ex]south west},
	ForestGreen!17,
	bend right = 1cm,
	"\scriptsize Equals"]
\ar[from = 5-1,
	to = 5-3,
	start anchor = {[xshift = -0.25ex, yshift = 0.5ex]south east},
	end anchor = {[xshift = 0.25ex, yshift = 0.5ex]south west},
	blue!17,
	bend right = 0.4cm,
	"\scriptsize Exceeds"]
\ar[from = AC,
	to = DF,
	gray,
	Leftrightarrow,
	"\rotatebox{90}{\tiny Same}"',
	"\rotatebox{90}{\tiny Ratio}"]
\end{tikzcd}
\end{tcolorbox}

\flushleftright%
As noted, similar annotations can be added immediately to every diagram throughout Book V, and they can illustrate every step in every argument. We could go through a number of these, but perhaps it is better to show how these annotations can be added to a particularly tricky proposition from Book VI.

\section{Proposition VI.19: Duplicate Ratios and Similar Triangles}

In this proposition, Euclid proves that similar triangles are in the duplicate ratio of their bases.

\subsection{Compounded Ratios}
Before commenting on Euclid's proposition, we should pause to consider the sense of duplicate ratios. Writing about allegories (categories of relations) Freyd and Schedrov \citep{freyd-scedrov1990} define the binary operation of ``COMPOSITION, RS ($x$ RS $y$ iff there exists a $z$ such that $x$ R $z$ and $z$ S $y$)'' (p. 195), where $x$ R $y$ is the notation for $x$ has relation R to $y$. Using this notation, when the ratios R and S are the same ratio, so that $x$ R $z$ and $z$ R $y$, then $x$ RR $y$. In this case, Euclid says that the ratio R has been duplicated, or rather, that $x$ has, to $z$ the duplicate ratio to the ratio it has to $y$.

Thus, for example, nine has the same ratio to six as six has to four---namely, in least terms, the ratio that three has to two. Therefore we say that nine has, to four, a duplicate ratio to that three has to two. We can symbolize this with a diagram:

\centering%
\begin{tikzcd}[column sep = 1cm, ampersand replacement = \&, row sep = 1.5cm]
3
\ar[r,
	""' {name = 65A}]
\ar[d,
	Leftrightarrow,
	gray,
	densely dotted,
	"\text{analogous}" sloped,
	"\text{quantities}"' sloped]
\&
	2
	\ar[dr,
		Leftrightarrow,
		gray,
		densely dotted]
	\&[-1cm]
	,
	\&[-1cm]
		3
		\ar[dl,
			Leftrightarrow,
			gray,
			densely dotted]
		\ar[r, ""' {name = 65B}]
		\&
			2\ar[d,
				Leftrightarrow,
				gray,
				densely dotted,
				"\text{analogous}" sloped,
				"\text{quantities}"' sloped]\\
9
\ar[rr, "" {name = 3630}]
\ar[rrrr,
	bend right,
	"{$(3\to 2)\;\circ\; (3\to 2),$}\\\footnotesize as it were"' {align = center, font = \tiny}]
\&\&
	6
	\ar[rr, "" {name = 3025}]
	\&\&
	4
\ar[from = 65B,
	to = 3025,
	Leftrightarrow,
	gray,
	"\text{\small same}" sloped,
	"\text{\small ratio}"' sloped]
\ar[from = 65A,
	to = 3630,
	Leftrightarrow,
	gray,
	"\text{\small same}" sloped,
	"\text{\small ratio}"' sloped]
\end{tikzcd}

\flushleftright
As we can see, this is a form of function composition, but a strange form: We act as if we can compose the ratio $3\to 2$ with itself, even though the domain and codomain do not match.

But the strangeness is actually a more general feature of relations. Consider, for example, the biblical story of Genesis. In this story, Adam is the first man created, and so he has no father. On the other hand, his son Abel is murdered before he can have children. Thus Adam has no father, and his son Abel has no son. Nevertheless, we can meaningfully ask ``What is the relation formed by composing the relation Adam has to Abel with itself?'' The answer, of course, is grandfather to grandson. Adam has the relation of father to son to Abel, and when the ratio of father to son is composed with itself, the resulting ratio is grandfather to grandson. We can even write this composite relation similarly to the mathematical example above.

\centering%
\begin{tikzcd}[column sep = 3cm, ampersand replacement = \&, row sep = 2.5cm, math mode = false, font = \scriptsize]
Adam
\ar[r,
	"father to son" {font = \tiny},
	""' {name = 65A}]
\ar[d,
	Leftrightarrow,
	gray,
	densely dotted,
	"analogous" {font = \tiny, sloped},
	"persons"' {font = \tiny, sloped}]
\&
	Abel
	\ar[dr,
		Leftrightarrow,
		gray,
		densely dotted,
		"analogous" {pos = 0.4, font = \tiny, sloped},
		"persons"' {font = \tiny, sloped}]
	\&[-3cm]
	,
	\&[-3cm]
		Adam
		\ar[dl,
			Leftrightarrow,
			gray,
			densely dotted,
			"analogous" {pos = 0.4, font = \tiny, sloped},
			"persons"' {font = \tiny, sloped}]
		\ar[r,
			"father to son" {font = \tiny},
			""' {name = 65B}]
		\&
			Abel\ar[d,
				Leftrightarrow,
				gray,
				densely dotted,
				"analogous" {font = \tiny, sloped},
				"persons"' {font = \tiny, sloped}]\\
$A$
\ar[rr,
	"father to son"' {font = \tiny},
	"" {name = 3630}]
\ar[rrrr,
	bend right,
	"grandfather\\to grandson" {align = center, font = \tiny},
	"Adam's relation to Abel\\composed with itself"' {align = center, font = \tiny}]
\&\&
	$B$
	\ar[rr,
		"father to son"' {font = \tiny},
		"" {name = 3025}]
	\&\&
	$C$
\ar[from = 65B,
	to = 3025,
	Leftrightarrow,
	gray,
	"\small same" sloped,
	"\small relation"' sloped]
\ar[from = 65A,
	to = 3630,
	Leftrightarrow,
	gray,
	"\small same" sloped,
	"\small relation"' sloped]
\end{tikzcd}

\flushleftright
In this example, it is accidental to fatherhood that it is in Adam pointing toward Abel. And the relation itself that happens to be in Adam, namely, fatherhood, or more specifically, father-to-son, can be composed with itself. But if the relation of Adam to Abel were, for some reason, prototypical then we would speak of the composition of Adam's relation to Abel with itself---even though focusing on these particular related subjects, Adam has no father and Abel no son, and so the composition is impossible. This is like saying that the ratio of nine to four is the ratio of three to two composed with itself. Since three and two are relatively prime, they are prototypical, and the relation itself is named after these numbers.\footnote{\label{wellDefined}Acerbi \citep{acerbi2011}, page 136, argues that Greek mathematicians assumed ratio composition was well-defined, but never proved it. I offer a quick demonstration that it is well defined in Appendix \ref{CompositionAppendix}.%
}

On the other hand, if it were essential to our story that we consider a first human (\textit{Adam}) and a murdered childless child (\textit{Abel}) who is the same sex as the parent, but the relation they have is accidental---they could be father and son or mother and daughter---then again we could speak of the composition of that ratio with itself, even though we would not be able to name it. This is the situation in this Proposition: The bases of the triangles have \textit{some} relation, and that relation can be composed with itself---in the sense that ``fatherhood'' can be composed with itself. But the specific relation is itself accidental to the Proposition. We therefore cannot name it, except as the relation of the base of one triangle to the base of the other, and speak of the composition of that ratio with itself.

From at least Euclid through to Newton, compounding ratios was treated as a form of \textit{addition}---the two ratios are joined, tip-to-tail to form a new ratio.\footnote{See Sylla \citep{sylla1984} and Acerbi \citep{acerbi2011}. Though, as Acerbi notes, it isn't quite correct to say Greek mathematicians treated it as the \textit{operation }of addition. Greek mathematics used strictly additive language, but compounding ratios was not an operation: A ratio was described as a compound composed of other ratios, but they did not compound two ratios to form a compound ratio. The opposite operation, subtraction of two ratios, was, however an operation, and was described with standard subtraction language. We would describe it as multiplying and dividing ratios---since, numerically, it corresponds to multiplication, and since $f\circ f$ is written as $f^2$, not $2f$ still less two $f$'s (which is probably closer to the ancient conception)---but I would caution readers not pre-judge against Euclid and Newton here. Indeed, their language still survives in music theory. When we refer to ratios as musical intervals---as a fifth or an octave, etc.---we still call compound intervals (e.g., two fifths) the sum, not the product, of the intervals. Indeed, the original example shows, numerically, the composition of two fifths to make one major ninth.

Since on a piano the key of C is most natural and intervals paradigmatically start from C, we can even ask what interval is formed when the interval of C to G is composed with itself, even though the domain and codomain don't match. To find the composite interval we need to find two intervals each of which is the same as a C to G interval (that is, we need two fifths) with a common middle term, e.g.,

\centering%
\begin{tikzcd}[ampersand replacement = \&, math mode = false, column sep = 2cm]
$C$
\ar[r, "\textit{perfect fifth}" {font = \tiny}]
\ar[rr,
	bend right,
	"\textit{two fifths,}"{font = \tiny},
	"\textit{aka, major ninth}"' {font = \tiny}]
\&
	$G$
	\ar[r, "\textit{perfect fifth}" {font = \tiny}]
	\&
		$D$.
\end{tikzcd}

\flushleftright%
Again, note that when we write the intervals, in Pythagorean style, as ratios, and label the ratios rather than writing ``same ratio'', this diagram becomes almost the same as one that started the chapter:

\centering%
\begin{tikzcd}[ampersand replacement = \&, math mode = false, column sep = 2cm]
$9$
\ar[r, "\textit{3-to-2}" {font = \tiny}]
\ar[rr,
	bend right,
	"\textit{two 3-to-2's,}"{font = \tiny},
	"\textit{aka, 9-to-4}"' {font = \tiny}]
\&
	$6$
	\ar[r, "\textit{3-to-2}" {font = \tiny}]
	\&
		$4$.
\end{tikzcd}

\flushleftright%
} Therefore, counting the two identical ratios that make it up, we say the composite ratio is the \textit{duplicate} of the ratios that make it up.

\subsection{Proposition VI.19}
In Proposition VI.19, Euclid demonstrates that similar triangles are in a duplicate ratio of their sides. He thus sets out two similar triangles; and specifies one pair of equal angles and the analogous sides around the angle.

\centering%
\begin{tikzpicture}[commutative diagrams/every diagram]
\coordinate (B) at (0,0);
\coordinate (A) at (60:2.71828);
\coordinate (C) at (3.14159,0);

\coordinate (E) at (6.14159,0);
\coordinate (D) at ($(E) +(60:1.5336)$);
\coordinate (F) at ($(E) +(1.77245,0)$);

\coordinate (G) at (1.0472,0);

\draw (A) to node[midway] (AB) {}
	  (B) to node[midway] (BC) {}
	  (C) to
	  cycle;

\draw (D) to node[midway, inner sep = 2] (DE) {}
	  (E) to node[midway, inner sep = 2, pos = 0.45] (EF) {}
	  (F) to
	  cycle;

%\draw (G) -- (A);

\node[anchor = south, inner sep = 2] at (A) {\tiny $A$};
\node[anchor = north, inner sep = 2] at (B) {\tiny $B$};
\node[anchor = north, inner sep = 2] at (C) {\tiny $C$};
\node[anchor = south, inner sep = 2] at (D) {\tiny $D$};
\node[anchor = north, inner sep = 2] at (E) {\tiny $E$};
\node[anchor = north, inner sep = 2] at (F) {\tiny $F$};
%\node[anchor = north, inner sep = 2] at (G) {\tiny $G$};

\path[commutative diagrams/.cd, every arrow, every label]%
	(AB)
		edge[out = 15,
			 in = 60,
			 gray,
			 "" {name = ABBC, inner sep = 1}]
	(BC)
	(DE)
		edge[out = 15,
			 in = 60,
			 gray,
			 "" {name = DEEF, inner sep = 1}]
	(EF)
	(ABBC)
		edge[commutative diagrams/equal,
			 out = 45,
			 in = 45,
			 white,
			 line width = 1.5mm,
			 distance = 2.5cm]
	(DEEF)
	(ABBC)
		edge[commutative diagrams/Leftrightarrow,
			 out = 45,
			 in = 45,
			 gray,
			 "\footnotesize Same" {sloped, pos = 0.43},
			 "\footnotesize Ratio" {sloped, pos = 0.57},
			 distance = 2.5cm]
	(DEEF);

\begin{scope}
\path[clip] (A) -- (B) -- (C);
\draw[blue] (B) circle (0.27);
\end{scope}

\begin{scope}
\path[clip] (D) -- (E) -- (F);
\draw[blue] (E) circle (0.27);
\end{scope}

\end{tikzpicture}

\flushleftright
After setting out the two triangles, Euclid prepares the triangles for the proof by presenting the base of the larger triangle as divided at a point $D$ so that, as is $BC$ to $EF$ so is $EF$ to $BG$, and similarly, he presents $AG$ as drawn.\footnote{This is the section called the construction. But in \textit{The Ethics of Geometry}, Lachterman \citep{Lachterman1989} argues it should be called the \textit{preparation}. Moreover, Acerbi \citep{acerbi2021}, p. 5, \textit{et passim}) argues that Euclid's perfect passive imperatives should be understood as presenting a static being, not as commanding actions vaguely in the past. This seems to fit with Avicenna \citep{avicenna-meta} who, in the \textit{Metaphysics of the Healing} (III.4), says that when mathematicians describe the motion of a point, it is for the sake of leading the imagination, but the line through which the point moves is, in fact, prior to the motion.}

\centering%
\begin{tikzpicture}[commutative diagrams/every diagram]
\coordinate (B) at (0,0);
\coordinate (A) at (60:2.71828);
\coordinate (C) at (3.14159,0);

\coordinate (E) at (6.14159,0);
\coordinate (D) at ($(E) +(60:1.5336)$);
\coordinate (F) at ($(E) +(1.77245,0)$);

\coordinate (G) at (1,0);

\draw (A) to node[midway] (AB) {}
	  (B) to 
	  (C) to
	  cycle;

\path (B) to node[midway] (BG) {} (G);

\draw (D) to node[midway] (DE) {}
	  (E) to node[midway, pos = 0.5] (EF) {}
	  (F) to
	  cycle;

\draw (G) -- (A);

\draw[gray, decorate, decoration = {brace, mirror}]%
	 ($(B) +(0,-0.1)$) to
	 	node[midway, anchor = north, inner sep = 0, yshift = -6] (BC) {}
	 ($(C) +(0,-0.1)$);

\node[anchor = south, inner sep = 2] at (A) {\tiny $A$};
\node[anchor = east, inner sep = 2] at (B) {\tiny $B$};
\node[anchor = west, inner sep = 2] at (C) {\tiny $C$};
\node[anchor = south, inner sep = 2] at (D) {\tiny $D$};
\node[anchor = north, inner sep = 2] at (E) {\tiny $E$};
\node[anchor = north, inner sep = 2] at (F) {\tiny $F$};
\node[anchor = south west, inner sep = 2] at (G) {\tiny $G$};

\path[commutative diagrams/.cd, every arrow, every label]%
	(BC)
		edge[out = -45,
			 in = -135,
			 gray,
			 "" {name = BCEF, inner sep = 1}]
	(EF)
	(EF)
		edge[out = 135,
			 in = 45,
			 white,
			 -,
			 line width = 1.25mm]
	(BG)
	(EF)
		edge[out = 135,
			 in = 45,
			 gray,
			 "" {name = EFBG, inner sep = 1}]
	(BG)
	(BCEF)
		edge[commutative diagrams/Leftrightarrow,
			 gray!75,
			 "\tiny Same" sloped,
			 "\tiny Ratio"' sloped]
	(EFBG);

\end{tikzpicture}

\flushleftright%
Note that this means that $CB$ has, to $BG$, the duplicate ratio of its ratio to $EF$: if we treat the diagram above as a commutative diagram, this is immediately clear---though Euclid will, of course, discuss it explicitly later.

Before Euclid can use these ratios, however, he needs to alternate the original ratios.\footnote{%
Alternating the analogy is transposing the square. So for example, fifteen is to ten as six is to four.

\centering%
\begin{tikzcd}[ampersand replacement = \&]
15
\ar[r,
	""' {name = 1510}]
\ar[d,
	gray,
	Leftrightarrow,
	densely dotted]
\&
	10
	\ar[d,
		gray,
		Leftrightarrow,
		densely dotted]\\
6
\ar[r,
	"" {name = 64}]
\&
	4
\ar[from = 1510,
	to = 64,
	Leftrightarrow,
	gray,
	"\text{same}" sloped,
	"\text{ratio}"' sloped]
\end{tikzcd}%

\flushleftright%
When we alternate this analogy, we produce the analogy fifteen is to six as ten is to four:

\centering%
\begin{tikzcd}[ampersand replacement = \&]
15
\ar[d,
	"" {name = 156}]
\ar[r,
	gray,
	Leftrightarrow,
	densely dotted]
\&
	10
	\ar[d,
		""' {name = 104}]\\
6
\ar[r,
		gray,
		Leftrightarrow,
		densely dotted]
\&
	4
\ar[from = 156,
	to = 104,
	Leftrightarrow,
	gray,
	"\text{same}" sloped,
	"\text{ratio}"' sloped]
\end{tikzcd}

\flushleftright%
In V.16, Euclid proves that if all four magnitudes are homogeneous, then if the quantities are analogous, they will still be analogous after the analogy is alternated. That is, here, since fifteen is to ten as six is to four, fifteen is to \textit{six} as ten is to four.
}
\begin{tcolorbox}[%
	colbacktitle=white,
	coltitle = black,
	center title,
	fonttitle = \bfseries,
	colbacklower = white,
	colback = white,
	width = 1\linewidth,
	lower separated=true,
	boxrule=.2mm]
Since,

\centering%
\begin{tikzpicture}[commutative diagrams/every diagram]
\coordinate (B) at (0,0);
\coordinate (A) at (60:2.71828);
\coordinate (C) at (3.14159,0);

\coordinate (E) at (6.14159,0);
\coordinate (D) at ($(E) +(60:1.5336)$);
\coordinate (F) at ($(E) +(1.77245,0)$);

\coordinate (G) at (1.0472,0);

\draw (A) to node[midway] (AB) {}
	  (B) to node[midway] (BC) {}
	  (C) to
	  cycle;

\draw (D) to node[midway, inner sep = 2] (DE) {}
	  (E) to node[midway, inner sep = 2, pos = 0.45] (EF) {}
	  (F) to
	  cycle;

\draw[gray!50] (G) -- (A);

\node[anchor = south, inner sep = 2] at (A) {\tiny $A$};
\node[anchor = north, inner sep = 2] at (B) {\tiny $B$};
\node[anchor = north, inner sep = 2] at (C) {\tiny $C$};
\node[anchor = south, inner sep = 2] at (D) {\tiny $D$};
\node[anchor = north, inner sep = 2] at (E) {\tiny $E$};
\node[anchor = north, inner sep = 2] at (F) {\tiny $F$};
\node[anchor = north, inner sep = 2, gray!50] at (G) {\tiny $G$};

\path[commutative diagrams/.cd, every arrow, every label]%
	(AB)
		edge[out = 15,
			 in = 60,
			 gray,
			 "" {name = ABBC, inner sep = 1}]
	(BC)
	(DE)
		edge[out = 15,
			 in = 60,
			 gray,
			 "" {name = DEEF, inner sep = 1}]
	(EF)
	(ABBC)
		edge[commutative diagrams/equal,
			 out = 45,
			 in = 45,
			 white,
			 line width = 1.5mm,
			 distance = 2.5cm]
	(DEEF)
	(ABBC)
		edge[commutative diagrams/Leftrightarrow,
			 out = 45,
			 in = 45,
			 gray,
			 "\footnotesize Same" {sloped, pos = 0.43},
			 "\footnotesize Ratio" {sloped, pos = 0.57},
			 distance = 2.5cm]
	(DEEF);

\end{tikzpicture}
\tcblower%
Hence, alternately,

\centering%
\begin{tikzpicture}[commutative diagrams/every diagram]
\coordinate (B) at (0,0);
\coordinate (A) at (60:2.71828);
\coordinate (C) at (3.14159,0);

\coordinate (E) at (6.14159,0);
\coordinate (D) at ($(E) +(60:1.5336)$);
\coordinate (F) at ($(E) +(1.77245,0)$);

\coordinate (G) at (1.0472,0);

\draw (A) to node[midway] (AB) {}
	  (B) to node[midway] (BC) {}
	  (C) to
	  cycle;

\draw (D) to node[midway] (DE) {}
	  (E) to node[midway, pos = 0.5] (EF) {}
	  (F) to
	  cycle;

\draw[gray!50] (G) -- (A);

\node[anchor = south, inner sep = 2] at (A) {\tiny $A$};
\node[anchor = north, inner sep = 2] at (B) {\tiny $B$};
\node[anchor = north, inner sep = 2] at (C) {\tiny $C$};
\node[anchor = south, inner sep = 2] at (D) {\tiny $D$};
\node[anchor = north, inner sep = 2] at (E) {\tiny $E$};
\node[anchor = north, inner sep = 2] at (F) {\tiny $F$};
\node[anchor = north, inner sep = 2, gray!50] at (G) {\tiny $G$};

\path[commutative diagrams/.cd, every arrow, every label]%
	(AB)
		edge[out = 130,
			 in = 130,
			 gray,
			 distance = 3cm,
			 ""' {name = ABDE, pos = 0.65, inner sep = 1}]
	(DE)
	(BC)
		edge[out = -135,
			 in = -45,
			 gray,
			 "" {name = BCEF, inner sep = 1}]
	(EF)
	(ABDE)
		edge[commutative diagrams/Leftrightarrow,
			 out = -110,
			 in = 90,
			 gray,
			 "\tiny Same" {sloped, pos = 0.5},
			 "\tiny Ratio"' {sloped, pos = 0.5}]
	(BCEF);

\end{tikzpicture}
\end{tcolorbox}

Then, since ``same ratio'' is transitive, it follows that $AB$ is to $DE$ as $EF$ is to $BG$.

\centering%
\begin{tikzpicture}[commutative diagrams/every diagram]
\coordinate (B) at (0,0);
\coordinate (A) at (60:2.71828);
\coordinate (C) at (3.14159,0);

\coordinate (E) at (6.14159,0);
\coordinate (D) at ($(E) +(60:1.5336)$);
\coordinate (F) at ($(E) +(1.77245,0)$);

\coordinate (G) at (1,0);

\draw (A) to node[midway] (AB) {}
	  (B) to 
	  (C) to
	  cycle;

\path (B) to node[midway] (BG) {} (G);

\draw (D) to node[midway] (DE) {}
	  (E) to node[midway, pos = 0.5] (EF) {}
	  (F) to
	  cycle;

\draw (G) -- (A);

\draw[gray, decorate, decoration = {brace, mirror}]%
	 ($(B) +(0,-0.1)$) to
	 	node[midway, anchor = north, inner sep = 0, yshift = -6] (BC) {}
	 ($(C) +(0,-0.1)$);

\node[anchor = south, inner sep = 2] at (A) {\tiny $A$};
\node[anchor = east, inner sep = 2] at (B) {\tiny $B$};
\node[anchor = west, inner sep = 2] at (C) {\tiny $C$};
\node[anchor = south, inner sep = 2] at (D) {\tiny $D$};
\node[anchor = north, inner sep = 2] at (E) {\tiny $E$};
\node[anchor = north, inner sep = 2] at (F) {\tiny $F$};
\node[anchor = south west, inner sep = 2] at (G) {\tiny $G$};

\path[commutative diagrams/.cd, every arrow, every label]%
	(AB)
		edge[out = 130,
			 in = 130,
			 gray,
			 distance = 3cm,
			 ""' {name = ABDE, pos = 0.85, inner sep = 1},
			 ""' {name = ABDE2, pos = 0.5, inner sep = 1}]
	(DE)
	(BC)
		edge[out = -45,
			 in = -135,
			 gray,
			 "" {name = BCEF, inner sep = 1, pos = 0.45},
			 "" {name = BCEF2, inner sep = 2, pos = 0.6}]
	(EF)
	(EF)
		edge[out = 135,
			 in = 45,
			 white,
			 -,
			 line width = 1.25mm]
	(BG)
	(EF)
		edge[out = 135,
			 in = 45,
			 gray,
			 "" {name = EFBG, inner sep = 1},
			 ""' {name = EFBG2, pos = 0.6, inner sep = 1}]
	(BG)
	(BCEF)
		edge[commutative diagrams/Leftrightarrow,
			 gray!75,,
			 out = 90,
			 in = -80,
			 "\tiny Same" sloped,
			 "\tiny Ratio"' sloped]
	(EFBG)
	(ABDE)
		edge[commutative diagrams/equal,
			 white,
			 out = -100,
			 in = 90,
			 line width = 1.75mm]
	(BCEF2)
	(ABDE)
		edge[commutative diagrams/Leftrightarrow,
			 gray!75,
			 out = -100,
			 in = 90,
			 "\tiny Same" sloped,
			 "\tiny Ratio"' sloped]
	(BCEF2)
	(ABDE2)
		edge[commutative diagrams/phantom,
			 BrickRed!50,
			 out = -70,
			 in = 120,
			 "{\phantom{\tiny \begin{tabular}{@{}c@{}} Hence,\\[-0.5ex] Same\end{tabular}\vspace{-1ex}}}" {sloped, inner sep = 0.5, yshift = 1.5, fill = white, fill},
			 "{\tiny \begin{tabular}{@{}c@{}} Hence,\\[-0.5ex] Same\end{tabular}\vspace{-1ex}}" {sloped, inner sep = 0.5, yshift = 1},
			 "\tiny Ratio"' sloped]
	(EFBG2)
		(ABDE2)
		edge[commutative diagrams/Leftrightarrow,
			 BrickRed!50,
			 out = -70,
			 in = 120]
	(EFBG2);

\end{tikzpicture}

\flushleftright%
And though it is not immediately obvious, Euclid, in Proposition VI.15 Euclid used a series of analogies between sides and areas of triangles to prove that this, in turn, implies that triangle $ABG$ is equal to triangle $DEF$.

\centering%
\begin{tcolorbox}[%
	colbacktitle=white,
	coltitle = black,
	center title,
	halign upper=center,
	fonttitle = \bfseries,
	colbacklower = white,
	colback = white,
	width = 0.8\linewidth,
	lower separated=true,
	boxrule=.2mm]
\begin{tikzpicture}[commutative diagrams/every diagram]
\coordinate (B) at (0,0);
\coordinate (A) at (60:2.71828);
\coordinate (C) at (3.14159,0);

\coordinate (E) at (6.14159,0);
\coordinate (D) at ($(E) +(60:1.5336)$);
\coordinate (F) at ($(E) +(1.77245,0)$);

\coordinate (G) at (1,0);

\draw (A) to node[midway] (AB) {}
	  (B) to 
	  (G) to
	  cycle;
\draw[gray!50] (G) -- (C) -- (A);

\path (B) to node[midway] (BG) {} (G);

\draw (D) to node[midway] (DE) {}
	  (E) to node[midway, pos = 0.5] (EF) {}
	  (F) to
	  cycle;

\draw (G) -- (A);

\node[anchor = south, inner sep = 2] at (A) {\tiny $A$};
\node[anchor = north, inner sep = 2] at (B) {\tiny $B$};
\node[anchor = north, inner sep = 2, gray!50] at (C) {\tiny $C$};
\node[anchor = south, inner sep = 2] at (D) {\tiny $D$};
\node[anchor = north, inner sep = 2] at (E) {\tiny $E$};
\node[anchor = north, inner sep = 2] at (F) {\tiny $F$};
\node[anchor = north, inner sep = 2] at (G) {\tiny $G$};

\path[commutative diagrams/.cd, every arrow, every label]%
	(AB)
		edge[out = 130,
			 in = 130,
			 gray,
			 distance = 3cm,
			 ""' {name = ABDE, pos = 0.5, inner sep = 1}]
	(DE)
	(EF)
		edge[out = 135,
			 in = 45,
			 white,
			 -,
			 line width = 1.25mm]
	(BG)
	(EF)
		edge[out = 135,
			 in = 45,
			 gray,
			 ""' {name = EFBG, pos = 0.6, inner sep = 1}]
	(BG)
	(ABDE)
		edge[commutative diagrams/Leftrightarrow,
			 gray!75,
			 out = -70,
			 in = 120,
			 "\tiny Same" sloped,
			 "\tiny Ratio"' sloped]
	(EFBG);

\begin{scope}
\path[clip] (A) -- (B) -- (G);
\draw[blue] (B) circle (0.27);
\end{scope}

\begin{scope}
\path[clip] (D) -- (E) -- (F);
\draw[blue] (E) circle (0.27);
\end{scope}
\end{tikzpicture}
\tcblower
Hence,

\centering
\begin{tikzpicture}[commutative diagrams/every diagram]
\coordinate (B) at (0,0);
\coordinate (A) at (60:2.71828);
\coordinate (C) at (3.14159,0);

\coordinate (E) at (6.14159,0);
\coordinate (D) at ($(E) +(60:1.5336)$);
\coordinate (F) at ($(E) +(1.77245,0)$);

\coordinate (G) at (1,0);

\draw[gray!50] (G) -- (C) -- (A);
\filldraw[fill = Violet!50]%
	  (A) to node[midway] (AB) {}
	  (B) to 
	  (G) to node[pos = 0.52] (AG) {}
	  		 node[pos = 0.48] (AG2) {}
	  		 node[pos = 0.5] (AG3) {}
	  cycle;

\path (B) to node[midway] (BG) {} (G);

\filldraw[fill = ForestGreen!50]%
	  (D) to node[pos = 0.4689] (DE) {}
	  		 node[pos = 0.531] (DE2) {}
	  		 node[pos = 0.5] (DE3) {}
	  (E) to node[midway, pos = 0.5] (EF) {}
	  (F) to
	  cycle;

\draw (G) -- (A);

\node[anchor = south, inner sep = 2] at (A) {\tiny $A$};
\node[anchor = north, inner sep = 2] at (B) {\tiny $B$};
\node[anchor = north, inner sep = 2, gray!50] at (C) {\tiny $C$};
\node[anchor = south, inner sep = 2] at (D) {\tiny $D$};
\node[anchor = north, inner sep = 2] at (E) {\tiny $E$};
\node[anchor = north, inner sep = 2] at (F) {\tiny $F$};
\node[anchor = north, inner sep = 2] at (G) {\tiny $G$};

\path[commutative diagrams/.cd, every arrow, every label]%
	(AG3)
		edge[out = -8.67414,
			 in = 150,
			 commutative diagrams/equal,
			 white,
			 line width = 3.5]
	(DE3)
	(AG)
		edge[out = -8.67414,
			 in = 150,
			 commutative diagrams/rightharpoonup,
			 "\tiny Equal to" sloped]
	(DE)
	(AG2)
		edge[out = -8.67414,
			 in = 150,
			 commutative diagrams/leftharpoondown,
			 "\scalebox{-1}{\tiny Equal to}"' sloped]
	(DE2);

\end{tikzpicture}
\end{tcolorbox}

\flushleftright%
Now in VI.1, Euclid showed that triangles in the same parallels are to each other as their bases. Therefore triangle $ABC$ is to triangle $ABG$ as $BC$ is to $BG$---that is, $ABC$ has, to $ABG$, the duplicate ratio that $BC$ has to $EF$.

\centering%
\begin{tikzpicture}[commutative diagrams/every diagram]
\coordinate (B) at (0,0);
\coordinate (A) at (60:2.71828);
\coordinate (C) at (3.14159,0);

\coordinate (E) at (6.14159,0);
\coordinate (D) at ($(E) +(60:1.5336)$);
\coordinate (F) at ($(E) +(1.77245,0)$);

\coordinate (G) at (1,0);

\fill[blue!50] (A) -- (G) -- (C);
\fill[Violet!50] (A) -- (B) -- (G);

\draw[gray!25] (G) -- (A);

\draw (A) to node[midway] (AB) {}
	  (B) to 
	  (C) to
	  cycle;

\draw[gray!50] (F) to (D) to node[midway] (DE) {}
	  (E);
\draw (E) to node[pos = 0.5] (EF) {}
	  		 node[pos = 0.55, inner sep = 0] (EF3) {}
	  		 node[pos = 0.55, inner sep = 0] (EF4) {}
	  (F);

\draw[gray, decorate, decoration = {brace, mirror}]%
	 ($(B) +(0,-0.1)$) to
	 	node[pos = 0.45, anchor = north, inner sep = 0, yshift = -6] (BC) {}
	 	node[pos = 0.55, anchor = north, inner sep = 0, yshift = -6] (BC2) {}
	 ($(C) +(0,-0.1)$);
	 
\path[name path=BCc] ($(B)!0.5!(C)$) -- (A);
\path[name path=ABc] ($(B)!0.5!(A)$) -- (C);
\path [name intersections={of=BCc and ABc,by={ABCc}}];
\node[inner sep = 0] (ABCc2) at (ABCc) {};

\path[name path=BGc2] ($(B)!0.5!(G)$) -- (A);
\path[name path=ABc2] ($(B)!0.5!(A)$) -- (G);
\path [name intersections={of=BGc2 and ABc2,by={ABGc}}];

\coordinate[inner sep = 0] (EF2) at ($(C)!0.25!(E)$);

\draw[gray, decorate, decoration = {brace, mirror}]%
	 ($(B) +(0,-0.25)$) to
	 	node[pos = 0.6, anchor = north, inner sep = 0, yshift = -4] (BG) {}
	 	node[pos = 0.43, anchor = north, inner sep = 0, yshift = -6] (BG2) {}
	 ($(G) +(0,-0.25)$);

\node[anchor = south, inner sep = 2] at (A) {\tiny $A$};
\node[anchor = east, inner sep = 2] at (B) {\tiny $B$};
\node[anchor = west, inner sep = 2] at (C) {\tiny $C$};
\node[anchor = south, inner sep = 2, gray!50] at (D) {\tiny $D$};
\node[anchor = north, inner sep = 2] at (E) {\tiny $E$};
\node[anchor = north, inner sep = 2] at (F) {\tiny $F$};
\node[anchor = south west, inner sep = 2, gray!25] at (G) {\tiny $G$};

\draw
	(BC2)
		edge[out = -30,
			 in = -150,
			 gray,
			 -,
			 "" {name = BCEF2, inner sep = 1}]
	(EF2)
	(EF2)
		edge[out = 30,
			 in = 60,
			 white,
			 -,
			 distance = 1.75cm,
			 line width = 0.75mm]
	($(EF3) +(240:0.05)$)
	(EF2)
		edge[out = 30,
			 in = 60,
			 gray,
			 -,
			 distance = 1.75cm,
			 "\tiny Duplicate Ratio" {pos = 0.2, red!50, sloped, inner sep = 1, name = BCEF2}]
	($(EF3) +(240:0.05)$);

\path[commutative diagrams/.cd, every arrow, every label]%
	(BC)
		edge[out = -45,
			 in = 135,
			 white,
			 -,
			 line width = 0.75mm]
	(EF)
	(BC)
		edge[out = -45,
			 in = 135,
			 gray!60,
			 "\tiny One"' {pos = 0.67, red!50, sloped},
			 ""' {name = BCEF, inner sep = 1}]
	(EF)
	(EF)
		edge[out = -135,
			 in = -45,
			 gray!60,
			 "\tiny Two" {pos = 0.67, red!50, sloped},
			 ""' {name = EFBG, inner sep = 1}]
	(BG)
	(EF4)
		edge[out = -120,
			 in = -45,
			 gray,
			 "\tiny Duplicate Ratio"' {red!50, pos = 0.53, sloped, inner sep = 2}]
	(BG2)
	(ABCc2)
		edge[out = 30,
			 in = 150,
			 distance = 6cm,
			 black!75,
			 "\tiny Whole {\color{Plum}Tri}{\color{blue}angle}"',
			 "\tiny to its {\color{Plum}Part}",
			 ""' {pos = 0.2, inner sep = 1, name = ABCABG}]
	(ABGc)
	(BCEF)
		edge[commutative diagrams/Leftrightarrow,
			 gray!50,
			 out = -45,
			 in = 90,
			 "\tiny Same" sloped,
			 "\tiny Ratio"' sloped]
	(EFBG)
	(ABCABG)
		edge[commutative diagrams/Leftrightarrow,
			 gray,
			 out = 0,
			 in = 120,
			 "\tiny Same" {sloped, pos = 0.55},
			 "\tiny Ratio"' {sloped, pos = 0.55}]
	(BCEF2)
;
\end{tikzpicture}

\flushleftright%
But $ABG$ is equal to $DEF$, and therefore $ABG$ has, to $DEF$ the duplicate ratio that $BC$ has to $EF$.

\centering%
\begin{tcolorbox}[%
	colbacktitle=white,
	coltitle = black,
	center title,
	fonttitle = \bfseries,
	colbacklower = white,
	colback = white,
	width = 0.8\linewidth,
	boxrule=.2mm]
\hspace{-3.25cm}
\begin{tikzpicture}[commutative diagrams/every diagram]
\coordinate (B) at (0,0);
\coordinate (A) at (60:2.71828);
\coordinate (C) at (3.14159,0);

\coordinate (E) at (6.14159,0);
\coordinate (D) at ($(E) +(60:1.5336)$);
\coordinate (F) at ($(E) +(1.77245,0)$);

\coordinate (G) at (1,0);

\coordinate (B') at (0,4.5);
\coordinate (A') at ($(B') +(60:2.71828)$);
\coordinate (C') at ($(B') +(3.14159,0)$);

\coordinate (E') at ($(B') +(6.14159,0)$);
\coordinate (D') at ($(B') +($(E) +(60:1.5336)$)$);
\coordinate (F') at ($(B') +($(E) +(1.77245,0)$)$);

\coordinate (G') at ($(B') +(1,0)$);

\fill[blue!50] (A) -- (G) -- (C);
\fill[Violet!50] (A) -- (B) -- (G);

\draw[gray!25] (G) -- (A);

\draw (A) to node[midway] (AB) {}
	  (B) to 
	  (C) to
	  cycle;

\draw[gray!50] (F) to (D) to node[midway] (DE) {}
	  (E);
\draw (E) to node[pos = 0.5] (EF) {}
	  		 node[pos = 0.55, inner sep = 0] (EF3) {}
	  		 node[pos = 0.55, inner sep = 0] (EF4) {}
	  (F);

\draw[gray, decorate, decoration = {brace, mirror}]%
	 ($(B) +(0,-0.1)$) to
	 	node[pos = 0.45, anchor = north, inner sep = 0, yshift = -6] (BC) {}
	 	node[pos = 0.55, anchor = north, inner sep = 0, yshift = -6] (BC2) {}
	 ($(C) +(0,-0.1)$);

\fill[blue!50]%
	  (A') to
	  (B') to
	  (C') to
	  cycle;
	  
\fill[ForestGreen!50]%
	  (D') to
	  (E') to
	  (F') to
	  cycle;

\draw[gray!25] (A') -- (G');
	 
\draw%
	  (A') to
	  (B') to
	  (C') to
	  cycle;
	  
\draw%
	  (D') to
	  (E') to
	  (F') to
	  cycle;
	 
\path[name path=BCc] ($(B)!0.5!(C)$) -- (A);
\path[name path=ABc] ($(B)!0.5!(A)$) -- (C);
\path [name intersections={of=BCc and ABc,by={ABCc}}];
\node[inner sep = 0] (ABCc2) at (ABCc) {};

\path[name path=BGc2] ($(B)!0.5!(G)$) -- (A);
\path[name path=ABc2] ($(B)!0.5!(A)$) -- (G);
\path [name intersections={of=BGc2 and ABc2,by={ABGc}}];

%Top Triangles
\path[name path=BCc'] ($(B')!0.5!(C')$) -- (A');
\path[name path=ABc'] ($(B')!0.5!(A')$) -- (C');
\path [name intersections={of=BCc' and ABc',by={ABCc'}}];
\node[inner sep = 0] (ABCc2') at (ABCc') {};

\path[name path=EFc] ($(E')!0.5!(F')$) -- (D');
\path[name path=DEc] ($(D')!0.5!(E')$) -- (F');
\path [name intersections={of=EFc and DEc,by={DEFc}}];

\coordinate[inner sep = 0] (EF2) at ($(C)!0.25!(E)$);

\draw[gray, decorate, decoration = {brace, mirror}]%
	 ($(B) +(0,-0.25)$) to
	 	node[pos = 0.6, anchor = north, inner sep = 0, yshift = -4] (BG) {}
	 	node[pos = 0.43, anchor = north, inner sep = 0, yshift = -6] (BG2) {}
	 ($(G) +(0,-0.25)$);

\node[anchor = south, inner sep = 2] at (A) {\tiny $A$};
\node[anchor = east, inner sep = 2] at (B) {\tiny $B$};
\node[anchor = west, inner sep = 2] at (C) {\tiny $C$};
\node[anchor = south, inner sep = 2, gray!50] at (D) {\tiny $D$};
\node[anchor = north, inner sep = 2] at (E) {\tiny $E$};
\node[anchor = north, inner sep = 2] at (F) {\tiny $F$};
\node[anchor = south west, inner sep = 2, gray!25] at (G) {\tiny $G$};

\node[anchor = south, inner sep = 2] at (A') {\tiny $A$};
\node[anchor = east, inner sep = 2] at (B') {\tiny $B$};
\node[anchor = west, inner sep = 2] at (C') {\tiny $C$};
\node[anchor = south, inner sep = 2] at (D') {\tiny $D$};
\node[anchor = north, inner sep = 2] at (E') {\tiny $E$};
\node[anchor = north, inner sep = 2] at (F') {\tiny $F$};
\node[anchor = south west, inner sep = 2, gray!25] at (G') {\tiny $G$};

\draw
	(BC2)
		edge[out = -30,
			 in = -150,
			 gray,
			 -,
			 "" {name = BCEF2, inner sep = 1}]
	(EF2)
	(EF2)
		edge[out = 30,
			 in = 60,
			 white,
			 -,
			 distance = 1.75cm,
			 line width = 0.75mm]
	($(EF3) +(240:0.05)$)
	(EF2)
		edge[out = 30,
			 in = 60,
			 gray,
			 -,
			 distance = 1.75cm,
			 "\tiny Duplicate Ratio" {pos = 0.2, red!50, sloped, inner sep = 1, name = BCEF2}]
	($(EF3) +(240:0.05)$);

\path[commutative diagrams/.cd, every arrow, every label]%
	(BC)
		edge[out = -45,
			 in = 135,
			 white,
			 -,
			 line width = 0.75mm]
	(EF)
	(BC)
		edge[out = -45,
			 in = 135,
			 gray!60,
			 "\tiny One"' {pos = 0.67, red!50, sloped},
			 ""' {name = BCEF, inner sep = 1}]
	(EF)
	(EF)
		edge[out = -135,
			 in = -45,
			 gray!60,
			 "\tiny Two" {pos = 0.67, red!50, sloped},
			 ""' {name = EFBG, inner sep = 1}]
	(BG)
	(EF4)
		edge[out = -120,
			 in = -45,
			 gray,
			 "\tiny Duplicate Ratio"' {red!50, pos = 0.53, sloped, inner sep = 2}]
	(BG2)
	(ABCc2)
		edge[out = 30,
			 in = 150,
			 distance = 6cm,
			 black!75,
			 "\tiny Whole {\color{Plum}Tri}{\color{blue}angle}"',
			 "\tiny to its {\color{Plum}Part}",
			 ""' {pos = 0.2, inner sep = 1, name = ABCABG},
			 ""' {pos = 0.3, inner sep = 1, name = ABCABG2}]
	(ABGc)
	(ABCc2')
		edge[out = 30,
			 in = 150,
			 black!75,
			 "\tiny Triangle" {blue, pos = 0.33, sloped, anchor = mid, yshift = 6},
			 "\tiny to" {sloped, anchor = mid, yshift = 6},
			 "\tiny Triangle" {ForestGreen, pos = 0.67, sloped, anchor = mid, yshift = 6},
			 ""' {inner sep = 1, name = ABCDEF}]
	(DEFc)
	(BCEF)
		edge[commutative diagrams/Leftrightarrow,
			 gray!50,
			 out = -45,
			 in = 90,
			 "\tiny Same" sloped,
			 "\tiny Ratio"' sloped]
	(EFBG)
	(ABCABG)
		edge[commutative diagrams/Leftrightarrow,
			 gray,
			 out = 0,
			 in = 120,
			 "\tiny Same" {sloped, pos = 0.55},
			 "\tiny Ratio"' {sloped, pos = 0.55}]
	(BCEF2)
	(ABCDEF)
			edge[commutative diagrams/Leftrightarrow,
			 gray,
			 out = -90,
			 in = 45,
			 "\tiny Same" {sloped, pos = 0.55},
			 "\tiny Ratio"' {sloped, pos = 0.55}]
	(ABCABG2)
;
\end{tikzpicture}

\footnotesize Note: {\color{blue} triangle $ABC$} has the same ratio to {\color{ForestGreen}$DEF$} and to {\color{Plum}$ABG$}, since {\color{ForestGreen}$DEF$} is equal to {\color{Plum}$ABG$}.
\end{tcolorbox}

\flushleftright%
Which, as Euclid says, is what was to be shown.

\section{Conclusion}

These diagrams may seem to impose a new and foreign syntax on Euclid. However, similar diagrams were used by monochordists and cosmologists throughout the Mediterranean world till at least the Renaissance. For example, Vat.gr.221 76 \citep{BAV_VatGr221_Ptolemy}, a manuscript of book I of Ptolemy's \textit{Harmonica,} has the following two diagrams, transcribed, here, with familiar Arabic numerals.\footnote{I used AI to help with transcribing the Greek calligraphy in the top diagram---and confirm my reading of the bottom diagram. I also had it help with the words \textit{homalou diatonon} (at the bottom), which are a label on the whole diagram.}

\centering%
\begin{tikzcd}[math mode = false, column sep = 1.5cm]
12
\ar[-,
	r,
	out = -90, in = -90]
\ar[r,
	shift right = 1.5,
	phantom,
	"\textit{epi 11$^{\text{os}}$}" {font = \footnotesize}]
\ar[rrr,
	-,
	out = -110, in = -80,
	"\textit{epitritos}" {font = \footnotesize}]
&
	11
	\ar[-,
		r,
		out = -90, in = -90]
	\ar[r,
		shift right = 1.5,
		phantom,
		"\textit{epi 10$^{\text{os}}$}" {font = \footnotesize}]
	&
		10
		\ar[-,
			r,
			out = -90, in = -90]
		\ar[r,
			phantom,
			shift right = 1.5,
			"\textit{epi 9$^{\text{os}}$}" {font = \footnotesize}]
		&
			9
\end{tikzcd}

\begin{tikzcd}[column sep = 1.5cm]
24
\ar[-,
	r,
	out = -90, in = -90]
\ar[r,
	phantom,
	shift right = 1.5,
	"\textit{epi 8$^{\text{os}}$}" {font = \footnotesize}]
\ar[-,
	rrrr,
	out = -100, in = -80, end anchor = {[xshift = 0.25mm]south},
	"\textit{hemiolos diapente}" {font = \footnotesize},
	"\textit{homalou diatonon}"' {font = \footnotesize}]
&
	27
	\ar[-,
		r,
		out = -90, in = -90]
	\ar[r,
		phantom,
		shift right = 1.5,
		"\textit{epi 9$^{\text{os}}$}" {font = \footnotesize}]
	\ar[rrr,
		-,
		out = -90,
		in = -80,
		"\textit{epi 3$^{\text{os}}$ diatessaron}" {font = \footnotesize}]
	&
		30
		\ar[-,
			r,
			out = -90, in = -90]
		\ar[r,
			phantom,
			shift right = 1.5,
			"\textit{epi 10$^{\text{os}}$}" {font = \footnotesize}]
		&
			33
			\ar[-,
				r,
				out = -90, in = -90]
			\ar[r,
				phantom,
				shift right = 1.5,
				"\textit{epi 11$^{\text{os}}$}" {font = \footnotesize}]
			&
				36
\end{tikzcd}

\flushleftright
Here, \textit{epi $n^{\text{os}}$} (for example, \textit{epi $11^{\text{os}}$}) is an abbreviation for the standard Greek term for the \textit{ratio} $(n+1)\xrightarrow{to}n$, and \textit{hemiolos} is the ratio of three to two,\footnote{Like a hemiola in music.} and epitritos (abbreviated elsewhere \textit{epi $3^{\text{os}}$}) the ratio of four to three. The two words \textit{diatessaron} and \textit{diapente} are musical words: fourth and fifth. \textit{Homalou diatonon} is a label on the whole diagram, distinguishing the tuning system from other tuning systems.

The English names for numerical relations have all but become obsolete, so to translate this into English we have to use symbols like e.g., $9\xrightarrow{to}8$. Leaving off the label on the whole diagram, and moving the upper labels to the arcs, these can be written as:

\centering%
\begin{tikzcd}[column sep = 1.5cm]
12
\ar[-,
	r,
	out = -90, in = -90,
	"12\;\xrightarrow{to}\,11" {yshift = 2pt}]
\ar[rrr,
	-,
	out = -110, in = -80,
	"4\xrightarrow{to}3" {font = \footnotesize}]
&
	11
	\ar[-,
		r,
		out = -90, in = -90,
		"11\,\xrightarrow{to}\,10" {yshift = 2pt}]
	&
		10
		\ar[-,
			r,
			out = -90, in = -90,
			"10\;\xrightarrow{to}\;9" {yshift = 2pt}]
		&
			9
\end{tikzcd}

\begin{tikzcd}[column sep = 1.5cm]
24
\ar[-,
	r,
	out = -90, in = -90,
	"8\;\leftarrow\, 9" {yshift = 1.5pt}]
\ar[-,
	rrrr,
	out = -100, in = -80, end anchor = {[xshift = 0.25mm]south},
	"2\leftarrow 3; \text{ fifth}" {font = \footnotesize}]
&
	27
	\ar[-,
		r,
		out = -90, in = -90,
		"9\;\leftarrow\, 10" {yshift = 1.5pt}]
	\ar[rrr,
		-,
		out = -90,
		in = -80,
		"3\leftarrow 4; \text{ fourth}" {font = \footnotesize}]
	&
		30
		\ar[-,
			r,
			out = -90, in = -90,
			"10\;\leftarrow\, 11" {yshift = 1.5pt}]
		&
			33
			\ar[-,
				r,
				out = -90, in = -90,
				"11\;\leftarrow\, 12" {yshift = 1.5pt}]
			&
				36
\end{tikzcd}

\flushleftright
Though these arcs are not arrows---though relations form something like a dagger category, and the directionality is clear---it is very easy to read these as commutative diagrams showing relational composition. Indeed, it is difficult \textit{not} to read them as commutative diagrams showing relational composition.

Diagrams like this were ubiquitous throughout the Mediterranean and West Asian intellectual worlds, in books written in Latin, Greek and Arabic. For example, here is an example from a renaissance manuscript of Boethius's \textit{de Musica} \citep{Penn_LJS47_Boethius}.

\centering%
\includegraphics[width = 9cm]{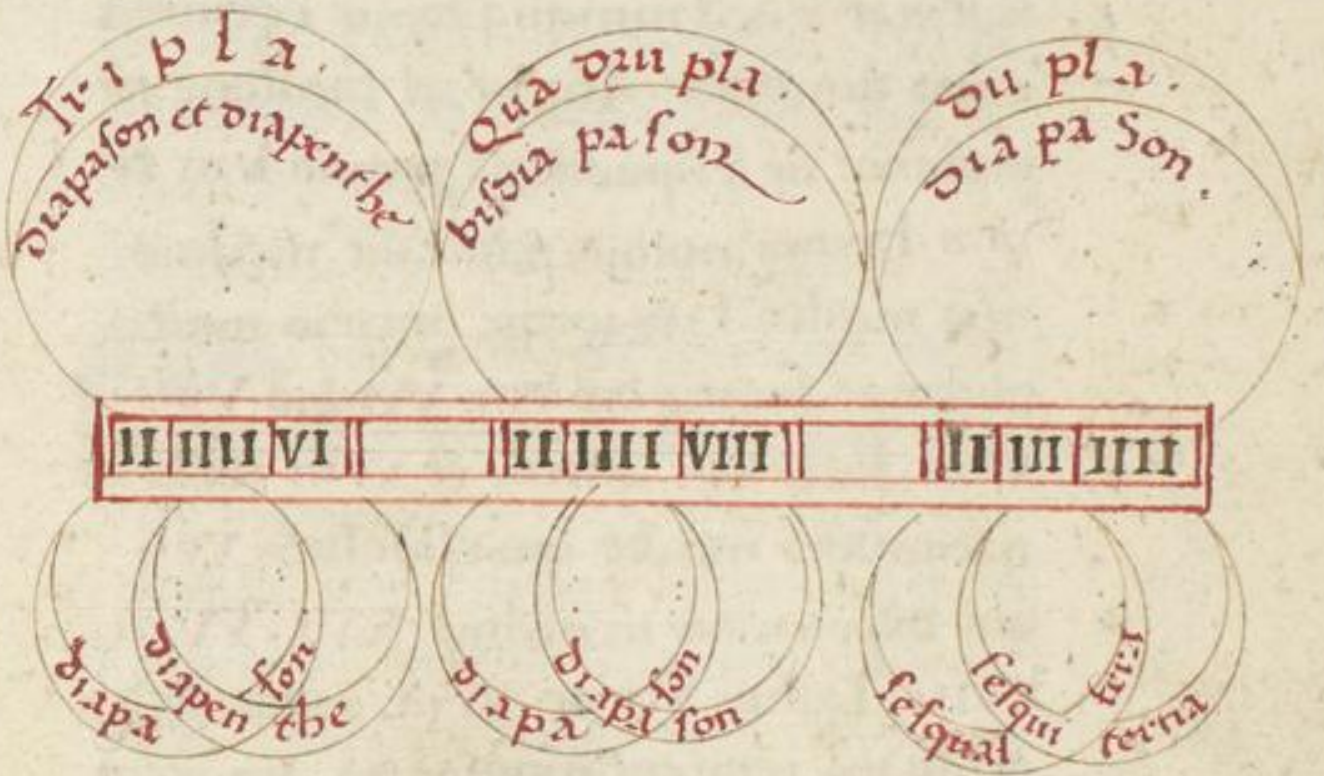}

\flushleftright%
Once again, numbers are written on a line, and the ratios or musical intervals are illustrated with arcs connecting the relevant numbers, and labeled with names of mathematical ratios\footnote{\textit{tripla} (triple), \textit{quadrupla} (quadruple), \textit{dupla} (double), \textit{sesquitertia} (four to three) and \textit{sesquialtera} (three to two)} and musical intervals.\footnote{\textit{diapason et diapenthe} (octave and fifth), \textit{bisdiapason} (twice octave), \textit{diapason} (octave) and \textit{diapenthe} (fifth).} The geometry of the two ends of the ratio arcs is different than before, but again, they show relational composition---for example, on the right side that the ratios of \( \text{II}\leftarrow \text{III}\leftarrow \text{IIII} \) compose to give the ratio \textit{double}, that is, \(\text{I} \leftarrow \text{II}.  \)

I have not surveyed this notation, and finding its historical provenance is difficult because it can differ slightly from manuscript to manuscript. Occasionally similar looking banners are used, in the same diagram, to show the relation of numbers and the function of the ratios. Other times, when distances on a monochord need labeled, banners the look like ones used elsewhere for ratios are used to label distances (numbers).\footnote{For examples of each of these, see Fogliano's \textit{Musica theorica} \protect\citep{fogliano1529musica}. Pages 38--42 of the scan show gorgeous woodcuts of a monochordist playing various intervals. Both the numbers and the intervals (ratios) are labeled with similar ribbons---though the number ribbons point to places on the monochord a movable bridge could be located, and so indicate measured lengths, whereas the intervals connect the number ribbons. On the other hand, page 59 shows a series of compounding musical intervals (that is, compounding ratios). Most of the banners connect numbers---and so show the ratio and interval---but there is a similar banner connecting two ratios that labels their functional place in the operations the diagrams are illustrating.} So it was not a universal diagrammatic language for expressing ratios, even in music. On the other hand, in the scholia to his translation of Psellus' mathematical works, Xylander \citep{psellos1556liber} very clearly uses arcs to illustrate ratios (and differences) so the the ten species of mean can be clearly distinguished.\footnote{He uses curvy arcs, like here, to indicate ratios, pointy arcs to indicate differences.} Thus at least some of the time, as in the diagrams shown here, the notation seems to have been used to illustrate ratios \textit{specifically}.\footnote{Some of this practice may even be similar to our own use of arrows, which sometimes are used for simple labels, and sometimes to show functions.} This suggests that even the categorical arrow notation I apply to Euclid in this paper may have important historical precedents---though I have never seen it used to illustrate a Euclidean diagram.

Mathematically, through these four sample Propositions, this intuitive categorical approach to Euclid has been shown capable of rigorously illuminating Euclid's arguments and diagrams, in a manner that is faithful to Euclid's own philosophical world, and does not take recourse in foreign concepts like fractions or symbolic algebra.

Needham \citep{needham2021} suggests that one of the reason Newton's method isn't followed is that he didn't have a symbol for equality. This is true enough, but if Newton's approach to ratios coincides with Barrow's \citep{barrow-ML} it seems Newton combined an intuitive approach to geometry with a relational calculus that could have inspired Lawvere---though oddly, an essentialist one. It is this relational calculus that is obscured by our colon notation for ratios, and that needs a rigorous and illuminating symbolism that is rooted in Newton's intuitive diagrammatic reasoning, and rigorously grounded in Euclidean style axioms. I propose that, in contexts where the mathematicians attempted to write in an Euclidean manner---e.g., in Newton---we adopt this or a similar relational calculus, added immediately to the diagrams. This practice would reconcile very different streams of modern mathematics: Needham's intuitive approach, Lawvere's relational calculus, and type theories---and it would open the door for demonstrative, axiomatic proofs that, precisely therefore, answer the question \textit{why}---this is, after all, the very goal of Aristotelian axiomatics. Though Needham (whose work I highly admire) does not generally follow the Newtonian and Euclidean \textit{type system,} and rewriting \textit{Differential Geometry} with a Euclidean type-system would be a monumental labor.

Second, for a pedagogical perspective, this approach to Euclid (and to geometry) perhaps opens the door for a more functions-first approach to algebra II. Specifically, students whose geometry curriculum included categorical diagrams like those present in this paper may be prepared, in algebra, to read diagrams like those briefly presented from Heron and Diophantus, and so to consider functions compositionally.

\appendix
\section{Allegories}\label{AppendixAllegories}
The critical differences between Euclidean ratios, and allegories in the sense of Freyd and Scedrov \citep{freyd-scedrov1990} are:

\begin{enumerate}
\item Euclid's theory is two-sorted, whereas Freyd and Schedrov give an arrows only definition of a category (and hence, of an allegory). Though see point \ref{equality} below.

\item Though Euclid could consider $R \cap S,$ he generally does not. Moreover when $R$ and $S$ are incompatible, Euclid would not say $R\cap S = \bot$ but just that there is no such relation. For example, ``double of AND triple of'' is nonsense, not the relation $\bot$.

\item In Freyd and Scedrov, relations are not types, and are all of the same type. Greek mathematics reasons about universal \textit{types}, not sets (though in the late middle ages the nominalists interpreted types as sets \citep{johnofstthomas1883logica}) or individuals in the world (similar to homotopy type theory \citep{hottbook}, p. 18, the type is an integral part of the thing---its very nature---though it is the very nature of the object \textit{in the world}). When the types are considered, they form a Julia-style type-hierarchy \citep{juliatypes2026}, in which concrete types cannot have any sub-types, and abstract types are only inhabited when instantiated through a concrete leaf-type. The posetal ordering on the relations required for their forming an allegory is type inclusion, written $<:$ in Julia. Thus, for example, ratio of three-to-two is a concrete type of ratio that can be in and to different pairs of quantities. On the other hand, 3-to-2 $<:$ greater-than, so greater-than is an abstract type of ratio. This means that greater-than inhabits real-world ratios in and through their intrinsic typal nature, which is 3-to-2, or 4-to-3, or any other concrete type.

\item Though abstract types of ratio (in theory) can be composed, in general doing so would create an ugly infinite union type. For example, using the Pythagorean nomenclature, the immediate parent of 3-to-2 in the type hierarchy is called superparticular. If we attempt to compose superparticular with superparticular, we would get an ugly infinite union type: dupla-sesquiquartanORdoubleORsesquialteraOR... (again, using the Pythagorean terminology). This is extremely unnatural and is avoided entirely. Though he is not then explicitly treating it as a ratio, Euclid very occasionally (e.g., I.18) composes greater-than with greater-than to form greater-than. And the definition of \textit{same ratio} seen above, involves a sort of trivial composition (trivial, in the sense that they compose to the highest genus, namely, ratio). Concrete leaf-type ratios, however, compose unambiguously to form concrete ratios (see Appendix \ref{CompositionAppendix}), and their composition is an important feature of Euclidean mathematics, and of Greek mathematics more generally (and of music theory).

\item Because ratios are types that are independent of the specific quantity or type of quantity they are in and to, a number of quantities can have the same relation. For example, a triangle may be the double of two different squares. Most notably, this is also true of the identity relation, equality: One triangle may be equal to two different squares, and those squares are then equal to each other.

\item\label{equality} Notably, this means even the identity relation, equality, is, in any Proposition, from and to \textit{different} objects. But since they are types, ratio composition is independent of objects, in the sense that the composition of two leaf-type ratios does not depend on objects that happen to have the ratios. (This was assumed by Greek mathematicians \citep{acerbi2011}, but can be shown in a straight-forward manner from Propositions 5.11 and V.22. Again, Appendix \ref{CompositionAppendix}.) Seen from this perspective, V.7 can be read as establishing, using a notation and exposition very unnatural to Euclid, that, for any leaf-type ratio \(r\), and \(e\) the ratio of equality, \(e \circ r = r\) and \(r \circ e = r\). 

\centering%
\begin{tcolorbox}[%
	colbacktitle=white,
	coltitle = black,
	center title,
	halign upper=center,
	halign lower=center,
	fonttitle = \bfseries,
	colbacklower = white,
	colback = white,
	colframe=white,
	width = 0.95\linewidth,
	lower separated=false,
	sidebyside,
	boxrule=.2mm]
\begin{tikzcd}[ampersand replacement = \&]
A
\ar[dd,
	leftarrow,
	shift right = 0.5,
	bend right = 0.25cm,
	"e"']
\ar[dd,
	shift left = 0.5,
	bend left = 0.25cm,
	"e"]
\ar[dr,
	"r"]\\
\&	C\\
B
\ar[ur,
	"r"']
\end{tikzcd}%

Equal magnitudes ($A$, $B$) have the same ratio ($r$) to the same ($C$).
\tcblower
%Second Column
\begin{tikzcd}[ampersand replacement = \&]
A
\ar[dd,
	leftarrow,
	shift right = 0.5,
	bend right = 0.25cm,
	"e"']
\ar[dd,
	shift left = 0.5,
	bend left = 0.25cm,
	"e"]
\ar[dr,
	leftarrow,
	"r"]\\
\&	C\\
B
\ar[ur,
	leftarrow,
	"r"']
\end{tikzcd}\\
The same ($C$) has the same ratio~($r$) to equal magnitudes ($A$, $B$).
\end{tcolorbox}

\flushleftright%
Furthermore, V.9 can be read as establishing that if \(s\circ r = r\) then \(s = e\) and if \(r \circ s = r\) then \(s = e\).

\centering%
\begin{tcolorbox}[%
	colbacktitle=white,
	coltitle = black,
	center title,
	halign upper=center,
	halign lower=center,
	fonttitle = \bfseries,
	colbacklower = white,
	colback = white,
	colframe=white,
	width = 0.95\linewidth,
	lower separated=false,
	sidebyside,
	boxrule=.2mm]
\begin{tikzcd}[ampersand replacement = \&]
A
\ar[dd,
	leftarrow,
	shift right = 0.5,
	bend right = 0.25cm,
	"s^{\circ}"']
\ar[dd,
	shift left = 0.5,
	bend left = 0.25cm,
	"s"]
\ar[dr,
	"r"]\\
\&	C\\
B
\ar[ur,
	"r"']
\end{tikzcd}%

If two magnitudes ($A$, $B$) have the same ratio ($r$) to the same ($C$), they are equal ($s=e=s^{\circ}$).
\tcblower
%Second Column
\begin{tikzcd}[ampersand replacement = \&]
A
\ar[dd,
	leftarrow,
	shift right = 0.5,
	bend right = 0.25cm,
	"s^{\circ}"']
\ar[dd,
	shift left = 0.5,
	bend left = 0.25cm,
	"s"]
\ar[dr,
	leftarrow,
	"r"]\\
\&	C\\
B
\ar[ur,
	leftarrow,
	"r"']
\end{tikzcd}\\
If the same ($C$) has the same ratio~($r$) to two magnitudes ($A$, $B$), they are equal ($s = e=s^{\circ}$).
\end{tcolorbox}

\flushleftright%
This shows that within the system of leaf-type relations, equality has the role of the unique identity relation, even though, from the perspective of objects, it is always from one and to another. Though it would no longer be true to Euclid and wouldn't help us explain his arguments---since he always grounds ratios in quantities---we could therefore say that if we treat the empty relation as a relation, Euclid's ratios form an arrows only allegory with one identity (one object)---though see the caveat about the law of modularity just below.
\end{enumerate}

Since $\cap$ is seldom used, Freyd and Scedrov's law of modularity is unimportant (though true if restricted to relations $R$, $S$ where $R \cap S$ is well-defined). I nevertheless call this system like an allegory because it is relational, and makes relatively extensive use of composition, which is understood in the manner Freyd and Scedrov explain in their introduction to and motivation of allegories. (Though there, as motivation, they adopt a notation that looks more set-theoretic than the notation they finally adopt, and Euclid's relations should \textit{not} be understood in a set-theoretic manner).

Also, as a warning: Though we can say, e.g. ``double $<:$ multiple $<:$ greater $<:$ inequality $<:$ ratio,'' moving up from most proper species to highest genus of ratio, Euclid defines a \textit{different} partial ordering on relations such that double $<$ triple, even though both are incomparable leaf-nodes in the genus-species hierarchy. These two partial orders should be kept distinct, which is easy enough when working through Euclid, but may be hard when thinking of his relations as forming an allegory. It is the genus-species relation that allows his system to be like an allegory.

\section{Ratio Composition}\label{CompositionAppendix}

Greek mathematicians found \( \left( A\xrightarrow{to}B\right)\circ \left(C\xrightarrow{to}D\right) \) by finding $E$, $F$, and $G$ such that $A$ is to $B$ as $E$ is to $F$, and $C$ is to $D$ as $F$ is to $G$, as I did in the body of the paper. But they did not show that this is well defined---that is, that the composition is independent of the choice of $E$, $F$, and $G$.

Specifically, we take $E$, $F$ and $G$ as above so that $E$ is to $F$ as $A$ is to $B$ and $F$ is to $G$ as $C$ is to $D$. We then define (to speak in contemporary language) \( \left( A\xrightarrow{to}B\right)\circ \left(C\xrightarrow{to}D\right) \) as \( \left( E\xrightarrow{to}F\right)\circ \left(F\xrightarrow{to}G\right) \), that is, as \(E\xrightarrow{to}G\). If we select different $H$, $K$, and $L$ similarly, then $E$ is to $G$ as $H$ is to $L$, and so the composition is independent of the quantities that happen to have the ratio.

This can be argued in detail: select $E$, $F$, $G$, $H$, $K$, and $L$ so that the following diagram holds, 

\centering%
\begin{tikzcd}[column sep = 0cm, ampersand replacement = \&, row sep = 1.25cm]
E
\ar[rr,
	""' {name = EF}]
\ar[d,
	Leftrightarrow,
	gray,
	densely dotted]
\&[1cm]\&
	F
	\ar[rr,
		""' {name = FG}]
	\ar[dl,
		Leftrightarrow,
		gray,
		densely dotted]
	\ar[dr,
		Leftrightarrow,
		gray,
		densely dotted]
	\&\&[1cm]
		G
		\ar[d,
			Leftrightarrow,
			gray,
			densely dotted]\\
A
\ar[r,
	"" {name = AB1},
	""' {name = AB2}]
\ar[d,
	Leftrightarrow,
	gray,
	densely dotted]
\&
	B
	\ar[dr,
		Leftrightarrow,
		gray,
		densely dotted]
	\&
	,
	\&
		C
		\ar[r,
			"" {name = CD1},
			""' {name = CD2}]
		\ar[dl,
			Leftrightarrow,
			gray,
			densely dotted]
		\&
			D
			\ar[d,
				Leftrightarrow,
				gray,
				densely dotted]\\
H
\ar[rr,
	"" {name = HK}]
\&\&
	K
	\ar[rr,
		"" {name = KL}]
	\&\&
		L.
\ar[from = EF,
	to = AB1,
	Leftrightarrow,
	gray]
\ar[from = FG,
	to = CD1,
	Leftrightarrow,
	gray]
\ar[from = AB2,
	to = HK,
	Leftrightarrow,
	gray]
\ar[from = CD2,
	to = KL,
	Leftrightarrow,
	gray]
\end{tikzcd},

\flushleftright%
then by V.11, transitivity of ``same ratio'', $E$ is to $F$ as $H$ is to $K$ and similarly $F$ is to $G$ as $K$ is to $L$:

\centering%
\begin{tikzcd}[ampersand replacement = \&, row sep = 1.25cm]
E
\ar[r,
	""' {name = EF}]
\ar[d,
	Leftrightarrow,
	gray,
	densely dotted]
\&
	F
	\ar[r,
		""' {name = FG}]
	\ar[d,
		Leftrightarrow,
		gray,
		densely dotted]
	\&
		G
		\ar[d,
			Leftrightarrow,
			gray,
			densely dotted]\\
H
\ar[r,
	"" {name = HK}]
\&
	K
	\ar[r,
		"" {name = KL}]
	\&
		L.
\ar[from = EF,
	to = HK,
	Leftrightarrow,
	gray]
\ar[from = FG,
	to = KL,
	Leftrightarrow,
	gray]	
\end{tikzcd}

\flushleftright%
Hence, by V.22, $E$ is to $G$ as $H$ is to $L$, 

\centering%
\begin{tikzcd}[ampersand replacement = \&, row sep = 1.75cm]
E
\ar[r,
	""' {name = EF}]
\ar[rr,
	bend right,
	""' {name = EG}]
\ar[d,
	Leftrightarrow,
	gray,
	densely dotted]
\&
	F
	\ar[r,
		""' {name = FG}]
	\&
		G
		\ar[d,
			Leftrightarrow,
			gray,
			densely dotted]\\
H
\ar[r,
	"" {name = HK}]
\ar[rr,
	bend left,
	"" {name = HL}]
\&
	K
	\ar[r,
		"" {name = KL}]
	\&
		L.
\ar[from = EG,
	to = HL,
	Leftrightarrow,
	gray]	
\end{tikzcd}

\flushleftright%
And so \(\left( A\to B \right) \circ \left(C\to D \right) \) is defined unambiguously as either \( E\to G \) or \(H\to L\). Since, at the level of leaf-types, Euclid's ``category'' is codiscrete within each genus of quantity, and we are, here, considering leaf-type ratios, there is no need to check that the composition is well-defined when there is a shared-middle term.

\bibliography{EuclidPaper}

\end{document}